\documentclass[twoside,11pt,reqno]{amsart} 
\usepackage{amsmath,amssymb,amscd,mathrsfs,stmaryrd,wasysym,latexsym}

\makeatletter

\newtheorem{Theorem}{Theorem}[section]
\newtheorem{Corollary}[Theorem]{Corollary}
\newtheorem{Proposition}[Theorem]{Proposition}
\newtheorem{Lemma}[Theorem]{Lemma}

\theoremstyle{definition}
\newtheorem{Definition}[Theorem]{Definition}
\newtheorem{Remark}[Theorem]{Remark}
\newtheorem{Problem}[Theorem]{Problem}
\newtheorem{Example}[Theorem]{Example}
\newtheorem{Conjecture}[Theorem]{Conjecture}

\numberwithin{equation}{section}

\newbox\squ  
\setbox\squ=\hbox{\vrule width.6pt
   \vbox{\hrule height.6pt width.4em\kern1ex
         \hrule height.6pt}%
   \vrule width.6pt}

\def\lexeq{\leq_{\operatorname{lex}}}
\def\lex{<_{\operatorname{lex}}}

\newcommand{\V}{\mathbb V}

\def\op{\operatorname{op}}

\def\ad{\operatorname{ad}}

\def\e{e}

\def\Rep#1{\operatorname{Rep}(#1)}

\def\rep#1{\underline{\operatorname{Re}}\!\operatorname{p}(#1)}
\def\proj#1{\underline{\operatorname{Pro}}\!\operatorname{j}(#1)}
\def\Proj#1{\operatorname{Proj}(#1)}
\def\Mod#1{#1\!\operatorname{-Mod}}
\def\id{\operatorname{id}}
\def\col{\operatorname{col}}
\def\row{\operatorname{row}}

\def\db{\delta}
\def\dd{\eps}

\def\defect{\operatorname{def}}
\def\C{{\mathbb C}}
\def\Q{{\mathbb Q}}
\def\Z{{\mathbb Z}}

\def\onto{{\twoheadrightarrow}}
\def\into{{\hookrightarrow}}

\def\0{{\bar 0}}
\def\1{{\bar 1}}
\def\pr{{\operatorname{pr}}}
\def\infl{{\operatorname{infl}}}
\def\aff{{\operatorname{aff}}}

\def\St{{\mathscr{T}}}
\def\T{{\mathtt T}}
\def\Laurent{\mathscr A}
\def\Stab{{\mathtt S}}
\def\U{{\mathtt U}}
\def\Par{{\mathscr P}^\kappa}
\def\RPar{{\mathscr{RP}}^\kappa}

\newcommand{\CL}{{\mathcal C}}
\newcommand{\SCL}{{\mathcal S}}
\def\qdim{{\operatorname{qdim}}}
\def\hom{\underline{\operatorname{Hom}}}
\def\HOM{{\operatorname{HOM}}}
\def\Hom{{\operatorname{Hom}}}

\def\End{{\operatorname{End}}}
\def\END{{\operatorname{END}}}

\def\Ind{{\operatorname{Ind}}}

\def\res{{\operatorname{res}}}
\def\Res{{\operatorname{Res}}}

\def\cha{{\operatorname{char}\,}}

\def\soc{{\operatorname{soc}\:}}

\def\CH{{\operatorname{ch}_q\:}}
\def\wt{{\operatorname{wt}}}
\def\height{{\operatorname{ht}}}

\def\cont{{\operatorname{cont}}}

\def\bi{\text{\boldmath$i$}}
\def\bj{\text{\boldmath$j$}}

\def\sign{\operatorname{sgn}}

\def\eps{{\varepsilon}}
\def\phi{{\varphi}}
\def\emptyset{{\varnothing}}
\def\ga{{\gamma}}
\def\Ga{{\Gamma}}
\def\la{{\lambda}}
\def\La{{\Lambda}}
\def\de{{\delta}}
\def\De{{\Delta}}
\def\al{{\alpha}}
\def\be{{\beta}}
\def\si{{\sigma}}
\def\g{{\mathfrak g}}
\def\h{{\mathfrak h}}
\def\sl{{\mathfrak sl}}

\def\nslash{\:\notslash\:}

\def\underbar{\mathpalette\@underbar}
\def\@underbar#1#2{\settowidth{\@tempdimb}{$#1#2$}\@tempdimb=0.8\@tempdimb
                   \ooalign{$#1#2$\crcr%
                         \hfil\rule[-.5mm]{\@tempdimb}{.4pt}\hfil}}

\newdimen\hoogte    \hoogte=14pt    
\newdimen\breedte   \breedte=14pt   
\newdimen\dikte     \dikte=0.5pt    

\newenvironment{young}{\begingroup
       \def\vr{\vrule height0.8\hoogte width\dikte depth 0.2\hoogte}
       \def\fbox##1{\vbox{\offinterlineskip
                    \hrule height\dikte
                    \hbox to \breedte{\vr\hfill##1\hfill\vr}
                    \hrule height\dikte}}
       \vbox\bgroup \offinterlineskip \tabskip=-\dikte \lineskip=-\dikte
            \halign\bgroup &\fbox{##\unskip}\unskip  \crcr }
       {\egroup\egroup\endgroup}
\def\diagram#1{\relax\ifmmode\vcenter{\,\begin{young}#1\end{young}\,}\else%
              $\vcenter{\,\begin{young}#1\end{young}\,}$\fi}
              
              \newdimen\Hoogte    \Hoogte=12pt    
\newdimen\Breedte   \Breedte=12pt   
\newdimen\Dikte     \Dikte=0.5pt    

\newenvironment{Young}{\begingroup
       \def\vr{\vrule height0.8\Hoogte width\Dikte depth 0.2\Hoogte}
       \def\fbox##1{\vbox{\offinterlineskip
                    \hrule height\Dikte
                    \hbox to \Breedte{\vr\hfill##1\hfill\vr}
                    \hrule height\Dikte}}
       \vbox\bgroup \offinterlineskip \tabskip=-\Dikte \lineskip=-\Dikte
            \halign\bgroup &\fbox{##\unskip}\unskip  \crcr }
       {\egroup\egroup\endgroup}
\def\Diagram#1{\relax\ifmmode\vcenter{\,\begin{Young}#1\end{Young}\,}\else%
              $\vcenter{\,\begin{Young}#1\end{Young}\,}$\fi}

\begin{document}

\title[Symmetric Groups and Hecke Algebras]{Representation Theory of Symmetric Groups and Related Hecke Algebras}
\author{Alexander Kleshchev}

\begin{abstract}
We survey some fundamental trends in representation theory of symmetric groups and related objects which became apparent in the last fifteen years. The emphasis is on connections with Lie theory via categorification. We present results on branching rules and crystal graphs, decomposition numbers and canonical bases, graded representation theory, connections with cyclotomic and affine Hecke algebras, Khovanov-Lauda-Rouquier algebras, category ${\mathcal O}$, $W$-algebras,  \dots
\end{abstract}
\thanks{{\em 2000 Mathematics Subject Classification:} 20C30, 20C08, 17B37, 20C20, 17B67.}
\thanks{Supported in part by the NSF grant 
DMS-0654147. The paper was completed while the author was visiting the Isaac Newton Institute for Mathematical Sciences in Cambridge, U.K. The author thanks the Institute for hospitality and support.}
\address{Department of Mathematics, University of Oregon, Eugene, USA.}
\email{klesh@uoregon.edu}
\maketitle


\section{Introduction}\label{SIntro}

The {\em symmetric group} $\Sigma_d$ on $d$ letters is a classical, fundamental, deep, well-studied, and much-loved mathematical object. It has natural connections with combinatorics, group theory, Lie theory, geometry, topology, \dots, natural sciences.  
Since the work of Frobenius in the end of the nineteenth century, {\em representation theory}\, of symmetric groups has developed into a large and important area of mathematics. 

In this expository article we survey {\em some}\, fundamental trends in representation theory of symmetric groups and related objects which became apparent in the last fifteen years. The emphasis is on connections with Lie theory via categorification. We present results on branching rules and crystal graphs, decomposition numbers and canonical bases, graded theory, connections with cyclotomic and affine Hecke algebras, Khovanov-Lauda-Rouquier  algebras, category ${\mathcal O}$, $W$-algebras, etc.

The main problem of representation theory is to understand irreducible modules. Let us mention up front that over fields of positive characteristic the dimensions of irreducible $\Sigma_d$-modules are not known. Trying to approach this problem of {\em modular representation theory} by induction lead us to studying restrictions of irreducible modules 
from $\Sigma_d$ to $\Sigma_{d-1}$, which resulted in a discovery of {\em modular branching rules} \cite{KBrI}--\cite{KBrIV}. 

These branching rules turned out to provide  only  partial information on dimensions of irreducible modules. However the underlying subtle combinatorics lead to a discovery by Lascoux, Leclerc and Thibon \cite{LLT} of some surprising and deep connections  between representation theory of $\Sigma_d$ and representations of {\em quantum Kac-Moody algebras}. Roughly speaking, it turned out that the modular branching rule corresponds to the {\em crystal graph} in the sense of Kashiwara of the basic module over a certain Kac-Moody algebra $\g$. This observation turned out to be a beginning of an exciting development which continues to this day. 

Since the work of Dipper and James \cite{DJ} it has been known that representation theory of $\Sigma_d$ over a field $F$ of characteristic $p$ resembles representation theory of the corresponding {\em Iwahori-Hecke algebra} over the complex field $\C$ at a $p$th root of unity. 

In fact it makes sense to work more generally from the beginning with the algebra $H_d=H_d(F,\xi)$ over an arbitrary field $F$ with a parameter $\xi\in F^\times$, which is given by generators $$T_1,\dots,T_{d-1}$$  and relations 
\begin{align}
\label{QECoxeterQuad}
T_r^2 &= (\xi -1) T_r +\xi\qquad (1\leq r<d),\\
T_rT_{r+1}T_r&=T_{r+1}T_rT_{r+1}\qquad(1\leq r<d-1),\label{QCoxeterClose}\\
T_rT_s&=T_sT_r \qquad(1\leq r,s<d,\ |r-s|>1).\label{QCoxeterFar}
\end{align}
Denote by $e$ the smallest positive integer such that 
$$1+\xi+\dots+\xi^{e-1} = 0,$$
setting $e := 0$
if no such integer exists. We refer to $e$ as the {\em quantum characteristic}. 
The Kac-Moody algebra $\g$ which we have alluded to above is 
$$\g=
\left\{
\begin{array}{ll}
\widehat{\sl}_e(\C) &\hbox{if $e>0$;}\\
\sl_\infty(\C) &\hbox{if $e=0$.}
\end{array}
\right.
$$

Note that when $\xi=1$, we have $H_d=F\Sigma_d$. 
If the field $F$ has characteristic $p>0$, then the relationship between the symmetric group algebra $H_d(F,1)=F\Sigma_d$ and the Hecke algebra $H_d(\C,e^{2\pi i/p})$   can be made precise using a reduction modulo $p$ procedure, which to every irreducible $H_d(\C,e^{2\pi i/p})$-module associates an $F\Sigma_d$-module. 

Even though reductions modulo $p$ of irreducible modules over the Hecke algebra are not always irreducible, {\em James' Conjecture} \cite{JamConj} predicts that they are in the {\em James region}. At any rate, irreducible modules over the Hecke algebra $H_d(\C,e^{2\pi i/p})$ can be considered as good `approximations' of irreducible modules over the symmetric group algebra $F\Sigma_d$. It turns out that the story originating in the Lascoux-Leclerc-Thibon paper leads to a rather satisfactory understanding of at least the irreducible modules over $H_d(\C,e^{2\pi i/p})$. 

To be more precise, Lascoux, Leclerc and Thibon conjectured a very precise connection between the {\em canonical bases} of modules over affine Kac-Moody algebras $\g$ in the sense of Lusztig \cite{Lubook} and Kashiwara \cite{Ka1}--\cite{Kas} on the one hand, and {\em projective indecomposable modules} over the Iwahori-Hecke algebras $H_d(\C,e^{2\pi i/p})$ on the other. Lascoux, Leclerc and Thibon also conjectured an explicit combinatorial algorithm for computing {\em decomposition numbers}, that is the multiplicities of the irreducible $H_d(\C,e^{2\pi i/p})$-modules in the corresponding {\em Specht modules}. 

It is easy to see that knowing decomposition numbers is sufficient for computing the dimensions and even characters of irreducible modules. The Lascoux-Leclerc-Thibon algorithm yields certain polynomials with non-negative coefficients which, when evaluated at $1$, conjecturally compute the decomposition numbers for $H_d(\C,e^{2\pi i/p})$. 

Building on powerful geometric results of Kazhdan-Lusztig \cite{KL}
and Ginzburg \cite[Chapter 8]{CG}, Ariki \cite{Ariki} has proved the  conjecture of Lascoux, Leclerc and Thibon, thus giving us a good understanding of modules over the complex Iwahori-Hecke algebras at roots of unity. A proof was also announced, but not  published, by Grojnowski. Later on Varagnolo and Vasserot \cite{VV1} proved a similar theorem for Schur algebras. 

A nice feature of Ariki's work is that he gets his results for a class of algebras more general than Iwahori-Hecke algebras. These algebras, known as {\em cyclotomic Hecke algebras} or {\em Ariki-Koike algebras}, were discovered independently by Cherednik \cite{Ch}, Brou\'e-Malle \cite{BM}, and Ariki-Koike \cite{AK}. They can be thought of as the Hecke algebras of complex reflection groups of types $G(\ell,1,d)$ in the Shephard-Todd classification. 

The cyclotomic Hecke algebras are denoted $$H_d^\La=H_d^\La(F,\xi),$$ where $d$ is a non-negative integer,  $F$ is the ground field, $\xi\in F^\times$ is a parameter, and $\La$ is a dominant integral weight for the Kac-Moody algebra $\g$ defined above. For the case where $\La$ is the fundamental dominant weight $\La_0$, we have $H_d^{\La_0}(F,\xi)=H_d(F,\xi)$, and so we incorporate the Iwahori-Hecke algebras and the group algebras of the symmetric groups naturally into the more general class of cyclotomic Hecke algebras.

A useful way to think of the connection between cyclotomic Hecke algebras and Lie theory is in terms of the idea of {\em categorification}, which goes back to I. Frenkel. It turns out that the  finite dimensional modules over the cyclotomic Hecke algebras $H_d^\La(F,\xi)$ for all $d\geq 0$ 
categorify the irreducible highest weight module $V(\La)$ over $\g$. 

This statement can be made much more precise, connecting various important stories in representation theory of cyclotomic Hecke algebras to important invariants of the module $V(\La)$. For example:
\begin{enumerate}
\item[(1)] the action of the Chevalley generators of $\g$ corresponds to the functors of $i$-induction and $i$-restriction on the categories of modules over the cyclomotic Hecke algebras \cite{Abook,LLT,Gr}. 
\item[(2)] the weight spaces of $V(\La)$ correspond to the blocks of the cyclotomic Hecke algebras \cite{Abook,Gr}; 
\item[(3)] the crystal graph of $V(\La)$ corresponds to the socle branching rule for the cyclotomic Hecke algebras \cite{LLT,KBrII,KBrIII,MiM,Gr,Abranch}; 
\item[(4)] the Shapovalov form on $V(\La)$ corresponds to the Cartan pairing on the Grothendieck group of modules over the cyclotomic Hecke algebras \cite{Gr}; 
\item[(5)] the action of the Weyl group of $\g$ on $V(\La)$ corresponds to certain derived equivalences between blocks conjectured by Rickard and constructed by Chuang and Rouquier \cite{CR};
\item[(6)] elements of a standard spanning set of $V(\La)$ coming from the construction of $V(\La)$ in terms of a higher level Fock space correspond to the classes of the Specht modules over cyclotomic Hecke algebras \cite{Abook}; 
\item[(7)] provided $F=\C$, elements of the dual canonical basis in $V(\La)$ correspond to the classes of the irreducible $H_d^\La(\C,\xi)$-modules \cite{Ariki,Abook}; 
\item[(8)] provided $F=\C$, elements of the canonical basis in $V(\La)$ correspond to the classes of projective indecomposable $H_d^\La(\C,\xi)$-modules \cite{Ariki,Abook}. 

\end{enumerate}

One thing which remains unexplained in the picture described above is the role of the {\em quantum group}. The categorification by Ariki and Grojnowski is only a categorification of $V(\La)$ as a module over $\g$, not over the {\em quantized enveloping algebra}\, $U_q(\g)$. On the other hand, appearance of canonical bases evaluated at $q=1$ suggests that the quantized enveloping algebra $U_q(\g)$, which so far remained invisible, should be relevant. So the picture seems to be incomplete unless one actually categorifies a $q$-analogue of $V(\La)$. A standard way of doing this is to find an appropriate grading on the cyclotomic Hecke algebras and then consider {\em graded representation theory}, with the action of the parameter $q$ on the Grothendieck group corresponding to the `grading shift' on modules. 

The existence of important well-hidden gradings on the blocks of cyclotomic Hecke algebras, and in particular group algebras of symmetric groups, has been predicted by Rouquier \cite{R1} and Turner \cite{T}. Recently Brundan and the author \cite{BKyoung} were able to construct such gradings. More precisely, we construct an explicit isomorphism between the cyclotomic Hecke algebras and certain  cyclotomic Khovanov-Lauda-Rouquier algebras defined independently by Khovanov-Lauda \cite{KL1,KL2} and Rouquier \cite{Ro}. The Khovanov-Lauda-Rouquier algebras are naturally $\Z$-graded, so combining with our isomorphism, we obtain an explicit grading on the blocks of cyclotomic Hecke algebras. 

In \cite{BKW} we then grade Specht modules, which allows us to define  {\em graded decomposition numbers}. Finally, in \cite{BKllt}, we prove that, for the cyclotomic Hecke algebras over $\C$, these graded decomposition numbers are precisely the {\em polynomials} coming from the conjecture of Lascoux, Leclerc, and Thibon (generalized to the cyclotomic case). 

We also use graded representation theory of the cyclotomic Hecke algebras to categorify $V(\La)$ as a module over $U_q(\g)$, and to obtain graded analogues of the results (1)--(8) described above. This categorification result (except for (6)) has been also announced by Rouquier in a more general setting of Khovanov-Lauda-Rouquier algebras of general type. Related results on canonical bases of $U^-_q(\g)$ have been obtained by Varagnolo and Vasserot \cite{VV3}.

The gradings also allowed Brundan and the author to define the {\em $q$-characters} of modules over cyclotomic Hecke algebras and to determine the $q$-characters of Specht modules. As a consequence, we obtain a graded dimension formula for the blocks of cyclotomic Hecke algebras. 
Finding $q$-characters of irreducible modules of symmetric groups can be considered the main problem of its representation theory. This problem is equivalent to finding the corresponding graded decomposition numbers.

The paper is organized as follows: 

\noindent
\ref{SIntro}. Introduction

\noindent
\ref{SMO}. Main Objects

\ref{SSGrFPar}. Ground field and parameters

\ref{SSGRT}. Graded representation theory

\ref{SSSGIHA}. Symmetric groups and Iwahori-Hecke algebras

\ref{SHD}. Homogeneous generators

\ref{SSDynkin}. Weights and roots

\ref{SSHP}. Homogeneous presentation

\ref{SSCHA}. Affine Hecke algebras

\ref{SSCHANew}. Cyclotomic Hecke algebras

\ref{SSBlocks}. Blocks

\ref{SSqCh}. The main problem

\ref{SSAKLRA}. Affine Khovanov-Lauda-Rouquier algebras

\noindent
\ref{SComb}. Combinatorics

\ref{SSPar}. Partitions and Young diagrams

\ref{SSTab}. Tableaux

\ref{SSDeg}. Degree of a standard tableau

\ref{SSCrGraph}. Good nodes and restricted multipartitions

\noindent
\ref{SSMPCharZerLev1}. Solution of the Main Problem for type ${A_\infty}$ at level $1$

\noindent
\ref{SCell}. Cyclotomic Hecke algebras as cellular algebras

\ref{SSCellGen}. Review of cellular algebras

\ref{SSCelCyc}. Cellular structures on cyclotomic Hecke algebras

\ref{SSSMIM}. Specht modules and irreducible modules

\ref{SSBlAg}. Blocks again

\noindent
\ref{SGrMod}. Graded modules over cyclotomic Hecke algebras

\ref{SSGIrrMod}. Graded irreducible modules

\ref{SSHBSM}. Homogeneous bases of Specht modules

\ref{SSGBRSM}. Graded branching rule for Specht modules

\ref{SSGDB}. Graded dimension of a block

\ref{SSGrCellStr}. Graded cellular structure on cyclotomic Hecke algebras

\noindent
\ref{SGIndResBr}. Graded induction, restriction, and  branching rules

\ref{SSAIndRes}. Affine induction and restriction

\ref{SSAIIndIRes}. Affine $i$-induction and $i$-restriction

\ref{SSADP}. Affine divided powers

\ref{sir}. Cyclotomic $i$-induction and $i$-restriction

\ref{sdivp}. Cyclotomic divided powers

\ref{SSBrBrRules}. Graded branching rules for irreducible modules

\noindent
\ref{SQG}. Quantum groups

\ref{SSAF}. The algebra $\mathbf f$

\ref{SSQEAUG}. The quantized enveloping algebra $U_q(\g)$

\ref{md}. The module $V(\La)$

\ref{SSFock}. Fock spaces 

\ref{ssb}. Canonical bases

\noindent
\ref{SCat}. Categorifications

\ref{sklth}. Categorification of $\mathbf f$

\ref{sc}. Categorification of $V(\La)$

\ref{SSMBSM}. Monomial bases and Specht modules

\ref{SSCBGDN}. Canonical bases and graded decomposition numbers

\ref{SSAlg}. An algorithm for computing decomposition numbers for $H_d(\C,\xi)$

\noindent
\ref{SJC}. Reduction Modulo $p$ and James' Conjecture

\ref{SRed}. Realizability over prime subfields

\ref{SSRedNew}. Reduction modulo $p$

\ref{sgam}. Graded adjustment matrices

\ref{SSJamConj}. James Conjecture

\noindent
\ref{SSBr}. Some other results

\ref{SSCR}. Blocks of symmetric groups: Brou\'e's conjecture and Chuang-Rouquier equivalences

\ref{SSBrLab}. Branching and labeling of irreducible modules

\ref{sex}. Extremal sequences

\ref{SSMoreBr}. More on branching for symmetric groups

\ref{SSMull}. Mullineux Involution

\ref{SSHLSWD}. Higher level Schur-Weyl duality, $W$-algebras, and category $\mathcal O$

\ref{SSProj}. Projective representations

\ref{SSASP}. Problems on symmetric groups related to Aschbacher-Scott program

\section{Main Objects}\label{SMO}

\subsection{Ground field and parameters}\label{SSGrFPar}
Let $F$ be an algebraically closed field, and 
$\xi \in F^\times$ be an invertible element. Denote by $e$  the {\em quantum characteristic}, i.e. the smallest positive integer such that 
$$1+\xi+\dots+\xi^{e-1} = 0,$$
setting $e := 0$
if no such integer exists. 
For example, if $\xi=1$, then $e=\operatorname{char} F$. If $\xi\neq 1$ then $\xi$ is a primitive $e$th root of unity if $e>0$, and $\xi$ is {\em generic} if $\e=0$. 
Define $$I:=\Z/e\Z.$$ 
For $i\in I$, we have a well-defined element $\nu(i)$ of $F$ defined as follows:
\begin{equation}\label{ENu}
\nu(i):= 
\left\{
\begin{array}{ll}
i &\hbox{if $\xi=1$;}\\
\xi^i &\hbox{if $\xi\neq1$.}
\end{array}
\right.
\end{equation}
Throughtout the paper $q$ is an indeterminate, and $\Laurent:=\Z[q,q^{-1}]$. As usual, set 
\begin{equation}\label{EGI}
[n]:=\frac{q^n-q^{-n}}{q-q^{-1}},\quad [n]!:=[n][n-1]\dots[1],\quad 
\left[
\begin{matrix}
 n   \\
 m
\end{matrix}
\right]
:=\frac{[n]!}{[n-m]![m]!}.
\end{equation}

\subsection{Graded representation theory}\label{SSGRT}
Later on in this article we explain how to grade symmetric group algebras and more generally cyclotomic Hecke algebras and advocate the idea of studying their {\em graded representation theory}. Since all the `usual' irreducible modules over finite dimensional $\Z$-graded algebras are gradable, by studying graded irreducible modules we `do not lose any information' but actually gain an additional insight. 

To explain this precisely, let, more generally, $H$ be a $\Z$-graded $F$-algebra, and  $\Mod{H}$ denote the abelian category of all graded left $H$-modules, with 
morphisms being {\em degree-preserving} module homomorphisms, which we denote by $\Hom$.
Let $\Rep{H}$ denote
the abelian subcategory of all
finite dimensional graded $H$-modules and
 $\Proj{H}$ denote the additive subcategory  of 
all finitely generated projective graded $H$-modules. 

Denote the corresponding Grothendieck groups by
 $[\Rep{H}]$ and $[\Proj{H}]$, respectively.
We view these as $\Laurent$-modules via 
\begin{equation}\label{EAMod}
q^m[M]:=[M\langle m\rangle],
\end{equation}
where $M\langle m\rangle$ denotes the module obtained by 
shifting the grading up by $m$:
\begin{equation}\label{obvious}
M\langle m\rangle_n=M_{n-m}.
\end{equation}
Given $f = \sum_{n \in \Z} f_n q^n 
\in \Z_{\geq 0}[[q,q^{-1}]]$ and $M\in\Mod{H}$,
we write $$f \cdot M:=\bigoplus_{n \in \Z} M \langle n \rangle^{\oplus f_n}.$$

For $n \in \Z$, we let
$$
\Hom_H(M, N)_n := \Hom_H(M \langle n \rangle, N)
= \Hom_H(M, N \langle -n \rangle)
$$
denote the space of all homomorphisms
that are homogeneous of degree $n$,
i.e. they map $M_i$ into $N_{i+n}$ for each $i \in \Z$.
Set
$$
\HOM_H(M,N) := \bigoplus_{n \in \Z} \Hom_H(M,N)_n,
\quad
\END_H(M) := \HOM_H(M,M).
$$
There is a canonical  {\em Cartan pairing}
$$
\langle.,.\rangle:[\Proj{H}]\times[\Rep{H}] \rightarrow \Laurent,
\quad 
\langle[P],[M]\rangle := \qdim\ \HOM_H(P,M),
$$
where $\qdim \, V$ denotes $\sum_{n \in \Z} q^n \dim V_n$
for any finite dimensional 
graded vector space $V$.
Note that the Cartan pairing is {\em sesquilinear}, i.e. anti-linear in the first argument and  linear in the second.

We denote the category of  finite dimensional
ungraded $H$-modules (resp.\  finitely generated projective
ungraded $H$-modules) by $\rep{H}$ (resp.\ $\proj{H}$),
with Grothendieck group
 $[\rep{H}]$ (resp.\ $[\proj{H}]$).
We denote homomorphisms in these categories by $\hom$.
Given a graded module $M$, we write $\underline{M}$ for the
ungraded module obtained from it by forgetting the grading.
For $M, N \in \Rep{H}$, we have that
\begin{equation}\label{Ehom}
\hom_H(\underline{M}, \underline{N}) = \underline{\HOM_H(M,N)}.
\end{equation}
Informally speaking, the following standard lemmas show that in studying graded representation theory, we do not lose any information compared to the ungraded representation theory, but actually gain an additional insight. 

\begin{Lemma}[{\cite[Theorem 4.4.6, Remark 4.4.8]{NO}}]\label{no1}
If $M$ is any finitely generated graded $H$-module,
the radical of $\underline{M}$ is a graded submodule of $M$.
\end{Lemma}

\begin{Lemma}[{\cite[Theorem 4.4.4(v)]{NO}}]\label{no2}
If $L \in \Rep{H}$ is irreducible then
$\underline{L} \in \rep{H}$ is irreducible too.
\end{Lemma}

\begin{Lemma}[{\cite[Theorem 9.6.8]{NO}, \cite[Lemma 2.5.3]{BGS}}]\label{no3}
Assume that $H$ is finite dimensional.
If $K \in \rep{H}$ is irreducible,
then there exists an irreducible $L \in \Rep{H}$ 
such that
$\underline{L} \cong K$.
Moreover, $L$ is unique up to isomorphism
and grading shift.
\end{Lemma}

Given $M, L \in \Rep{H}$ with $L$ irreducible, 
we write $[M:L]_q$ for the {\em $q$-composition multiplicity},
i.e. 
$$[M:L]_q := \sum_{n \in \Z} a_n q^n,$$ 
where $a_n$ is the multiplicity
of $L\langle n\rangle$ in a graded composition series of $M$.
In view of Lemma~\ref{no2}, we recover the ordinary composition multiplicity
$[\underline{M}:\underline{L}]$ from
$[M:L]_q$ on setting $q$ to $1$.

\subsection{Symmetric groups and Iwahori-Hecke algebras}\label{SSSGIHA}
Always, $\Sigma_d$ is the symmetric group on $d$ letters with transpositions $(r,s)$ and  simple transpositions 
$$s_r:=(r,r+1)\qquad (1\leq r<d).$$ 
Denote by $F\Sigma_d$ the group algebra of $\Sigma_d$ over the ground field $F$.

The {\em Iwahori-Hecke algebra} of $\Sigma_d$ with the parameter $\xi$  is the $F$-algebra $H_d=H_d(F,\xi)$ given by generators $T_1,\dots,T_{d-1}$  and the relations (\ref{QECoxeterQuad})--(\ref{QCoxeterFar}). We will normally use the short version of the notation $H_d$ rather than $H_d(F,\xi)$, with the understanding that the ground field $F$ and the parameter $\xi$ are fixed. Only when $F$ and $\xi$ are not clear from the context or when we have more than one pair  $(F,\xi)$ in play, will we specify the field and the parameter explicitly.


If $\xi=1$ then $H_d$ is identified with $F\Sigma_d$ so that the generator $T_r$ corresponds to the simple transposition $s_r$ for each $1\leq r<d$. Thus the main object of our interest, the symmetric group algebra, is incorporated into the family of algebras $H_d$ depending on the fixed parameter $\xi\in F^\times$.

Define the {\em Jucys-Murphy elements} 
$L_1,\dots,L_d\in H_d$:\begin{equation}\label{EMurphy}
L_r:=
\left\{
\begin{array}{ll}
(1,r)+(2,r)+\dots+(r-1,r) &\hbox{if $\xi=1$;}\\
\xi^{1-r}T_{r-1}\dots T_2T_1T_1T_2\dots T_{r-1} &\hbox{if $\xi\neq 1$.}
\end{array}
\right.
\qquad(1\leq r\leq d).
\end{equation}
It is well-known and easy to check that the Jucys-Murphy elements commute, see e.g. \cite{Jucys2,Jucys3,MurphyYoung,MathasB}. The {\em Gelfand-Zetlin subalgebra} is the commutative subalgebra $\langle L_1,\dots,L_d\rangle\subset H_d$ generated by the Jucys-Murphy elements. 

Okounkov and Vershik \cite{OV} (cf. also \cite{DG}) have advocated the idea of studying  representation theory of $H_d$ by exploiting  the Gelfand-Zetlin subalgebra as a `Cartan subalgebra'. In particular, one should study the corresponding `weight spaces' in $H_d$-modules. The following comes from   \cite[Lemma 4.7]{Gr} and \cite[Lemma 7.1.2]{Kbook}.

\begin{Lemma}\label{LWInt} 
Let $M$ be a finite dimensional $H_d$-module. Then all eigenvalues of $L_1,\dots,L_d$ in $M$ are of the form $\nu(i)$ for $i\in I$. 
\end{Lemma}

Let $\bi=(i_1,\dots,i_d)\in I^d$, and $M$ be a finite dimensional $H_d$-module. Define the {\em $\bi$-weight  space} of $M$ as follows:
$$
M_\bi=\{v\in M\mid (L_r-\nu(i_r))^Nv=0\ \text{for $N\gg0$ and $r=1,\dots,d$}\}.
$$
By Lemma~\ref{LWInt}, we have a {\em weight space decomposition}:
$$
M=\bigoplus_{\bi\in I^d} M_\bi.
$$

\subsection{Homogeneous generators}\label{SHD}
Using the weight space decomposition of the left regular $H_d$-module, one gets a system of orthogonal idempotents
\begin{equation}\label{EIdempotents}
\{e(\bi)\mid \bi\in I^d\}
\end{equation}
in $H_d$, almost all of which are zero, such that
$$
\sum_{\bi\in I^d}e(\bi)=1,
$$
and 
$$
e(\bi)M=M_\bi\qquad(\bi\in I^d)
$$
for any finite dimensional $H_d$-module $M$, cf. \cite{MurphyId,BKyoung,Ma}.

Now define a  family of (nilpotent) elements $y_1,\dots,y_r\in H_d$ via:
\begin{equation}\label{EY}
y_r:=
\left\{\begin{array}{ll}
\sum_{\bi\in I^d}(1-\xi^{-i_r} L_r)e(\bi)&\text{if $\xi \neq 1$}\\
\sum_{\bi\in I^d}(L_r-i_r) e(\bi)&\text{if $\xi=1$}
\end{array}\right.
\qquad(1\leq r\leq d).
\end{equation}

In \cite{BKyoung}, for every $\bi\in I^d$ and $1\leq r<d$, we define explicitly power series    
$
P_r(\bi), Q_r(\bi)\in F[[y_r,y_{r+1}]]
$
such that $Q_r(\bi)$ has non-zero constant term. 
As $y_r$'s are nilpotent in $H_d$, we can interpret $P_r(\bi)$ and $Q_r(\bi)$ as  elements of $H_d$, with $Q_r(\bi)$ being invertible. The precise form of these elements is not going to be important---we just mention that there is some freedom in choosing $Q_r(\bi)$ and refer the interested reader to \cite[sections 3.3, 4.3]{BKyoung} for details. 
Set
\begin{equation}\label{QQCoxKL}
\psi_r:=
\sum_{\bi\in I^d}(T_r+P_r(\bi))Q_r(\bi)^{-1}e(\bi) \qquad(1\leq r<d). 
\end{equation}

The main result of \cite{BKyoung} claims that $H_d$ is generated by the elements 
\begin{equation}\label{EKLGens}
\{e(\bi)\:|\: \bi\in I^d\}\cup\{y_1,\dots,y_{d}\}\cup\{\psi_1, \dots,\psi_{d-1}\}
\end{equation}
and describes defining relations between these generators. This presentation turns out to  yield a hidden grading on $H_d$ which plays a fundamental role. Existence of such gradings was conjectured by Rouquier \cite[Remark 3.11]{R1} and Turner \cite{T}.

In order to describe the graded presentation of $H_d$, we need 
some rudimentary Lie theoretic notation introduced in the next subsection. 
At first, this Lie-theoretic terminology will play a purely notational or combinatorial role. However, it will gradually become clear that connections with Lie theory hinted at here are  deep and natural.

\subsection{Weights and roots}\label{SSDynkin}
Let $\Ga$ be the quiver with vertex set  $I$,
and a directed edge from $i$ to $j$ if $j  = i+1$.
Thus $\Gamma$ is the quiver of type $A_\infty$ if $e=0$
or $A_{e-1}^{(1)}$ if $e > 0$, with a specific orientation:
\begin{align*}
A_\infty&:\qquad\cdots \longrightarrow-2\longrightarrow -1 \longrightarrow 0 \longrightarrow 1 \longrightarrow 
2\longrightarrow \cdots\\
A_{e-1}^{(1)}&:\qquad0\rightleftarrows 1
\qquad
\begin{array}{l}
\\
\,\nearrow\:\:\:\searrow\\
\!\!2\,\longleftarrow\, 1
\end{array}
\qquad
\begin{array}{rcl}\\
0&\!\rightarrow\!&1\\
\uparrow&&\downarrow\\
3&\!\leftarrow\!&2
\end{array}
\qquad
\begin{array}{l}
\\
\:\nearrow\quad\searrow\\
\!4\qquad\quad \!1\\
\nnwarrow\quad\quad\,\sswarrow\\
\:\:3\leftarrow 2
\begin{picture}(0,0)
\put(-152.5,41){\makebox(0,0){0}}
\put(-13.5,53.5){\makebox(0,0){0}}
\end{picture}
\end{array}
\qquad \cdots
\end{align*}
The corresponding {\em Cartan matrix} 
$(a_{i,j})_{i, j \in I}$ is defined by
\begin{equation}\label{ECM}
a_{i,j} := \left\{
\begin{array}{rl}
2&\text{if $i=j$},\\
0&\text{if $i \nslash j$},\\
-1&\text{if $i \rightarrow j$ or $i \leftarrow j$},\\
-2&\text{if $i \rightleftarrows j$}.
\end{array}\right.
\end{equation}
Here the symbol 
$i \nslash j$
indicates that $j \neq i, i \pm 1$.

Following \cite{Kac}, let $(\h,\Pi,\Pi^\vee)$ be a realization of the Cartan matrix $(a_{ij})_{i,j\in I}$, so we have the simple roots 
$$\{\al_i\mid i\in I\},$$ 
the fundamental dominant weights 
$$\{\La_i\mid i\in I\},$$ and the normalized invariant form $(\cdot,\cdot)$ such that
$$
(\al_i,\al_j)=a_{ij}, \quad (\La_i,\al_j)=\de_{ij}\qquad(i,j\in I).
$$
Let $P_+$ be the set of dominant integral weights, and 
$$
Q_+ := \bigoplus_{i \in I} \Z_{\geq 0} \alpha_i
$$ 
denote the positive part of the root lattice. For $\alpha \in Q_+$, we write $\height(\alpha)$ for the {\em height of $\al$}, i.e. the sum of its 
coefficients when expanded in terms of the $\alpha_i$'s.

In this paper we will always work with a fixed positive integer $l$, referred to as the {\em level},  and an ordered $l$-tuple  
\begin{equation}\label{EKappa}
\kappa=(k_1,\dots,k_l)\in I^l.
\end{equation}
We will also need the corresponding dominant weight $\La$ (of level $l$) defined as 
\begin{equation}\label{ELa}
\La=\La(\kappa):=\La_{k_1}+\dots+\La_{k_l}\in P_+.
\end{equation}

\subsection{Homogeneous presentation}\label{SSHP}
Now we can state the main result of \cite{BKyoung}: 

\begin{Theorem} \label{TBK}
The algebra $H_d$ is generated by the elements 
(\ref{EKLGens}) subject 
only to the following relations for $\bi,\bj\in I^d$ and all 
admissible $r, s$:
\begin{equation}
e(\bi) e(\bj) = \de_{\bi,\bj} e(\bi);\label{R0}
\end{equation}
\begin{equation}
\sum_{\bi \in I^d} e(\bi) = 1;\label{R1}
\end{equation}
\begin{equation}
y_r e(\bi) = e(\bi) y_r;
\end{equation}
\begin{equation}
\psi_r e(\bi) = e(s_r{ }\bi) \psi_r;\label{R2PsiE}
\end{equation}
\begin{equation}
\label{R3Y}
y_r y_s = y_s y_r;
\end{equation}
\begin{equation}
\label{R3YPsi}
\psi_r y_s  = y_s \psi_r\qquad\text{if $s \neq r,r+1$};
\end{equation}
\begin{equation}
\psi_r y_{r+1} e(\bi) = 
\left\{
\begin{array}{ll}
(y_r\psi_r+1)e(\bi) &\hbox{if $i_r=i_{r+1}$},\\
y_r\psi_r e(\bi) \hspace{8mm}&\hbox{if $i_r\neq i_{r+1}$};
\end{array}
\right.
\label{R6}
\end{equation}
\begin{equation}
y_{r+1} \psi_re(\bi) =
\left\{
\begin{array}{ll}
(\psi_r y_r+1) e(\bi) &\hbox{if $i_r=i_{r+1}$},\\
\psi_r y_r e(\bi)  \quad&\hbox{if $i_r\neq i_{r+1}$};
\end{array}
\right.
\label{R5}
\end{equation}
\begin{equation}
\psi_r^2e(\bi) = 
\left\{
\begin{array}{ll}
0&\text{if $i_r = i_{r+1}$},\\
e(\bi)&\text{if $i_r \nslash i_{r+1}$},\\
(y_{r+1}-y_r)e(\bi)&\text{if $i_r \rightarrow i_{r+1}$},\\
(y_r - y_{r+1})e(\bi)&\text{if $i_r \leftarrow i_{r+1}$},\\
(y_{r+1} - y_{r})(y_{r}-y_{r+1}) e(\bi)\!\!\!&\text{if $i_r \rightleftarrows i_{r+1}$};
\end{array}
\right.
 \label{R4}
 \end{equation}
\begin{equation}
\psi_r \psi_s = \psi_s \psi_r\qquad\text{if $|r-s|>1$};\label{R3Psi}
\end{equation}
 \begin{equation}
\psi_{r}\psi_{r+1} \psi_{r} e(\bi)
=
\left\{\begin{array}{ll}
(\psi_{r+1} \psi_{r} \psi_{r+1} +1)e(\bi)\quad&\text{if $i_{r+2}=i_r \rightarrow i_{r+1}$},\\
(\psi_{r+1} \psi_{r} \psi_{r+1} -1)e(\bi)\quad&\text{if $i_{r+2}=i_r \leftarrow i_{r+1}$},\\
\big(\psi_{r+1} \psi_{r} \psi_{r+1} -2y_{r+1}
\\\qquad\:\quad +y_r+y_{r+2}\big)e(\bi)
\quad&\text{if $i_{r+2}=i_r \rightleftarrows i_{r+1}$},\\
\psi_{r+1} \psi_{r} \psi_{r+1} e(\bi)&\text{otherwise}.
\end{array}\right.
\label{R7}
\end{equation}
\begin{equation}
y_1^{\de_{i_1,0}}e(\bi)=0;\label{ERCyc}
\end{equation}
\end{Theorem}

Theorem~\ref{TBK}, as well as its generalizations and refinements given in Theorems~\ref{TBKCyc} and \ref{TBKBlocks}, establish an isomorphism between the Hecke algebras we are interested in and cyclotomic Khovanov-Lauda-Rouquier algebras to be discussed later on in this article. 

A remarkable feature of the given presentation is that it does not contain the parameter $\xi$. Rather, $\xi$ comes in indirectly---it determines $e$, which in turn determines the `Lie type' $\Gamma$.  

Note that if $\xi=1$ and $\operatorname{char}F=p>0$, then $e=p$. On the other hand, if $F=\C$ and $\xi=e^{2\pi i/p}$ then again $e=p$. So the `Lie type' is the same, and we get the same relations in these two cases (but over different fields). This observation will be used in section~\ref{SSRedNew} to define a reduction modulo $p$ procedure.

Another important feature of our presentation is that it is obviously homogeneous with respect to the following grading:

\begin{Corollary}
There is a unique $\Z$-grading on
$H_d$ such that
$$
\deg(e(\bi))= 0,\quad \deg(y_r)=2, \quad\deg(\psi_r e(\bi))=-a_{i_r,i_{r+1}}
$$
for all admissible $r$ and $\bi$.
\end{Corollary}

\subsection{Affine Hecke algebras}\label{SSCHA}
Let $H_d^\aff$ denote
the {\em affine Hecke algebra} over the ground field $F$ associated to $\Sigma_d$
if $\xi \neq 1$, or its rational
degeneration if $\xi=1$ \cite{Drinfeld,Lu}.
Thus, if $\xi\neq 1$, then 
$H_d^\aff$ is the $F$-algebra generated by 
$$T_1,\dots,T_{d-1},X_1^{\pm 1}, \dots, X_d^{\pm 1}$$
subject only to the relations (\ref{QECoxeterQuad})--(\ref{QCoxeterFar}) and the relations 
\begin{eqnarray}
X_r^{\pm1} X_s^{\pm1}=X_s^{\pm1} X_r^{\pm1}\qquad(1\leq r,s\leq d), \\ 
X_rX_r^{-1}=1\qquad(1\leq r\leq d),\\
T_r X_r T_r = \xi X_{r+1}\qquad(1\leq r< d),\\
T_r X_s = X_{s}T_r \qquad(1\leq r<d,1\leq s\leq d,s \neq r,r+1).
\end{eqnarray}
If $\xi = 1$, then 
$H_d^\aff$ is the $F$-algebra generated by 
$$T_1,\dots,T_{d-1}, X_1,\dots,X_d$$ subject only to the relations (\ref{QECoxeterQuad})--(\ref{QCoxeterFar}) and the relations:
\begin{eqnarray}
X_rX_s=X_sX_r\qquad(1\leq r,s\leq d),
\\
\label{EDAHA}
T_r X_{r+1} = X_r T_r + 1\qquad(1\leq r< d),\\
\label{EDAHA2}
 T_r X_s = X_s T_r \qquad(1\leq r<d,1\leq s\leq d,s \neq r,r+1).
\end{eqnarray}

One motivation for introducing affine Hecke algebras  is as follows. Consider for example the case $\xi=1$. In section~\ref{SSSGIHA}, we have mentioned the idea of using the Gelfand-Zetlin subalgebra $\langle L_1,\dots,L_d\rangle\subset H_d$ as a `Cartan subalgebra' of $H_d$. One problem with this approach is that the Gelfand-Zetlin subalgebra is in general rather complicated. It would be much nicer to `free' the generators $L_1,\dots,L_d$ and consider algebraically independent commuting elements $X_1,\dots,X_d$ instead. This has to be done `outside' of $H_d$, so we should form the tensor product $H_d\otimes F[X_1,\dots,X_d]$. However, we want to preserve the relations which we had between the Jucys-Murphy elements $L_1,\dots,L_d$ and the standard generators $T_1,\dots,T_{d-1}$ of $H_d$. The definition  of Jucys-Murhy elements implies the relations $T_r L_{r+1} = L_r T_r + 1$ and $T_r L_s = L_s T_r$ for $s \neq r,r+1$. This explains the relations (\ref{EDAHA}) and (\ref{EDAHA2}). 

\subsection{Cyclotomic Hecke algebras}\label{SSCHANew}
We now introduce the main class of algebras we are going to work with 

\begin{Definition}\label{DCHA}
{\rm 
Let $\La\in P_+$ be a dominant weight as in (\ref{ELa}). 
The {\em cyclotomic Hecke algebra} $H_d^\La=H_d^\La(F,\xi)$ 
is the quotient 
\begin{equation}\label{ECHA}
H_d^\La := H_d^\aff \Big/ \big\langle \,\textstyle\prod_{i\in I}(X_1-\nu(i))^{(\La,\al_i)}\,\big\rangle= H_d^\aff\Big/\langle \,\textstyle\prod_{m=1}^l (X_1-\nu(k_m))\,\big\rangle.
\end{equation}
}
\end{Definition}

Let us make it explicit again that the algebra $H_d^\La$ depends on the ground field $F$, the `rank' $d\in\Z_{\geq 0}$, the parameter $\xi$ (which determines $e$ and the `Lie type' $\Gamma$), and the dominant weight $\La$ for Lie type $\Gamma$. In this paper we will mainly work in the generality of cyclotomic Hecke algebras.

The algebra $H_d^\La$ can be thought of as the Hecke algebra of the complex reflection group of type $G(l,1,d)$. 
If the weight $\La$ is of level $l$, we say that $H_d^\La$ is a  cyclotomic Hecke algebra of {\em level $l$}. 
The algebra $H_d$ appears as a cyclotomic Hecke algebra of level $1$:
\begin{equation}\label{EL1}
H_d\cong H_d^{\La_i}\qquad(i\in I).
\end{equation}
This is an easy consequence of the Basis Theorem~\ref{TBasisCyc} below.

It is easy to see that there exists an antiautomorphism $*$ of $H_d^\La$ defined on the generators by
\begin{equation}\label{EStar}
*:H_d^\La\to H_d^\La,\quad T_r\mapsto T_r,\ X_s\mapsto X_s\qquad(1\leq r<d,1\leq s\leq d).
\end{equation}

The presentation from Theorem~\ref{TBK} generalizes to the whole class of cyclotomic Hecke algebras. Namely in \cite{BKyoung} we construct explicit elements 
\begin{equation}\label{EKLGensCyc}
\{e(\bi)\:|\: \bi\in I^d\}\cup\{y_1,\dots,y_{d}\}\cup\{\psi_1, \dots,\psi_{d-1}\}
\end{equation}
of $H_d^\La$ and prove:

\begin{Theorem}\label{TBKCyc}
The algebra $H_d^\La$ is generated by the elements (\ref{EKLGensCyc}) 
subject only to the relations (\ref{R0})--(\ref{R7}) 
and one  additional relation 
\begin{equation}\label{ECyclotomic}
y_1^{(\La,\al_{i_1})}e(\bi)=0\qquad(\bi=(i_1,\dots,i_d)\in I^d).
\end{equation}
\end{Theorem}

Again, as in level $1$, we can now get a grading on our algebra:

\begin{Corollary}
There is a unique $\Z$-grading on
$H_d^\La$ such that
\begin{equation}\label{EDeg}
\deg(e(\bi))= 0,\quad \deg(y_r)=2, \quad\deg(\psi_r e(\bi))=-a_{i_r,i_{r+1}}
\end{equation}
for all admissible $r$ and $\bi$.
\end{Corollary}

Just like for $H_d$, a surprising feature of the presentation given in Theorem~\ref{TBKCyc} is that it does not contain the parameter $\xi$.  One corollary of this can already be stated here; this observation will also be used in section~\ref{SSRedNew} to define a reduction modulo $p$ procedure for cyclotomic Hecke algebras.

\begin{Corollary}
Suppose that $F$ is of characteristic zero.
Then the algebra
$H^\La_d$ for $\xi$ not a root of unity
is isomorphic to the
algebra $H^\La_d$ for $\xi=1$.
In other words, the
  cyclotomic Hecke algebra for generic $\xi$
is isomorphic to its rational degeneration.
\end{Corollary}

Theorem~\ref{TBKCyc} shows that the algebra $H_d^\La$ possesses a {\em graded} anti-automorphism
\begin{equation}\label{star}
\circledast:H_d^\La \rightarrow H_d^\La,\quad e(\bi)\mapsto e(\bi),\ y_r\mapsto y_r,\  \psi_s\mapsto \psi_s
\end{equation}
for all admissible $r,s$ and $\bi$. 
We write $x^{\circledast}$ for the image of the element $x\in H_d^\La$ under $\circledast$. 
Using this we introduce a {\em graded duality} 
on $\Rep{H^\La_d}$,
mapping a module $M$ to 
$$M^\circledast := \HOM_F(M, F)$$ with the
action defined by 
$$(xf)(m) = f(x^\circledast m)\qquad(m\in M,\ f\in M^\circledast,\ x\in H_d^\La).$$

Let $w\in \Sigma_d$. Pick any reduced decomposition $w=s_{r_1}\dots s_{r_\ell}$ in $\Sigma_d$, and define 
\begin{equation}\label{ETW}
T_w:=T_{r_1}\dots T_{r_\ell}\in H_d^\La.
\end{equation}
By Matsumoto's Theorem on reduced decompositions, the element $T_w$ is well-defined. Moreover:

\begin{Theorem}\label{TBasisCyc} {\rm \cite{AK}, \cite[Theorem 7.5.6]{Kbook}} 
$$\{T_wX_1^{a_1}\dots X_d^{a_d}\mid w\in \Sigma_d,\ 0\leq a_1,\dots,a_d< l\}$$ 
is a basis of the cyclotomic Hecke algebra $H_d^\La$ of level $l$. 
In particular, 
$$
\dim H_d^\La=l^d d!.
$$
\end{Theorem}


\subsection{Blocks}\label{SSBlocks}
The group $\Sigma_d$ acts on the left
on the set 
$ I^d$ by place
permutation. 
The $\Sigma_d$-orbits on $I^d$ are the sets
\begin{equation*}
I^\alpha := \{\bi=(i_1,\dots,i_d) \in I^d\:|\:\alpha_{i_1}+\cdots+\alpha_{i_d} = \alpha\}
\end{equation*}
parametrized by all $\alpha \in Q_+$ of height $d$. 
 Let 
\begin{equation}\label{bal}
e_{\alpha} := \sum_{\bi \in I^\alpha} e(\bi) \in H^\La_d.
\end{equation}
As a consequence of \cite{LM} or \cite[Theorem~1]{cyclo},
$e_{\alpha}$ is either zero or it is a primitive central idempotent
in $H^\La_d$.
Hence the algebra
\begin{equation}\label{fe}
H^\La_\alpha := e_{\alpha} H^\La_d
\end{equation}
is either zero or it is a single {\em block}
of the algebra $H^\La_d$. For $h\in H_d^\La$ let us write $h$ again for the element 
$he_\al\in H_\al^\La$. Then we get generators 
\begin{equation}\label{EBlockGenerators}
\{e(\bi)\mid \bi\in I^\al\}\cup\{y_1,\dots,y_d\}\cup\{\psi_1,\dots,\psi_{d-1}\}
\end{equation}
for $H_\al^\La$. 

The presentation of Theorem~\ref{TBK} can be refined to produce a graded presentation of the individual blocks $H_\alpha^\La$:

\begin{Theorem} \label{TBKBlocks}
{\rm \cite{BKyoung}} 
The block algebra $H_\al^\La$ is generated by the elements (\ref{EBlockGenerators})  subject 
only to the  relations (\ref{R0})--(\ref{R7}) and (\ref{ECyclotomic}) 
for all 
$\bi,\bj\in I^\al$ and all 
admissible $r, s$. In particular, 
there is a unique $\Z$-grading on
$H_\al^\La$ such that
\begin{equation}\label{EBlockDeg}
\deg(e(\bi))= 0,\quad \deg(y_r)=2, \quad\deg(\psi_r e(\bi))=-a_{i_r,i_{r+1}}
\end{equation}
for all $\bi\in I^\al$ and all admissible $r$.
\end{Theorem}

We point out that Theorem~\ref{TBKBlocks} holds even for those $\al$ for which $H_\al^\La=0$. However, it is difficult to see from our presentation when $H_\al^\La=0$.

\subsection{\boldmath The main problem}\label{SSqCh}
Let ${\mathscr C}$ be the free $\Laurent$-module on $I^d$. If $M\in\Rep{H_d^\La}$, then the {\em $q$-character} of $M$ is  the formal expression 
\begin{equation}\label{qch}
\CH M := 
\sum_{\bi \in I^\alpha} (\qdim\ e(\bi) M) \cdot \bi\in{\mathscr C}.
\end{equation}
At combinatorial level this goes back 
to Leclerc \cite{Lec}.

The main problem in representation theory of $H_d^\La$ can now be stated:

\vspace{2 mm}
\noindent
{\bf Main Problem.} Classify irreducible graded modules over $H_d^\La$ and describe their $q$-characters. 

\vspace{2mm}

The following theorem is established in {\cite[Theorem 3.17]{KL1}}, going back to Bernstein at the ungraded level.

\begin{Theorem}\label{ch}
The map $$
\CH:[\Rep{R}] \rightarrow
{\mathscr C}, \quad [M] \mapsto \CH M
$$ 
is injective.
\end{Theorem}


\subsection{Affine Khovanov-Lauda-Rouquier algebras}\label{SSAKLRA}
If we drop the cyclotomic relation 
$$\prod_{i\in I}(X_1-\nu(i))^{(\La,\al_i)}=0$$ 
from the definition of the cyclotomic Hecke algebra $H_d^\La$, we get the affine Hecke algebra $H_d^\aff$. On the other hand, let us drop the cyclotomic relation (\ref{ECyclotomic}) from the presentation of $H_d^\La$ obtained in Theorem~\ref{TBKCyc}. Then we get the affine Khovanov-Lauda-Rouquier algebra:

\begin{Definition}
{\rm 
The {\em affine Khovanov-Lauda-Rouquier algebra} $H_d^\infty$ (of type $\Gamma$) is the graded $F$-algebra given abstractly by generators 
\begin{equation}\label{EAffKLGens}
\{e(\bi)\:|\: \bi\in I^d\}\cup\{y_1,\dots,y_{d}\}\cup\{\psi_1, \dots,\psi_{d-1}\}
\end{equation}
and relations (\ref{R0})--(\ref{R7}), with the degrees of the generators  given by the formulas (\ref{EDeg}). 
}
\end{Definition}

The affine Khovanov-Lauda-Rouquier algebras were introduced independently by Khovanov-Lauda \cite{KL1,KL2} and Rouquier \cite{Ro}. The algebra $H_d^\infty$ is {\em not} isomorphic to the affine Hecke algebra $H_d^\aff$, although the two are isomorphic  after some  completions, cf. \cite{Ro}.

As in the cyclotomic case, for each $\al\in Q_+$ with $\height(\al)=d$, define the idempotent $e_\al$ by the formula (\ref{bal}). 
It follows from \cite[Corollary 2.11]{KL1} that 
$e_{\alpha}$ is a primitive central idempotent
in $H_d^\infty$.
Hence the algebra
$$
H^\infty_\alpha := e_{\alpha} H^{\infty}_d
$$ 
 is a single {\em block}
of the algebra $H^\infty_d$. The defining presentation of $H_d^\infty$ can be refined to produce a graded presentation of the blocks $H_\al^\infty$:

\begin{Theorem} \label{TKLBlocks}
The block algebra $H_\al^\infty$ is generated by elements 
\begin{equation}\label{EAffBlockGenerators}
\{e(\bi)\mid \bi\in I^\al\}\cup\{y_1,\dots,y_d\}\cup\{\psi_1,\dots,\psi_{d-1}\}
\end{equation}
subject 
only to the  relations (\ref{R0})--(\ref{R7}) for 
$\bi,\bj\in I^\al$ and all 
admissible $r, s$. In particular, 
there is a unique $\Z$-grading on
$H_\al^\infty$ such that 
$$\deg(e(\bi))= 0,\quad \deg(y_r)=2, \quad\deg(\psi_r e(\bi))=-a_{i_r,i_{r+1}}
$$
for all $\bi\in I^\al$ and all admissible $r$.
\end{Theorem}

All algebras $H_\al^\infty$ are non-zero. In fact, unlike their cyclotomic quotients, they have a nice `PBW-type' bases \cite[Theorem 2.5]{KL1}. 

There is a duality denoted $\#$
on $\Proj{H^\infty_\alpha}$ mapping a projective module $P$
to 
\begin{equation}\label{ESharp}
P^\# := \HOM_{H^\infty_\alpha}(P, H^\infty_\alpha),
\end{equation}
with the action defined by 
$$(xf)(p) = f(p) x^\circledast\qquad(x\in H^\infty_\al,f\in  \HOM_{H^\infty_\alpha}(P, H^\infty_\alpha), p\in P).$$

If, instead of the quiver $\Gamma$, chosen in section~\ref{SSDynkin}, one considers any simply laced quiver $Q$ without loops with the set of vertices $I$, generators (\ref{EAffKLGens}) and relations (\ref{R0})--(\ref{R7}) define the affine Khovanov-Lauda-Rouquier algebra $R_d(Q)$ of type $Q$. 

If there are no cycles in the underlying graph of $Q$, then the affine Khovanov-Lauda-Rouquier algebras obtained using different orientations are all isomorphic,  but if there are cycles this is not the case, cf. \cite{KL2}. For non-simply-laced $Q$, the algebras $R_d(Q)$ are not much harder to define, see \cite{KL2} and \cite{Ro}. Note that one non-simply-laced case, namely $A_1^{(1)}$, already occurs in Theorems~\ref{TBK} and \ref{TBKCyc}---it corresponds to  $e=2$.

\section{Combinatorics}\label{SComb}

In this section we fix our notation concerning multipartitions and related combinatorial objects. 
Throughout the section we continue working with a fixed positive integer $l$, referred to as the level, a tuple $\kappa=(k_1,\dots,k_l)$ as in (\ref{EKappa}), and the corresponding dominant weight $\La=\La_{k_1}+\dots+\La_{k_l}$ as in (\ref{ELa}). The reader, who is only interested in the case of symmetric groups or Iwahori-Hecke algebras, can assume that $\La=\La_0$, in particular $l=1$, and all $l$-multipartitions appearing below will become usual partitions.

\subsection{Partitions and Young diagrams}\label{SSPar}
An {\em $l$-multipartition} of $d$
is an ordered $l$-tuple of partitions 
$\mu = (\mu^{(1)} , \dots,\mu^{(l)})$
such that $\sum_{m=1}^l |\mu^{(m)}|=d$. 
We call $\mu^{(m)}$ the {\em $m$th component} of $\mu$. 
The set of all $l$-multipartitions of $d$ is denoted $\Par_d$, and we put 
$$\Par:=\bigsqcup_{d\geq 0}\Par_d.$$ 
This notation might seem excessive at the moment as the set $\Par$ really only depends on $l$, and not on $\kappa$. However, once we start considering residues of nodes of multipartitions, dependence on $\kappa$ will become crucial. 

The {\em Young diagram} of the multipartition $\mu = (\mu^{(1)} , \dots,\mu^{(l)})\in \Par$ is 
$$
\{(a,b,m)\in\Z_{>0}\times\Z_{>0}\times \{1,\dots,l\}\mid 1\leq b\leq \mu_a^{(m)}\}.
$$
The elements of this set
are called the {\em nodes of $\mu$}. More generally, a {\em node} is an element of $\Z_{>0}\times\Z_{>0}\times \{1,\dots,l\}$. 
Usually, we identify the multipartition $\mu$ with its 
Young diagram and visualize it as a column vector of Young diagrams. 
For example, $((3,1),\emptyset,(4,2))$ is the Young diagram
$$
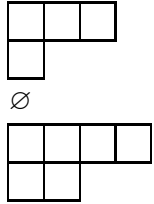
\begin{array}{l}
\diagram{&&\cr \cr}\\\emptyset\\\diagram{&&&\cr &\cr}
\end{array}
$$ 

To each node $A=(a,b,m)$ 
we associate its {\em residue}, which is defined to be the following element of $I=\Z/e\Z$: 
$$\res\, A:=k_m+(b-a)\pmod{e}.$$ 

We refer to the nodes of residue $i$ as the {\em $i$-nodes}.
Define the {\em residue content of $\mu$} to be
\begin{equation}\label{EContent}
\cont(\mu):=\sum_{A\in\mu}\al_{\res A} \in Q_+.
\end{equation}
Denote
$$
\Par_\al:=\{\mu\in\Par\mid \cont(\mu)=\al\}\qquad(\al\in Q_+).
$$

A node $A\in\mu$ is called {\em removable (for $\mu$)}\, if $\mu\setminus \{A\}$ has a shape of a multipartition. A node $B\not\in\mu$ is called {\em addable (for $\mu$)}\, if $\mu\cup \{B\}$ has a shape of a multipartition. We use the notation
$$
\mu_A:=\mu\setminus \{A\},\qquad \mu^B:=\mu\cup\{B\}.
$$
For the example above, the removable nodes are $(1,3,1)$, $(2,1,1)$, $(1,4,3)$, $(2,2,3)$, and addable nodes are $(1,4,1)$, $(2,2,1)$,  $(3,1,1)$, $(1,1,2)$, $(1,5,3)$, $(2,3,3)$, $(3,1,3)$, in order from top to bottom in the diagram.

Let $\mu,\nu\in \Par_d$. We say that $\mu$ {\em dominates} $\nu$, written $\mu\unrhd\nu$, if 
$$
\sum_{a=1}^{m-1}|\mu^{(a)}|+\sum_{b=1}^c\mu_b^{(m)}\geq 
\sum_{a=1}^{m-1}|\nu^{(a)}|+\sum_{b=1}^c\nu_b^{(m)}
$$
for all $1\leq m\leq l$ and $c\geq 1$.
In other words, $\mu$ is obtained from $\nu$ by moving nodes up in
the diagram.

Let $\lex$ denote the {\em lexicographic ordering} on partitions, so
for partitions $\la,\mu \in \mathscr{P}$ we have that
$\mu\lex\la$ if and only if $\mu_1=\la_1,\dots,\mu_{a-1}=\la_{a-1}$, 
and $\mu_a < \la_a$ for some $a \geq 1$. 
We extend this notion to $l$-multipartitions: for 
$\la, \mu \in \Par$ we have that $\mu\lex\la$ if and only 
$\mu^{(1)} = \la^{(1)},\dots,\mu^{(m-1)}=\la^{(m-1)}$, and
$\mu^{(m)} \lex \la^{(m)}$ for some $1 \leq m \leq l$. 
This total order refines the dominance ordering on $\Par$ 
in the sense
that $\mu \lhd \la$ implies $\mu \lex \la$.

\subsection{Tableaux}\label{SSTab}
Let $\mu=(\mu^{(1)},\dots,\mu^{(l)})\in\Par_d$. 
A {\em $\mu$-tableau} 
$$\T=(\T^{(1)},\dots,\T^{(l)})$$ is obtained from the diagram of $\mu$ by 
inserting the integers $1,\dots,d$ into the nodes, allowing no repeats. 
The tableaux $\T^{(m)}$ are 
called the {\em components} of $\T$.
To each tableau $\T$ we associate its {\em residue sequence}
\begin{equation}\label{EResSeq}
\bi^\T=(i_1,\dots,i_d)\in I^d,
\end{equation}
where $i_r$ is the residue of the node occupied by 
$r$ in $\T$ ($1\leq r\leq d$). 

A $\mu$-tableau $\T$ is {\em row-strict} (resp. {\em column-strict}) if its entries increase from left to right (resp. from top to bottom) along the rows (resp. columns) of each component of $\T$. 
A $\mu$-tableau $\T$ is {\em standard} if it is row- and column-strict. 
The set of all standard $\mu$-tableaux will be denoted by $\St(\mu)$.  The group $\Sigma_d$ acts on the set of $\mu$-tableaux on the left by its action on
the entries.

Let $\T^\mu$ be the $\mu$-tableau in which the numbers $1,2,\dots,d$ appear in order along the successive rows, working from top to bottom.
Set 
\begin{equation}\label{EBIMu}
\bi^\mu:=\bi^{\T^\mu}.
\end{equation} 
If $\T$ is a $\mu$-tableau, then $w_\T\in \Sigma_d$ is defined from
\begin{equation}\label{EWT}
w_\T  \T^\mu=\T.
\end{equation}

\begin{Example}
{\rm 
Let $\mu=((3,1),\emptyset,(4,2))$, $e=3$, and $k_1=0,k_2=1,k_3=1$. The following are examples of standard $\mu$-tableaux:
$$
\T^\mu=
\begin{array}{l}
\diagram{1&2&3\cr 4\cr}\\\emptyset\\\diagram{5&6&7&8\cr 9&10\cr}
\end{array}
\qquad\qquad \T=\begin{array}{l}
\diagram{2&5&6\cr3 \cr}\\\emptyset\\\diagram{1&4&9&10\cr 7&8\cr}\end{array}
$$
Then 
\begin{align*}
\bi^\mu&=(0,1,2,2,1,2,0,1,0,1), \\
\bi^\T&=(1,0,2,2,1,2,0,1,0,1),\\
w_\T&=(1\:2\:5)(3\:6\:4)(7\: 9)(8\:10).
\end{align*} 
}
\end{Example}

\vspace{1 mm}
Let $\mu\in\Par_d$. Recalling the Bruhat order `$\leq$' on $\Sigma_d$,
the {\em Bruhat order} on $\St(\mu)$ is defined as follows: 
\begin{equation}\label{EBruhat}
\Stab\leq\T\quad\text{if and only if}\quad w_\Stab\leq w_\T.
\end{equation}

\subsection{Degree of a standard tableau}\label{SSDeg}
Let $\mu\in\Par$, $i\in I$, $A$ be a removable $i$-node and $B$ be an addable $i$-node of 
$\mu$. We set
\begin{equation}\label{EDMUA}
\begin{split}
d_A(\mu):= &\#\{\text{addable $i$-nodes of $\mu$ strictly below $A$}\}
\\
&-\#\{\text{removable $i$-nodes of $\mu$ strictly below  $A$}\};
\end{split}
\end{equation} 
\begin{equation}\label{EDMUB}
\begin{split}
d^B(\mu):=&\#\{\text{addable $i$-nodes of $\mu$ strictly above $B$}\}\\
&-\#\{\text{removable $i$-nodes of $\mu$ strictly above $B$}\};
\end{split}
\end{equation} 
\vspace{.5 mm}
\begin{equation}\label{EDKWeight}
d_i(\mu):=\#\{\text{addable $i$-nodes of $\mu$}\}
-\#\{\text{removable $i$-nodes of $\mu$}\}.
\end{equation} 
It is easy to see \cite[Lemma 3.11]{BKW} that:
\begin{equation}\label{EDI}
d_i(\mu)=(\La-\al,\al_i)\qquad(\mu
\in \Par_\al).
\end{equation}
Finally, for $\al\in Q_+$, define the {\em defect} of $\al$ (relative to $\La$) to be
\begin{equation}\label{defdef}
\defect(\al)=(\La,\al)-(\al,\al)/2.
\end{equation}

Given $\mu \in \Par_d$ and $\T \in \St(\mu)$,
the {\em degree} of $\T$ is defined in \cite[section~3.5]{BKW} inductively as follows. If $d=0$, then 
$\T$ is the empty tableau $\emptyset$, and we set
$\deg(\T):=0$.
Otherwise, let $A$ be the node occupied by $d$ in $\T$. Remove this node to get the standard tableaux $\T_{\leq(d-1)}\in\St(\mu_A)$, and set
\begin{equation}\label{EDegTab}
\deg(\T):=d_A(\mu)+\deg(\T_{\leq(d-1)}).
\end{equation}

\subsection{Good nodes and restricted multipartitions}\label{SSCrGraph}

Let $\mu\in\Par$. 
Set
\begin{equation}\label{EWt}
\wt(\mu) = \La - \cont(\mu)\qquad(\mu\in\Par).
\end{equation}
Further, given  $i \in I$, let $A_1,\dots,A_n$ denote the addable and removable $i$-nodes
of $\mu$ ordered so that $A_m$ is above $A_{m+1}$ for each
$m=1,\dots,n-1$.
Consider the sequence $(\sigma_1,\dots,\sigma_n)$
where $\sigma_r = +$ if $A_r$ is addable or $-$ if $A_r$ is removable.
If we can find $1 \leq r < s \leq n$ such that $\sigma_r=  -$,
$\sigma_s = +$ and $\sigma_{r+1}=\cdots=\sigma_{s-1} = 0$
then replace $\sigma_r$ and $\sigma_s$ by $0$.
Keep doing this until we are left with a sequence
$(\sigma_1,\dots,\sigma_n)$ in which no $-$ appears to the left of
a $+$. This
is called the {\em reduced $i$-signature} of $\mu$. 

Now, let $(\sigma_1,\dots,\sigma_n)$ be the reduced $i$-signature of $\mu$. Set 
\begin{align*}
\eps_i(\mu) = \#\{r=1,\dots,n\:|\:\sigma_r = -\},\qquad
&\phi_i(\mu) = \#\{r=1,\dots,n\:|\:\sigma_r = +\}.
\end{align*}

If $\eps_i(\mu) > 0$, let $r$ index the leftmost $-$ in the reduced $i$-signature. In this case the removable node $A_r$ is called the {\em good} $i$-node of $\mu$ and we set 
$$\tilde e_i \mu = \mu_{A_r}.$$ 
Similarly, if $\phi_i(\mu) >0 $ let $s$ index the rightmost $+$ in the reduced $i$-signature. The addable node $A_s$ is called the {\em cogood} $i$-node for $\mu$, and  we set 
$$\tilde f_i \mu = \mu^{A_s}.$$ 

It is easy to check that a removable node $A$ is good for $\mu$ if and only if $A$ is cogood for $\mu_A$. Similarly, an addable node $B$ is cogood for $\mu$ if and only if $B$ is good for $\mu^B$. Hence 
\begin{equation}
\mu=\tilde f_i\nu \quad\text{if and only if}\quad \nu=\tilde e_i\mu
\qquad\qquad(\mu,\nu\in\Par). 
\end{equation}


Define the set $\RPar$ of the {\em ($\kappa$-)restricted multipartitions} to be the set of all multipartitions $\mu\in\Par$ which can be obtained from the empty multipartition $\emptyset \in\Par$ by applying a series of $\tilde f_i$'s:
\begin{equation}\label{ERPar}
\begin{split}
\RPar&=\{\tilde f_{i_1}\dots \tilde f_{i_m}\emptyset\mid i_1,\dots,i_m\in I\}\\
&=\{\mu\in\Par\mid \tilde e_{j_1}\dots\tilde e_{j_m}\mu=\emptyset\ \text{for some $j_1,\dots,j_m\in I$}\}. 
\end{split}
\end{equation}
Also for $\alpha \in Q_+$ set
\begin{equation*}
\RPar_\alpha := \RPar \cap \Par_\alpha.
\end{equation*}

\begin{Example}\rm\label{c3}
(i) Suppose that $e > 0$ and $l=1$.
Then $\RPar$
is the set of all
{\em $e$-restricted partitions}, that is, 
partitions $\mu$ such that 
$\mu_a - \mu_{a+1} < e$ for all $a \geq 1$.

(ii)
Suppose that $e=0$ and 
$k_1 \geq \cdots \geq k_l$.
Then $\RPar$
consists of all $l$-multipartitions $\mu$
such that $\mu_{a+k_m - k_{m+1}}^{(m)}
\leq \mu_a^{(m+1)}$ for all $m=1,\dots,l-1$ and $a \geq 1$;
see \cite[(2.52)]{BKariki} or \cite{Vaz}. 
Such multi-partitions are closely related to combinatorics of standard tableaux, cf. \cite{BKariki,BKrep}

(iii) Suppose that $e=0$ and 
$k_1 \leq \cdots \leq k_l$.
Then $\RPar$
consists of all $l$-multipartitions $\mu$
such that $\mu_a^{(m)}\leq
\mu_a^{(m+1)}+k_{m+1} - k_{m}$ for all $m=1,\dots,l-1$ and $a \geq 1$;
see \cite[(2.53)]{BKariki}.

(iv) For more examples see \cite{ArikiJac,ArKrTs}. 
\end{Example}

In general a non-recurrent description of the set $\RPar$ is not known.

\section{Solution of the Main Problem for type ${A_\infty}$ at level $1$}\label{SSMPCharZerLev1}
In this small section we present what can be considered an ideal solution of the Main Problem from section~\ref{SSqCh} for the case where $e=0$ and $l=1$. In view of (\ref{EL1}) this covers two classical cases:  representation theory of symmetric groups in characteristic $0$ and representation theory of Iwahori-Hecke algebra of the symmetric group with generic parameter (both algebras are semisimple). Theorem~\ref{TL1AInf}(i) in particular gives the formulas for the generators of $H_d^\La$ acting on a basis of an arbitrary irreducible $H_d^\La$-module.  The reader might note how unusually easy these formulas are. For examples the only scalars that ever appear are $0$ and $1$. 

So assume throughout the section that $e=0$, and $\La=\La_i$ is a dominant weight of level $1$. Then the set $\Par_d$ is just the set of usual partitions of $d$, with residue of the box in row $a$ and column $b$ being $b-a+i\in I=\Z$. Define the vector spaces
$$
S( \mu):=\bigoplus_{\T\in \St( \mu)}F\cdot v_\T\qquad( \mu\in\Par_d).
$$
with basis elements labeled by the standard $ \mu$-tableaux.

\begin{Theorem}\label{TL1AInf}
Let $e=0$ and $l=1$. 
\begin{enumerate}
\item[{\rm (i)}] Let $ \mu\in\Par_d$. The formulas
\begin{align*}
e(\bj)v_\T&=\de_{\bi^\T,\bj}v_\T \qquad (\bj\in I,\ \T\in  \St(\mu)),\\
y_r v_\T&=0\qquad (1\leq r\leq d,\ \T\in  \St(\mu)),\\
\psi_rv_\T&=
\left\{
\begin{array}{ll}
v_{s_r\T} &\hbox{if $s_r\T$ is standard,}\\
0 &\hbox{otherwise;}
\end{array}
\right.
\quad(1\leq r<d,\ \T\in  \St(\mu))
\end{align*}
define an action of $H_d^\La$ on $S(\mu)$, under which $S(\mu)$ is an  irreducible $H_d^\La$-module. 
\item[{\rm (ii)}] $\{S( \mu)\mid \mu\in\Par_d\}$ is a complete and irredundant set of irreducible $H_d^\La$-modules up to isomorphism. 
\item[{\rm (iii)}] Let us consider $S( \mu)$ as a graded vector space by declaring that it is concentrated in degree $0$. Then $$\{S( \mu)\langle n\rangle\mid \mu\in\Par_d,\ n\in \Z\}$$
is a complete and irredundant set of irreducible graded $H_d^\La$-modules up to isomorphism. 
\item[\rm (iv)] Let $\mu\in\Par_d$ and $n\in\Z$. Then 
$$
\CH S( \mu)\langle n\rangle=q^n\sum_{\T\in  \St( \mu)}\bi^\T.
$$   
\item[(v)] Let $\mu,\nu\in\Par_d$, $\T\in \St(\mu)$, and $\Stab\in \St(\nu)$. Then $\bi^\T\neq\bi^\Stab$, unless $\mu=\nu$ and $\T=\Stab$.   
\end{enumerate}
\end{Theorem}

We point out that the basis $\{v_\T\mid\T\in \St(\mu)\}$ of the irreducible module $S(\mu)$  is essentially the same as {\em Young's seminormal basis}---this is explained in detail in \cite[section 5]{BKyoung}. However, the action of the standard generators $T_r$ on Young's seminormal basis is given by some well-known and rather tricky formulas, cf. for example \cite[section 25]{Jbook} or \cite[Theorem 2.3.1]{Kbook} for the symmetric group case, or \cite{H,Ram,W} for the Hecke algebras case. It is therefore somewhat miraculous that our new homogeneous generators act by the easy formulas of Theorem~\ref{TL1AInf}(i), which moreover are the {\em same} formulas for the symmetric groups and Hecke algebras! 

Standard results on characters of symmetric groups such as the Murnaghan-Nakayama formulas follow easily from Theorem~\ref{TL1AInf}, as was noticed first by Okounkov and Vershik \cite{OV}, see also \cite[section 2.3]{Kbook}.

Another remark to be made here is that all graded irreducible modules turn out to be {\em pure}, i.e. concentrated in one degree. This is definitely special for the situation where $e=0$ and $l=1$ and does not hold in general. In some sense this phenomenon `explains' why representations of symmetric groups in characteristic zero are so `easy to understand'. The purity happens to be closely connected to the property described in Theorem~\ref{TL1AInf}(v), and has many other equivalent descriptions, cf. \cite{KR} and references therein.

One could consider more generally pure representations of Khovanov-Lauda-Rouquier algebras of any type, and they too turn out to be amenable to a beautiful combinatorial description similar to the one in Theorem~\ref{TL1AInf}---this is done in \cite{KR}. In particular, one can define shapes and skew shapes corresponding to any simply-laced type $\Gamma$. 

Then the (appropriately defined) standard tableaux of a given skew shape label a basis of an irreducible module, on which the generators act by the formulas similar to the ones in Theorem~\ref{TL1AInf}(i). 
Moreover, if our skew shape is a shape, there is a {\em hook formula for any Lie type} which computes the dimension of the corresponding irreducible module! We refer the reader to \cite{KR} for the details on that. 

On the other hand,  non-pure irreducible modules are much harder to understand. One way to approach them is via Specht modules, which will be discussed in the next section.

\section{Cyclotomic Hecke algebras as cellular algebras}\label{SCell}
Symmetric groups are known to possess a nice class of modules called Specht modules, see e.g. \cite{Jbook}. Specht modules are not irreducible in general, but they are very useful because they are defined over an arbitrary field, their dimensions and characters are well-understood, and irreducible modules arise as the simple heads of Specht modules corresponding to a certain class of partitions depending on the characteristic of the ground field. In that sense Specht modules resemble Verma modules or other  `standard' modules used in Lie theory. The goal of this section is to describe Specht modules and to explain how to grade them. 

One way to approach Specht modules and to generalize them from symmetric group algebras to cyclotomic Hecke algebras, is to use the framework of the theory of {\em cellular algebras} developed by Graham and Lehrer \cite{GL}, where Specht modules arise naturally as {\em cell modules}. 

\subsection{Review of cellular algebras} \label{SSCellGen}
Let $H$ be a finite dimensional associative $F$-algebra.

\begin{Definition}\label{DCA}
{\rm 
A {\em cell datum}\, for $H$ is a quadruple $({\mathscr P},\St,C,*)$, where
\begin{enumerate}
\item[$\bullet$] ${\mathscr P}$ is a partially ordered set;
\item[$\bullet$] $\St$ is a map from ${\mathscr P}$ to finite sets;
\item[$\bullet$] $C$ is a map $\bigsqcup_{\mu\in{\mathscr P}}\St(\mu)\times\St(\mu)\to H;$ 
if $\Stab,\T\in\St(\mu)$, denote $C(\Stab,\T)$ by $C^\mu_{\Stab,\T}$; 
\item[$\bullet$] $*$ is an anti-automorphism of the $F$-algebra $H$;
\end{enumerate}
and the following axioms hold:
\begin{enumerate}
\item[{\rm (i)}] $\{C^\mu_{\Stab,\T}\mid \mu\in {\mathscr P},\ \Stab,\T\in\St(\mu)\}$ is an $F$-basis of $H$;
\item[{\rm (ii)}] $(C^\mu_{\Stab,\T})^*=C^\mu_{\T,\Stab}$;
\item[{\rm (iii)}] if $\mu\in{\mathscr P}$ and $\Stab,\T\in\St(\mu)$, then for every element $h\in H$ we have 
$$
hC^\mu_{\Stab,\T}\equiv \sum_{\U\in\St(\mu)}r_h(\U,\Stab)C^\mu_{\U,\T}\pmod{H(>\mu)},
$$
where $r_h(\U,\Stab)$ is independent of $\T$, and $H(>\mu)$ denotes the $F$-span of the basis elements $\{C^\nu_{\Stab',\T'}\mid \nu>\mu,\ \Stab',\T'\in \St(\nu)\}$. 
\end{enumerate}
The basis $\{C^\mu_{\Stab,\T}\mid \mu\in {\mathscr P},\ \Stab,\T\in\St(\mu)\}$ is referred to as a {\em cellular basis} and the algebra $H$ with a cell datum is called a {\em cellular algebra}. 
}
\end{Definition}


Given $\mu\in {\mathscr P}$, define the {\em cell module} $S(\mu)$ be a free $F$-module with basis $\{C_\T\mid\T\in \St(\mu)\}$ and the $H$-action 
$$
hC_\T=\sum_{\Stab\in\St(\mu)}r_h(\Stab,\T)C_\Stab\qquad(h\in H,\ \T\in\St(\mu)),
$$
where $r_h(\Stab,\T)$ are the same as in Definition~\ref{DCA}(iii) (and are easily checked to be uniquely defined). It is an easy consequence of the axioms that the formula above does define an $H$-module structure on $S(\mu)$. 

For any $\mu\in {\mathscr P}$, let $\phi^\mu$ be the $F$-bilinear form on $S(\mu)$ such that 
$$
\phi^\mu(C_\Stab,C_\T)=r_{C^\mu_{\Stab,\T}}(\Stab,\Stab)\qquad(\Stab,\T\in\St(\mu)).
$$
Graham and Lehrer \cite[Proposition 2.4]{GL} prove that the form $\phi^\mu$ is symmetric and invariant in the following sense:
$$
\phi^\mu(hx,y)=\phi^\mu(x,h^*y)\qquad(h\in H,\ x,y\in S(\mu)).
$$
The invariance of $\phi^\mu$ implies that its radical $R(\mu)$ 
is an $H$-submodule of $S(\mu)$. Denote
$$
D(\mu):=S(\mu)/R(\mu)\qquad(\mu\in {\mathscr P}).
$$
In practice, it is often not easy to determine when $D(\mu)$ is non-zero. However, we have the following fundamental theorem on cellular algebras:

\begin{Theorem}\label{TGL} {\rm \cite{GL}} 
Let $H$ be a cellular algebra as above, and put $${\mathscr{RP}}:=\{\mu\in{\mathscr P}\mid D(\mu)\neq 0\}.$$ 
\begin{enumerate}
\item[{\rm (i)}] $
\{D(\mu)\mid \mu\in{\mathscr{RP}}\}
$
is a complete and irredundant set of irreducible $H$-modules up to isomorphism. 
\item[{\rm (ii)}] The head 
of $S(\mu)$ is isomorphic to $D(\mu)$ for any $\mu\in{\mathscr{RP}}$. 
\item[{\rm (iii)}] The multiplicity $$d_{\mu\nu}:=[S(\mu):D(\nu)]$$ is $0$ unless $\nu\leq \mu$; moreover $d_{\mu\mu}=1$ for any $\mu\in {\mathscr{RP}}$. 
\end{enumerate}
\end{Theorem}

The numbers $d_{\mu\nu}$ appearing in Theorem~\ref{TGL}(ii) are called the {\em decomposition numbers} and the matrix $$D=(d_{\mu\nu})_{\mu\in{\mathscr P},\nu\in{\mathscr{RP}}}$$ is called the {\em decomposition matrix} (corresponding to the given cellular structure). Theorem~\ref{TGL} shows that the decomposition matrix is unitriangular (although not square in general). 

Hence, knowing decomposition numbers, one can write the class $[D(\mu)]$ of an arbitrary irreducible $H$-module $D(\mu)$ in the Grothendieck group of finite dimensional $H$-modules as a linear combination of the classes of the cell modules $[S(\nu)]$ for $\nu\leq \mu$. In particular, one can compute the dimension of $D(\mu)$ in terms of the dimensions of cell modules which are given by $\dim S(\mu)=|\St(\mu)|$. 


\subsection{Cellular structures on cyclotomic Hecke algebras}\label{SSCelCyc}
Following the work of Murphy \cite{MurphyCell} in level $1$, Dipper, James and Mathas \cite{DJM} have exhibited cellular algebra structures on all cyclotomic Hecke algebras $H_d^\La$. Versions of these cellular structures for degenerate cyclotomic Hecke algebras are explained for example in \cite[section 6]{AMR}. 

To describe a Dipper-James-Mathas cellular structure on $H_d^\La$, 
first of all take 
$
{\mathscr P}$ to be $\Par_d$, the set of $l$-multipartitions of $d$ defined in section~\ref{SSPar}. Next, for any $\mu\in\Par_d$, we take $\St(\mu)$ to be the set of standard $\mu$-tableaux defined in section~\ref{SSTab}. Let `$*$' be the anti-automorphism defined in (\ref{EStar}).

To define a cellular basis, fix for the moment $\mu=(\mu^{(1)},\dots,\mu^{(l)})\in\Par_d$ with $d_m:=|\mu^{(m)}|$ for $m=1,\dots,l$. The {\em standard parabolic subgroup} $\Sigma_\mu\leq \Sigma_d$ is the row stabilizer of the leading tableaux $\T^\mu$ introduced in section~\ref{SSTab}. Recalling (\ref{ENu}) and (\ref{ETW}), set
$$
x_\mu:=\prod_{m=2}^l\prod_{r=1}^{d_1+\dots+d_{m-1}}\big(X_r-\nu(k_r)\big)\big(\sum_{w\in \Sigma_{\mu}} T_w\big).
$$
Let $\Stab,\T\in\St(\mu)$. Recalling \ref{EWT}, define
\begin{equation}\label{CellularBasis}
C^\mu_{\Stab,\T}:=T_{w_\Stab}x_\mu T_{w_\T^{-1}}.
\end{equation}

We now have

\begin{Theorem}\label{TDJM} {\rm \cite{DJM}} 
The datum defined above is a cell datum for $H_d^\La$.
\end{Theorem} 

Two comments on this theorem. First, we point out that the cell datum given in Theorem~\ref{TDJM} depends on the fixed tuple $\kappa=(k_1,\dots,k_l)$, and not just the dominant weight $\La=\La(\kappa)$, cf. (\ref{EKappa}), (\ref{ELa}). 

Second, the cellular basis appearing in (\ref{CellularBasis}) is in general {\em not}\, homogeneous with respect to our grading on the cyclotomic Hecke algebra. 
Very recently, Hu and Mathas \cite{HuMathas} constructed {\em homogeneous cellular bases} in cyclotomic Hecke algebras, see section~\ref{SSGrCellStr}.

\subsection{Specht modules and irreducible modules}\label{SSSMIM}
Theorem~\ref{TDJM} and the general theory of cellular algebras described in section~\ref{SSCellGen} yield the family of cell modules  
\begin{equation*}
\{S(\mu)\mid \mu\in\Par_d\},
\end{equation*}
for the cyclotomic Hecke algebra $H_d^\La$, which are usually called {\em Specht modules}. By definition, each Specht module $S(\mu)$ comes with the cellular basis 
$$\{C_\T\mid\T\in\St(\mu)\}.$$

We note that the definition of the Specht module depends on our fixed tuple~$\kappa$. We also point out that in level $1$ the Specht modules defined here are actually {\em dual} to the Specht modules used in the older literature such as \cite{Jbook} and \cite{DJ}.

By Theorem~\ref{TGL}, in order to classify the irreducible $H_d^\La$-modules, we need to be able to select  the multipartitions $\mu$ for which $D(\mu):=S(\mu)/R(\mu)$ is non-zero. In level $1$ this has been done by James \cite{Jbook} for $\xi=1$ and Dipper and James \cite{DJ} for $\xi\neq 1$. In higher levels however, this turns out to be an arduous task. At the moment the only approach known to work in full generality relies on Ariki's Categorification Theorem, and so it is proved rather late in the theory. Having said that, let us state the result here, as it  fits naturally into the general theory of cellular algebras. 

Recall the class $\RPar$ of restricted partitions introduced in section \ref{SSCrGraph}. 

\begin{Theorem}\label{TIrrLab}
Let $\mu\in \Par_d$. Then $D(\mu)\neq 0$ if and only if $\mu\in\RPar_d$. In particular,
$$
\{D(\mu)\mid \mu\in\RPar_d\}
$$
is a complete and irredundant set of irreducible $H_d^\La$-modules up to isomorphism. 
\end{Theorem}

Theorem~\ref{TIrrLab} first appears in \cite{Aclass} (although compare a combinatorial claim made in the  proof of \cite[Theorem~3.4(2)]{Aclass} with examples in \cite[section 3.5]{BKllt}). Another proof is given in \cite[Theorems 5.12 and 5.17]{BKllt}. 

There is a more elementary approach to the classification of irreducible $H_d^\La$-modules based on the socle  branching rules, see \cite{Gr} and \cite{Kbook}. In that approach, discussed in detail in section~\ref{SSBrLab}, one gets hold of irreducible modules over $H_d^\La$ inductively using the socles of the restrictions of the  irreducible modules to $H_{d-1}^\La$.  Unfortunately, it is far from obvious that both approaches lead to the same labeling of irreducible modules, although this is known to be true, see section~\ref{SSBrLab} for details. 

\begin{Remark}\label{RJamD}
{\rm 
Assume that $l=1$ and $\xi=1$. Then, by Example~\ref{c3}(i), the set $\RPar_d$ consists of the $e$-restricted partitions of $d$. Recall that the Specht module $S(\mu)$ used here is dual to the Specht module $S^\mu$ used in James' book \cite{Jbook}. To see the relation between our $D(\mu)$ and the irreducible modules $D^\mu$ from \cite{Jbook}, one uses \cite[Theorem 8.15]{Jbook} to deduce $$D^\mu\cong D(\mu^t)\otimes \sign,$$ where $\sign$ is the sign representation and $\mu^t$ is the partition obtained from $\mu$ by transposing with respect to the main diagonal. A similar relation holds in the case where $l=1$ and $\xi\neq 1$, with $D^\mu$ as in \cite{DJ}. 
}
\end{Remark}

\subsection{Blocks again}\label{SSBlAg}
 Now that we have explicit constructions of Specht modules and an explicit labeling of irreducible modules we can explain how these modules split into equivalence classes according to the blocks of $H_d^\La$. Recall from section~\ref{SSBlocks} that for any $\al\in Q_+$ of height $d$ we have defined a primitive central idempotent $e_\al$, which could be zero, and the corresponding block algebra $H_\al^\La=e_\al H_d^\La$.

\begin{Theorem}
Let $\mu\in\Par_d$, $\nu\in\RPar_d$, and $\al\in Q_+$ be of height $d$. 
\begin{enumerate}
\item[{\rm (i)}] We have $e_\al\neq 0$ if and only if there is $\la\in\Par_d$ such that $\cont(\la)=\al$.
\item[{\rm (ii)}] The Specht module $S(\mu)$ belongs to the block $H_\al^\La$ if and only if $\cont(\mu)=\al$. 
\item[{\rm (iii)}]  The irreducible module $D(\nu)$ belongs to the block $H_\al^\La$ if and only if $\cont(\nu)=\al$
\end{enumerate}
\end{Theorem}

The proof of this theorem involves identifying the blocks of Specht modules, which is easy to do since their formal characters are known. In fact, it suffices to know that $S(\mu)$ is in a fixed block, and then to determine only one weight of $S(\mu)$, which can be easily done, see for example \cite[Lemma 4.1]{BKW}. 

\section{Graded modules over cyclotomic Hecke algebras}\label{SGrMod}

In this section we begin to systematically take gradings into account. In particular, we explain how to grade irreducible modules and Specht modules in a consistent way. 

\subsection{Graded irreducible modules}\label{SSGIrrMod}
Since $H_d^\La$ is a finite dimensional graded algebra, it follows from general principles explained in section~\ref{SSGRT} that the irreducible modules $D(\mu)$ are gradable uniquely up to a grading shift. It has been pointed out by Khovanov and Lauda \cite[section~3.2]{KL1} that there is always a preferred choice of the shift which makes the modules self-dual with respect to the graded duality $\circledast$ defined in section~\ref{SSCHANew}. To be more precise:

\begin{Theorem}\label{clas} \cite[Theorem 4.11]{BKllt}
For each $\mu \in \RPar$, there exists  a unique grading on $D(\mu)$ which makes it into a graded $H_d^\La$-module such that  $$D(\mu)^\circledast \cong D(\mu).$$ Moreover, the modules
$$\{D(\mu)\langle m \rangle\:|\:\mu \in \RPar, m \in \Z\}$$ 
give a complete and irredundant set of irreducible 
graded $H_d^\La$-modules.
\end{Theorem}

Next, we would like to grade Specht modules so that the natural surjection $S(\mu)\onto D(\mu)$ is a homogeneous map. This cannot be deduced from general theory, and so Specht modules have to be graded `by hand'.

\subsection{Homogeneous bases of Specht modules}\label{SSHBSM}
Let us fix a reduced decomposition
for each element $w$ of the symmetric group $\Sigma_d$:  
$$w=s_{r_1}\dots s_{r_m},$$ 
which we refer to as a {\em preferred reduced 
decomposition} of $w$. We define the elements
$$
\psi_w:=\psi_{r_1}\dots\psi_{r_m}\in H_d^\La\qquad(w\in \Sigma_d).
$$
Since the `Coxeter relation' (\ref{R7}) for $\psi_r$'s could have an `error term', in general, $\psi_w$ depends on the choice of a preferred reduced decomposition of~$w$.

Let $\mu\in \Par_d$. We want to define a basis of the Specht module $S(\mu)$ which is  better adapted to the weight space decomposition and the grading than  the cellular basis $\{C_\T\mid \T\in\St(\mu)\}$. Recall the leading tableau $\T^\mu$ and the elements $w_T$ defined in section~\ref{SSTab}. 
For any $\mu$-tableau $\T$ define the vector
\begin{equation}\label{vT}
v_\T:=\psi_{w_\T}C_{\T^\mu}\in S(\mu) .
\end{equation}
For example, $v_{\T^\mu}=C_{\T^\mu}$. 
Just like the elements $\psi_{w}$, the vectors $v_\T$ in general depend  on the choice of preferred reduced decompositions in $\Sigma_d$. Recall the Bruhat order $\leq$ on $\St(\mu)$ defined in section~\ref{SSTab}.

The following theorem describes a connection between the vectors $v_\T$ and the vectors $C_\T$, and gives some nice properties of the vectors $v_\T$ similar to those that $C_\T$ are known to possess.

\begin{Theorem}\label{PVZ} {\rm \cite{BKW}} 
Let $\mu\in\Par_d$. Then 
\begin{enumerate}
\item[{\rm (i)}] For any $\mu$-tableau $\T$ we have 
$$v_\T=\sum_{\Stab\in\St(\mu),\ \Stab\leq \T}a_{\Stab}C_\Stab\qquad(a_\Stab\in F).$$ Moreover, if $\T$ is standard, then $a_\T\neq 0$. 
\item[{\rm (ii)}] $\{v_\T\mid \T\in \St(\mu)\}$ is a basis of  $S(\mu) $. Moreover, for any $\mu$-tableau $\T$, we have
$$
v_\T=\sum_{\Stab\in \St(\mu),\ \Stab\leq \T}b_\Stab v_\Stab\qquad(b_\Stab\in F).
$$
\end{enumerate}
\end{Theorem}

The first advantage of the vectors $v_\T$ over the vectors $C_\T$ is that they are actually weight vectors:

\begin{Lemma}\label{LVWeight} {\rm \cite{BKW}} 
Let $\mu\in\Par_d$ and $\T$ be a $\mu$-tableau. Then $v_\T$ is an element of the weight space $S(\mu)_{\bi^\T}=e(\bi^\T)S(\mu) $.  
\end{Lemma}

The second advantage of the vectors $v_\T$ is that they turn out to be homogeneous with respect to a grading of $S(\mu)$ as an $H_d^\La$-module. To explain this, recall the notion of the degree $\deg(\T)$ of a standard tableaux $\T$ introduced in section~\ref{SSDeg}. 
Define the {\em degree} of $v_\T $ to be 
$$
\deg(v_\T ):=\deg(\T).
$$
As $\{v_\T \mid \T \in \St(\mu)\}$ is a basis of $S(\mu) $, this makes $S(\mu) $ into a $\Z$-graded {\em vector space}. Since the vectors $v_\T $ in general depend on the choice of preferred reduced decompositions, our grading on $S(\mu)$ might also depend on it. However, the following theorem shows that it does not!

\begin{Theorem}\label{TGrIndep} {\rm \cite{BKW}} 
Let $\mu\in\Par_d$ and $\T\in\St(\mu)$. 
If 
$$w_\T =s_{r_1}\dots s_{r_m}=s_{t_1}\dots s_{t_m}$$ are two reduced decompositions of $w_\T$, then
$$
\psi_{r_1}\dots \psi_{r_m}C_{\T^\mu}-\psi_{t_1}\dots \psi_{t_m}C_{\T^\mu}=\sum_{\Stab\in\St(\mu),\ \Stab< \T,\ \bi^\Stab=\bi^\T,\ \deg(\Stab)=\deg(\T )} a_\Stab v_\Stab
$$
for some scalars $a_\Stab\in F$. In particular, our grading on $S(\mu) $ is independent of the choice of preferred reduced decompositions. 
\end{Theorem}

The next result shows moreover that our vector space grading 
makes $S(\mu)$ into a graded $H^\La_d$-module:

\begin{Theorem}\label{TMainSpecht} {\rm \cite{BKW}} 
Let $\mu\in\Par_d$ and $\T\in\St(\mu)$. 
For each admissible $r$, the vectors $y_r v_\T$ and $\psi_r v_\T$ 
are homogeneous, 
and we have that
\begin{align*}
e(\bi) v_\T &=\de_{\bi,\bi^\T }v_\T\qquad\qquad\qquad\qquad\quad\,(\bi \in I^d),\\ 
\deg(y_r v_\T )&=\deg(y_r)+\deg(v_\T )\qquad\qquad\:\:(1 \leq r \leq d),\\
\deg(\psi_rv_\T )&=\deg(\psi_re(\bi^\T ))+\deg(v_\T ) \qquad(1 \leq r < d).
\end{align*} 
In particular, our grading makes $S(\mu)$ into a graded
$H_d^\La$-module. 
\end{Theorem}

The $q$-characters of the graded Specht modules are not difficult to describe: 

\begin{Corollary}\label{CCHS}
Let $\mu\in\RPar_d$. Then 
\begin{equation}\label{CHSMu}
\CH S(\mu)=\sum_{\T\in\St(\mu)}q^{\deg \T}\,\bi^\T.
\end{equation}
\end{Corollary}

It would be interesting to find different descriptions of the right-hand side of \ref{CHSMu}, which is of course a refinement of the number of standard $\mu$-tableaux.

The explicit grading of the Specht module $S(\mu)$, coming from Theorem~\ref{TMainSpecht}, and the grading of the irreducible module $D(\mu)$, obtained in section~\ref{SSGIrrMod} using general principles, are consistent in the following sense:

\begin{Theorem}
 {\cite[Theorem 5.9]{BKllt}} 
Let $\mu\in\RPar_d$. Then the radical $R(\mu)$ of the Specht module $S(\mu)$ is a homogeneous submodule of $S(\mu)$, and there is a (homogeneous) isomorphism of graded modules $S(\mu)/R(\mu)\cong D(\mu)$. 
\end{Theorem}

\subsection{Graded branching rule for Specht modules}\label{SSGBRSM}
Once the Specht modules have been graded, it is natural to consider the {\em graded branching rule} for them. For the classical branching rule for Specht modules over cyclotomic Hecke algebras, see \cite[Proposition 1.9]{AM}, which for symmetric groups goes back to James  \cite{Jbook}. 

It is very important that the grading shifts $d_{A_m}(\mu)$ occurring in this branching rule are the same combinatorial quantities which arise in the  definition of the quantum Fock space, see section~\ref{SSFock}. This is an indication of the fact that the grading on the module category should correspond to the quantization of its Grothendieck group. 

To state the branching rule, note that the natural embedding
\begin{equation}\label{natemb}
H_{d-1}^\La \hookrightarrow H_{d}^\La
\end{equation}
sending $X_r$ to $X_r$ for $1\leq r \leq d-1$ and $T_s$ to $T_s$ for $1\leq s \leq d-2$
maps
\begin{align*}
e(\bi) &\mapsto \sum_{i \in I}
e(i_1,\dots,i_{d-1},i)\qquad(\bi=(i_1,\dots,i_{d-1}) \in I^{d-1}),\\
y_r &\mapsto y_r
\qquad\qquad\qquad\:\;\quad\qquad (1 \leq r \leq d-1),\\
\psi_s &\mapsto \psi_s
\qquad\qquad\qquad\:\,\qquad\quad(1 \leq s \leq d-2).
\end{align*}
Hence this embedding is a degree-preserving homomorphism of graded algebras.
Recall the integers $d_A(\mu)$ defined in (\ref{EDMUA}) and the notation 
$M\langle m\rangle$  
for a graded module $M$ with grading shifted up by $m$.

\begin{Theorem}\label{TBr} {\rm \cite{BKW}} 
Let $\mu\in\Par_d$, and $A_1,\dots,A_b$ be all the removable nodes of $\mu$ in order from bottom to top. Then the restriction of $S(\mu)$ to $H_{d-1}^\La$
has a filtration 
$$
\{0\}=V_0\subset V_1\subset\dots\subset V_b= S(\mu) 
$$
as a graded $H_{d-1}^\La$-module
such that $$V_m/V_{m-1}\cong S(\mu_{A_m})\langle d_{A_m}(\mu)\rangle$$ for all $1\leq m\leq b$. 
\end{Theorem}

\subsection{Graded dimension of a block}\label{SSGDB}
Here we present a combinatorial formula for 
the graded dimension of any block $H^\La_\alpha$ of the cyclotomic Hecke algebra $H_d^\La$. In fact, the next theorem provides us with even more delicate information:

\begin{Theorem}\label{gdim} \cite[Theorem 4.20]{BKllt}
For $\alpha \in Q_+$ and 
$\bi, \bj \in I^\alpha$, we have that
\begin{align*}
\qdim\ e(\bi) H^\La_\alpha e(\bj) 
=
\sum_{\substack{
\mu \in \Par\\ \Stab, \T \in \St(\mu) \\
\bi^\Stab = \bi, \bi^\T = \bj}}
q^{\deg(\Stab)+\deg(\T)}
= 
\sum_{\substack{
\mu \in \Par\\ \Stab, \T \in \St(\mu) \\
\bi^\Stab = \bi, \bi^\T = \bj}}
q^{2\defect(\alpha)-\deg(\Stab)-\deg(\T)}.
\end{align*}
\end{Theorem}

We point out that the original proof of the theorem has been rather indirect. In particular, it relied on the categorification results of \cite{BKllt}. Very recently, Hu and Mathas \cite{HuMathas} gave a more direct proof by exhibiting a homogeneous basis of each $\ e(\bi) H^\La_\alpha e(\bj)$. This is explained in the next section. 

\subsection{Graded cellular structure on cyclotomic Hecke algebras}\label{SSGrCellStr}
Here we describe the recent work of Hu and Mathas \cite{HuMathas}. Let $\mu\in\Par_d$ and $1\leq m \leq d$. Denote by $A(m)\in\mu$ the node occupied with $m$ in $T^\mu$. Let $\mu_{\leq m}\in\Par_m$ be the union of the nodes $A(1),A(2),\dots,A(m)$. Set $\deg_m(\mu):=d_{A(m)}(\mu_{\leq m})$ so that $\deg(T^\mu)=\sum_{m=1}^d\deg_m(\mu)$. 
Define the element of degree $2\deg(T^\mu)$: 
$$
y^\mu:=\prod_{m=1}^d y_m^{\deg_m(\mu)}\in H_d^\La\qquad(\mu\in\Par_d).
$$

The combinatorial data for the Hu-Mathas cellular structure on $H_d^\La$ is the same as  for the Dipper-James-Mathas cellular structure on $H_d^\La$, namely  
${\mathscr P}=\Par_d$ and $\St(\mu)$ is the set of standard $\mu$-tableaux. On the other hand, let `$*$' be the anti-automorphism `$\circledast$' defined in (\ref{star}).
Finally, for $\Stab,\T\in\St(\mu)$, define the homogeneous element of degree $\deg\Stab+\deg\T$: 
\begin{equation}\label{GradedCellularBasis}
\psi^\mu_{\Stab,\T}:=\psi_{w_\Stab}e(\bi^\mu)y^\mu \psi_{w_\T}^\circledast.
\end{equation}

\begin{Theorem}\label{THuMathas} {\rm \cite{HuMathas}} 
$\{\psi^\mu_{\Stab,\T}\mid \mu\in\Par_d,\ \Stab,\T\in\St(\mu)\}$ is a homogeneous cellular basis of $H_d^\La$.
\end{Theorem} 

Since the Hu-Mathas cellular basis is homogeneous, the corresponding cell modules are automatically graded. It is checked in \cite{HuMathas} that these graded cell modules are precisely the graded Specht modules introduced in \cite{BKW} and described in section~\ref{SSHBSM}. The graded dimension formula in Theorem~\ref{gdim} follows immediately from Theorem~\ref{THuMathas}. For other important applications of Theorem~\ref{THuMathas} we refer the reader to \cite{HuMathas}.

\section{Graded induction, restriction, and  branching rules}\label{SGIndResBr}
The functors of induction and restriction and their refinements and generalizations play a crucial role in representation theory of Hecke algebras. In this section we review their definitions and key properties. We first consider induction and restriction for affine Khovanov-Lauda-Rouquier algebras and then pass to the cyclotomic Hecke algebras. We pay special attention to the graded aspect of the story. 

The functors of $i$-induction and $i$-restriction refine induction and restriction and `connect' blocks of cyclotomic Hecke algebras of different ranks. So it is natural to consider the cyclotomic Hecke algebras of all ranks together. This motivates the introduction of the algebra
\begin{equation}\label{EHLa}
H^\La:=\bigoplus_{d\geq 0} H^\La_d=\bigoplus_{\al\in Q_+} H^\La_\al,
\end{equation}
so that 
\begin{equation}\label{ERepHLa}
\Rep{H^\La}=\bigoplus_{d\geq 0} \Rep{H^\La_d}=\bigoplus_{\al\in Q_+}\Rep{H^\La_\al},
\end{equation}
and
\begin{equation}\label{EProjHLa}
\Proj{H^\La}=\bigoplus_{d\geq 0} \Proj{H^\La_d}=\bigoplus_{\al\in Q_+}\Proj{H^\La_\al}.
\end{equation}
We also have for the Grothendieck groups:
\begin{equation}\label{EGrRepHLa}
[\Rep{H^\La}]=\bigoplus_{d\geq 0} [\Rep{H^\La_d}]=\bigoplus_{\al\in Q_+}[\Rep{H^\La_\al}],
\end{equation}
and
\begin{equation}\label{EGrProjHLa}
[\Proj{H^\La}]=\bigoplus_{d\geq 0} [\Proj{H^\La_d}]=\bigoplus_{\al\in Q_+}[\Proj{H^\La_\al}].
\end{equation}

\subsection{Affine induction and restriction}\label{SSAIndRes}
Given $\alpha, \beta \in Q_+$, we set $$
H^\infty_{\alpha,\beta} := H^\infty_\alpha \otimes 
H^\infty_\beta,$$ 
viewed as an algebra in the usual way.
We denote 
the outer tensor product of an $H^\infty_\alpha$-module $M$ and an $H^\infty_\beta$-module 
$N$ by $M \boxtimes N$.

There is an obvious injective (non-unital) algebra homomorphism 
$$H^\infty_{\alpha,\beta}\,\into\, H^\infty_{\alpha+\beta}$$
mapping $e(\bi) \otimes e(\bj)$ to $e(\bi\bj)$,
where $\bi\bj$ denotes the concatenation of the two sequences. The image of the identity
element of $H^\infty_{\alpha,\beta}$ under this map is
\begin{equation}\label{EEAlBe}
e_{\alpha,\beta}:= \sum_{\bi \in I^\alpha,\ \bj \in I^\beta} e(\bi\bj).
\end{equation}
Let $\Ind_{\alpha,\beta}^{\alpha+\beta}$ and $\Res_{\alpha,\beta}^{\alpha+\beta}$
denote the corresponding induction and restriction functors between the graded module categories: 
\begin{align}
\Ind_{\alpha,\beta}^{\alpha+\beta} &:= H^\infty_{\alpha+\beta} e_{\alpha,\beta}
\otimes_{H^\infty_{\alpha,\beta}} ?:\Mod{H^\infty_{\alpha,\beta}} \rightarrow \Mod{H^\infty_{\alpha+\beta}},\\
\Res_{\alpha,\beta}^{\alpha+\beta} &:= e_{\alpha,\beta} H^\infty_{\alpha+\beta}
\otimes_{H^\infty_{\alpha+\beta}} ?:\Mod{H^\infty_{\alpha+\beta}}\rightarrow \Mod{H^\infty_{\alpha,\beta}}.
\end{align}
The functor $\Res_{\alpha,\beta}^{\alpha+\beta}$ is just left multiplication by
the idempotent $e_{\alpha,\beta}$, so it is exact and sends finite dimensional modules to
finite dimensional modules. 
By \cite[Proposition 2.16]{KL1},
$e_{\alpha,\beta} H^\infty_{\alpha+\beta}$ is a graded free left $H^\infty_{\alpha,\beta}$-module of finite rank,
so $\Res_{\alpha,\beta}^{\alpha+\beta}$ also sends projectives to projectives.

The functor $\Ind_{\alpha,\beta}^{\alpha+\beta}$ is left adjoint to $\Res_{\alpha,\beta}^{\alpha+\beta}$,
so it 
sends projectives to projectives.
Finally $H^\infty_{\alpha+\beta} e_{\alpha,\beta}$ is a graded 
free right $H^\infty_{\alpha, \beta}$-module of finite rank, so
$\Ind_{\alpha,\beta}^{\alpha+\beta}$ sends finite dimensional modules to finite dimensional
modules.

\subsection{\boldmath Affine $i$-induction and $i$-restriction}\label{SSAIIndIRes}
For $i \in I$, let $P(i)$ denote the regular
representation of $H^\infty_{\alpha_i}$.
Define a functor
\begin{equation}\label{fri1}
\theta_i := \Ind_{\alpha,\alpha_i}^{\alpha+\alpha_i} (? \boxtimes P(i)):
\Mod{H^\infty_\alpha} \rightarrow \Mod{H^\infty_{\alpha+\alpha_i}}.
\end{equation}
This functor is exact, and it restricts to a functor
$\theta_i: \Proj{H^\infty_\alpha}
\rightarrow \Proj{H^\infty_{\alpha+\alpha_i}}$.

The functor $\theta_i$ possesses a right adjoint 
\begin{equation}\label{fri2}
\theta_i^*:= \HOM_{(H^\infty_{\alpha_i})'}(P(i), ?):\Mod{H^\infty_{\alpha+\alpha_i}}\rightarrow \Mod{H^\infty_\alpha},
\end{equation}
where $(H^\infty_{\alpha_i})'$ denotes the subalgebra $1 \otimes H^\infty_{\alpha_i}$
of $H^\infty_{\alpha,\alpha_i}$.
Equivalently, $\theta_i^*$ is defined by
multiplication by the idempotent $e_{\alpha,\alpha_i}$
followed by restriction to the subalgebra $H^\infty_\alpha = H^\infty_{\alpha} \otimes 1$
of $H^\infty_{\alpha,\alpha_i}$.
The functor $\theta_i^*$ also restricts to a well-defined functor
$
\theta^*_i: \Rep{H^\infty_{\alpha+\alpha_i}}
\rightarrow \Rep{H^\infty_{\alpha}}$.

\subsection{Affine divided powers}\label{SSADP}
In the case $\alpha = n \alpha_i$ for some $i \in I$,
the algebra
$H^\infty_\alpha$ is isomorphic to the nil-Hecke algebra, see  \cite[$\S$2.2]{KL1}. So it has 
a canonical representation on the polynomial algebra $F[y_1,\dots,y_n]$
such that each $y_r$ acts as multiplication by $y_r$ and each $\psi_r$ acts
as the {\em divided difference operator}
$$
\partial_r: f \mapsto \frac{{^{s_r}} f - f}{y_{r}-y_{r+1}}.
$$
Let $P(i^{(n)})$ denote the polynomial representation of $H^\infty_{n\alpha_i}$
viewed as a 
graded $H^\infty_{n\alpha_i}$-module with grading defined by
$$
\deg(y_1^{m_1} \cdots y_n^{m_n}) := 2m_1+\cdots+2m_n - \frac{1}{2}n(n-1).
$$
Denoting the left 
regular $H^\infty_{n\alpha_i}$-module
by $P(i^n)$, it is noted in \cite[$\S$2.2]{KL1} that
\begin{equation}\label{div}
P(i^n) \cong [n]! \cdot P(i^{(n)}),
\end{equation}
where the notation is as in (\ref{EGI}) and (\ref{EAMod}). 
In particular, $P(i^{(n)})$ is projective.

Now we 
generalize the definition of the functors $\theta_i$ and $\theta_i^*$:
for $i \in I$ and $n \geq 1$, set
\begin{align}\label{div1}
\theta_{i}^{(n)}&:= 
\Ind_{\alpha,n \alpha_i}^{\alpha+n\alpha_i} (? \boxtimes P(i^{(n)})):\Mod{H^\infty_\alpha} \rightarrow \Mod{H^\infty_{\alpha+n\alpha_i}},\\
(\theta_{i}^*)^{(n)}&:= 
\HOM_{(H^\infty_{n \alpha_i})'}(P(i^{(n)}), ?): \Mod{H^\infty_{\alpha+n\alpha_i}} \rightarrow \Mod{H^\infty_{\alpha}},\label{div2}
\end{align}
where $(H^\infty_{n\alpha_i})' := 1 \otimes H^\infty_{n\alpha_i} \subseteq H^\infty_{\alpha,n\alpha_i}$.
Both functors are exact, and induce $\Laurent$-module homomorphisms 
$\theta_i^{(n)}$ and $(\theta_i^*)^{(n)}$
on the various Grothendieck groups.
By the definitions (\ref{fri1})--(\ref{fri2}) and
transitivity of induction and restriction,
there are natural isomorphisms
$$
\theta_i^n
\cong \Ind_{\alpha,n\alpha_i}^{\alpha+n\alpha_i}(? \boxtimes P(i^n)),
\qquad(\theta_i^*)^n
\cong \HOM_{(H^\infty_{n\alpha_i})'}(P(i^n), ?).
$$
Hence (\ref{div}) implies that
the $n$th powers
of $\theta_i$ and $\theta_i^*$
decompose as
\begin{equation}\label{dividedform}
\theta_i^n \cong [n]! \cdot \theta_i^{(n)},
\qquad
(\theta_i^*)^n \cong [n]! \cdot (\theta_i^*)^{(n)}.
\end{equation}

\subsection{\boldmath Cyclotomic $i$-induction and $i$-restriction}\label{sir}
To deal with the cyclotomic $i$-induction and $i$-restriction, fix $\La\in P_+$. 
For any $i \in I$, the 
embedding 
$$H_\alpha^\infty = H_\alpha^\infty \otimes 1 
\hookrightarrow H^\infty_{\alpha,\alpha_i}\hookrightarrow H^\infty_{\alpha+\alpha_i}$$
factors through the quotients to induce a (not necessarily injective)
graded algebra homomorphism
\begin{equation}\label{iota1}
\iota_{\alpha,\alpha_i}:H_\alpha^\La \rightarrow H_{\alpha+\alpha_i}^\La.
\end{equation}
This homomorphism maps the identity element of $H_\alpha^\La$ to the idempotent
$e_{\alpha,\alpha_i} \in H_{\alpha+\alpha_i}^\La$.
Alternatively, one can define the idempotent $e_{\al,\al_i}\in H_{\al+\be}^\La$  using a special case of the formula (\ref{EEAlBe}). 

For any $i \in I$ and $\alpha \in Q_+$, let
$E_i$ and $F_i$ be the functors 
\begin{align}\label{gre}
E_i:=e_{\alpha,\alpha_i} H^\La_{\alpha+\alpha_i} \otimes_{H^\La_{\alpha+\alpha_i}} ?
:\Rep{H^\La_{\alpha+\alpha_i}} \rightarrow \Rep{H^\La_{\alpha}},\\
F_i:=
H^\La_{\alpha+\alpha_i} e_{\alpha,\alpha_i} \otimes_{H^\La_{\alpha}} ? \langle 1\!-\!(\La\!-\!\alpha,\alpha_i) \rangle
:\Rep{H^\La_\alpha} \rightarrow \Rep{H^\La_{\alpha+\alpha_i}},\label{grf}
\end{align}
interpreting the tensor products via (\ref{iota1}). 
By taking direct sums over all $\alpha \in Q_+$ we  consider $E_i$ and $F_i$ as functors on the category 
$
\Rep{H^\La},
$
defined in (\ref{ERepHLa}). 

It could have been more appropriate to use a more detailed notation $E_i^\La$ and $F_i^\La$, but since $\La$ is usually fixed, using short hand versions will not cause a problem. 

Note the grading shift $\langle 1\!-\!(\La\!-\!\alpha,\alpha_i) \rangle$, appearing in the definition of $F_i$. At this stage this seems contrived. The first indication that this shift is convenient and natural comes from Theorem~\ref{gg}. Another reason will become clear when we deal with categorifications.  

Also define a functor $K_i$ by letting
\begin{equation}
K_i:\Rep{H^\La_\alpha} \rightarrow \Rep{H^\La_\alpha}
\end{equation}
denote the degree shift functor 
$M\mapsto M \langle (\La-\alpha,\alpha_i) \rangle$. Taking the direct sums over all $\alpha \in Q_+$ we then get a functor on $
\Rep{H^\La}$.

\begin{Lemma}\label{cad}
The functors $E_i$ and $F_i$ are both
exact and send projectives to projectives. Furthermore, 
there is a canonical adjunction making
$$(F_i K_i \langle-1\rangle, E_i)$$ into an adjoint pair.
\end{Lemma}

There is an equivalent way to describe the
functors $E_i$ and $F_i$ which relates them to
the functors $\theta_i^*$ and $\theta_i$ from (\ref{fri1})--(\ref{fri2}).
To formulate this, we first introduce the
{\em inflation} and {\em truncation} functors
\begin{align*}
\infl&:\Rep{H^\La_\alpha} \rightarrow \Rep{H_\alpha^\infty},\\
\pr&:\Rep{H^\infty_\al} \rightarrow \Rep{H^\La_\alpha},
\end{align*}
where  for $M \in \Rep{H^\La_\alpha}$, we denote by $\infl\, M$ its pull-back through the natural
surjection $$H_\alpha^\infty \twoheadrightarrow H^\La_\alpha,$$
and for $N \in \Rep{H_\alpha^\infty}$ we write
$\pr\, N$ for $H^\La_\alpha \otimes_{H_\alpha^\infty} N$, which is the
largest graded quotient of $N$ that factors through to $H^\La_\alpha$-module. Like many other functors defined in this section, $\pr$ and $\infl$ of course depend on $\La$, but we again suppress $\La$ from the notation. 

We obviously have that
\begin{equation}\label{ididid}
\pr \circ \infl  = \id.
\end{equation}
Observe also that $(\pr, \infl)$ is an adjoint pair
in a canonical way.
Hence, $\pr$ sends projectives to projectives, 
so
it restricts to give an additive functor
\begin{equation}\label{prprpr}
\pr:\Proj{H^\infty_\alpha} \rightarrow \Proj{H^\La_\alpha}.
\end{equation}

Now we give another description of the functors $E_i$ and $F_i$:

\begin{Lemma}\label{sun} \cite{BKllt}
There are canonical isomorphisms of functors
$$E_i \cong \pr \circ \theta_i^* \circ \infl
\qquad\text{and}\qquad F_i K_i \langle -1 \rangle 
\cong \pr \circ \theta_i \circ \infl.$$
\end{Lemma}

\subsection{Cyclotomic divided powers}\label{sdivp}
Lemma~\ref{sun} also
makes it clear how to define
divided powers $E_i^{(n)}$ and $F_i^{(n)}$ of the functors
$E_i$ and $F_i$. For $n \geq 1$, set
\begin{align*}
E_i^{(n)} &:= \pr \circ (\theta_i^*)^{(n)} \circ \infl:
\Rep{H^\La_{\alpha+n\alpha_i}} \rightarrow \Rep{H^\La_\alpha},\\
F_i^{(n)} &:= \pr \circ \theta_i^{(n)} \circ \infl 
\langle 
 n^2 - n(\La-\alpha,\alpha_i) \rangle
:\Rep{H^\La_\alpha} \rightarrow \Rep{H^\La_{\alpha+n\alpha_i}},
\end{align*}
recalling (\ref{div1})--(\ref{div2}).
Again we use the same notation $E_i^{(n)}$ and $F_i^{(n)}$ for the 
direct sums of these functors over all $\alpha \in Q_+$.

\begin{Lemma}\label{divp}
There are isomorphisms 
$$E_i^n \cong [n]! \cdot E_i^{(n)}\qquad\text{and}\qquad
F_i^n \cong [n]! \cdot F_i^{(n)}.
$$
Hence
$E_i^{(n)}$ and $F_i^{(n)}$ are exact
and send projectives to projectives.
\end{Lemma}

To summarize, we have now defined the exact functors
\begin{equation}\label{EEFKRep}
E_i^{(n)},F_i^{(n)},K_i:\Rep{H^\La}\to\Rep{H^\La}\qquad(i\in I,n\in\Z_{>0}),
\end{equation}
and the additive functors
\begin{equation}\label{EEFKProj}
E_i^{(n)},F_i^{(n)},K_i:\Proj{H^\La}\to\Proj{H^\La}\qquad(i\in I,n\in\Z_{>0}),
\end{equation}

\subsection{Graded branching rules for irreducible modules}\label{SSBrBrRules}
As usual we work with a fixed $\La\in P_+$ of level $l$. The original branching rules for irreducible modules over symmetric groups  in characteristic $p$ were developed in \cite{KBrI}--\cite{KBrIV} and \cite{BKtr}. This case corresponds to $l=1$ and $\xi=1$. In characteristic zero this goes back to a 1908 paper of Schur \cite[p. 253]{Schur}. 

These branching rules were generalized in various directions, see for example \cite{BBr,Gr,GV,Abranch,DuRui, BKS,BOBr,BKHCl,Ku,Hu,Ts,WW,Sh,CR,AJL,EM}. The graded case is dealt with in \cite{BKllt}. Recall the reduced signatures, good nodes, and other related notions introduced in section~\ref{SSCrGraph}. We point out that all of these combinatorial notions depend on $\kappa$. 

Parts (i),(ii),(iv) of the next theorem go back to \cite{KBrIII,KBrIV}, and part (iii) to \cite{Gr}. 

\begin{Theorem}\label{gg} \cite[Theorem 4.12]{BKllt}
For any $\mu \in \RPar_\al$ and $i \in I$, we have:
\begin{itemize}
\item[(i)] $E_i D(\mu)$ is non-zero if and only if
$\eps_i(\mu) \neq 0$, in which case
$E_i D(\mu)$ has irreducible 
socle  isomorphic
to $D(\tilde e_i \mu)\langle \eps_i(\mu)-1\rangle$
and head isomorphic to
$D(\tilde e_i \mu)\langle 1-\eps_i(\mu)\rangle$.
\item[(ii)]
 $F_i D(\mu)$ is non-zero if and only if
$\phi_i(\mu) \neq 0$, in which case
$F_i D(\mu)$ has irreducible 
socle isomorphic to
 $D(\tilde f_i \mu) \langle \phi_i(\mu)-1\rangle$
and
head isomorphic to
 $D(\tilde f_i \mu) \langle 1-\phi_i(\mu)\rangle$.
\item[(iii)] In the Grothendieck group we have that
\begin{align*}
[E_i D(\mu)] 
&= 
[\eps_i(\mu)] \cdot [D(\tilde e_i \mu)]
+ \sum_{\substack{\nu \in \RPar_{\al-\al_i}\\
\eps_i(\nu) < \eps_i(\mu)-1}} u_{\nu,\mu;i}(q)\cdot
[D(\nu)],\\
[F_i D(\mu)] 
&= 
[\phi_i(\mu)] \cdot [D(\tilde f_i \mu)]
+ \sum_{\substack{\nu \in \RPar_{\al+\al_i} \\
\phi_i(\nu) < \phi_i(\mu)-1}} v_{\nu,\mu;i}(q)\cdot 
[D(\nu)],
\end{align*}
for some bar-invariant polynomials $u_{\nu,\mu;i}(q), v_{\nu,\mu;i}(q) \in 
\Z[q,q^{-1}]$ with non-negative coefficients.
(The first term on the right-hand side of these formulae
should be interpreted as zero if $\eps_i(\mu) = 0$ (resp.\
$\phi_i(\mu) = 0$).)
\item[(iv)] Viewing $F[x]$ as a graded algebra
by putting $x$ in degree 2, 
there are graded algebra isomorphisms
\begin{align*}
F[x] / (x^{\eps_i(\mu)}) &\cong
\END_{H^\La_{\al-\al_i}}(E_i D(\mu)),\\
F[x] / (x^{\phi_i(\mu)}) &\cong
\END_{H^\La_{\al+\al_i}}(F_i D(\mu)).
\end{align*}
\end{itemize}
\end{Theorem}

There are more results on branching rules known in level $1$ which we will comment on in section~\ref{SSBr}.


\section{Quantum groups}\label{SQG}
As mentioned in the introduction, graded representations of the cyclotomic Hecke algebras $H_d^\La$ categorify the irreducible module $V(\La)$ over the quantized enveloping algebra  $U_q(\g)$. In this section we explain some facts from the theory of quantum groups which are needed to state the categorification results precisely. 

\subsection{\boldmath The algebra $\mathbf f$}\label{SSAF}
It is convenient to first introduce the Lusztig's algebra $\mathbf f$ from \cite[$\S$1.2]{Lubook}
attached to the Cartan matrix (\ref{ECM}) over the field $\Q(q)$ (our $q$ is Lusztig's $v^{-1}$).
Thus $\mathbf f$ 
is the $\Q(q)$-algebra on generators $\theta_i\:(i \in I)$ subject only to
the quantum Serre relations
\begin{equation}\label{qserre}
(\ad_q \theta_i)^{1-a_{j,i}}(\theta_j)=0
\end{equation}
where  
\begin{equation}\label{adq}
(\ad_q x)^n (y):=\sum_{m=0}^n(-1)^m
\left[
\begin{matrix}
 n   \\
 m
\end{matrix}
\right]
x^{n-m}yx^m.
\end{equation}
There is a $Q_+$-grading 
$\mathbf f = \bigoplus_{\alpha \in Q_+} \mathbf f_\alpha$
such that $\theta_i$ is of degree $\alpha_i$.
The algebra $\mathbf f$ possesses a bar-involution
$-:\mathbf f \rightarrow \mathbf f$
that is anti-linear with respect to
the field automorphism sending $q$ to $q^{-1}$, such that
$\overline{\theta_i} = \theta_i$ for each $i \in I$.

If we equip $\mathbf f \otimes \mathbf f$ with algebra structure via the rule
$$
(x_1 \otimes x_2) (y_1 \otimes y_2)
= q^{-(\alpha,\beta)} x_1 y_1 \otimes x_2 y_2
$$
for $x_2 \in \mathbf f_{\alpha}$ and $y_1 \in \mathbf f_{\beta}$,
there is a $Q_+$-graded
comultiplication 
$m^*:\mathbf f  \rightarrow \mathbf f \otimes \mathbf f,$
which is
the unique algebra homomorphism
such that 
$$m^*: \theta_i \mapsto \theta_i \otimes 1 + 1 \otimes \theta_i\qquad (i \in I).$$
For $\alpha,\beta \in Q_+$, we let
$$
m_{\alpha,\beta}: \mathbf f_\alpha 
\otimes \mathbf f_\beta 
\rightarrow \mathbf f_{\alpha+\beta},\qquad
m^*_{\alpha,\beta}:\mathbf f_{\alpha+\beta} \rightarrow \mathbf f_\alpha 
\otimes \mathbf f_{\beta}
$$
denote the multiplication and comultiplication maps induced on individual 
weight spaces,
so $m = \sum m_{\alpha,\beta}$ is the multiplication on $\mathbf f$
and $m^* = \sum m_{\alpha,\beta}^*$.

Finally let ${_\Laurent} \mathbf f$
be the $\Laurent$-subalgebra of $\mathbf f$ generated by the quantum
divided powers $\theta_i^{(n)} := \theta_i^n / [n]!$.
The bar
involution induces an involution of ${_\Laurent}\mathbf f$, and also
the map $m^*$ restricts to a well-defined comultiplication
$m^*:{_\Laurent}\mathbf f \rightarrow {_\Laurent}\mathbf f \otimes {_\Laurent}\mathbf f$.

Let $\mathbf B = \bigsqcup_{\alpha \in Q_+} \mathbf B_\alpha$ 
be the canonical basis for $\mathbf f = \bigoplus_{\alpha \in Q_+}
\mathbf f_\alpha$, see 
\cite[$\S$14.4]{Lubook}.

\subsection{\boldmath The quantized enveloping algebra $U_q(\g)$}\label{SSQEAUG}
Let $\g$ be the Kac-Moody algebra corresponding to the Cartan matrix
(\ref{ECM}), so
$\g = \widehat{\mathfrak{sl}}_e(\C)$ if $e > 0$ 
and $\g = 
\mathfrak{sl}_\infty(\C)$ if $e = 0$.

Let $U_q(\g)$ be the quantized enveloping algebra of $\g$. So 
$U_q(\g)$
is the $\Q(q)$-algebra generated by the {\em Chevalley  
generators} $E_i,F_i,K_i^{\pm 1}$ for $i\in I$, subject 
only to the usual quantum Serre relations (for all admissible $i,j\in I$):
\begin{align}
K_iK_j&=K_jK_i,\quad \hspace{8.5mm}K_iK_i^{-1}=1,\\
K_iE_jK_i^{-1}&= q^{a_{i,j}}E_j,\quad \hspace{3mm}K_iF_jK_i^{-1}= q^{-a_{i,j}}F_j,\label{wight}\\
[E_i,F_j]&=\de_{i,j}\frac{K_i-K_i^{-1}}{q-q^{-1}},\label{mix}\\
(\ad_q E_i)^{1-a_{j,i}}(E_j)&=0\qquad\qquad\qquad(i\neq j),\label{serree}\\
\quad (\ad_q F_i)^{1-a_{j,i}}(F_j)&=0\qquad\qquad\qquad(i\neq j),\label{serref}
\end{align}
with $(\ad_q x)^n (y)$ as defined in (\ref{adq}). 
Let $U_q(\g)^-$ be the subalgebra of $U_q(\g)$ generated
by the $F_i$'s.

We consider $U_q(\g)$ as a Hopf algebra with respect to the coproduct given for all $i\in I$ as follows: 
$$
\Delta:\ K_i\mapsto K_i\otimes K_i,\quad E_i\mapsto E_i\otimes K_i+1\otimes E_i,\quad F_i\mapsto F_i\otimes 1+K_i^{-1}\otimes F_i.
$$

The {\em bar-involution} $-:U_q(\g) \rightarrow U_q(\g)$
is the anti-linear involution such that
$$
\overline{K_i} = K_i^{-1},\qquad
\overline{E_i} = E_i, \qquad \overline{F_i} = F_i.
$$
Given a $U_q(\g)$-module $V$, a
{\em compatible bar-involution} on $V$ means
an anti-linear involution
$-:V \rightarrow V$ such that $\overline{xv} = \overline{x}\, \overline{v}$
for all $x \in U_q(\g)$ and $v \in V$.

Also let $\tau:U_q(\g) \rightarrow U_q(\g)$ be the 
anti-linear anti-automorphism defined by
\begin{equation}\label{taudef}
\tau:\
 K_i \mapsto K_i^{-1}, \qquad E_i \mapsto q F_i K_i^{-1},
\qquad F_i \mapsto q^{-1} K_i E_i.
\end{equation}

Let $U_q(\g)_\Laurent$ denote Lusztig's $\Laurent$-form for $U_q(\g)$, which is
the $\Laurent$-subalgebra generated by the quantum divided
powers $E_i^{(n)} := E_i^n / [n]!$ and $F_i^{(n)} := F_i^n / [n]!$
for all $i \in I$ and $n \geq 1$.
The bar-involution, the comultiplication $\Delta$ and the anti-automorphism
$\tau$ descend to this $\Laurent$-form.

Comparing the relations of $U_q(\g)$ and $\mathbf f$, 
it follows easily that there
is an algebra homomorphism 
\begin{equation}\label{flatmap}
\mathbf f \rightarrow U_q(\g),
\qquad
x \mapsto x^\flat
\end{equation}
such that $\theta_i^\flat := q^{-1} F_i K_i=\tau^{-1}(E_i)$, 
and an algebra isomorphism 
\begin{equation}
\mathbf f \stackrel{\sim}{\rightarrow} U_q(\g)^-, \qquad x \mapsto x^-
\end{equation}
such that $\theta_i^- := F_i$.

\subsection{\boldmath The module $V(\La)$}\label{md}
Let $V(\La)$ denote the integrable highest weight
module for $U_q(\g)$ of highest weight $\La$, where $\La$
is the dominant integral weight fixed in (\ref{ELa}).
Fix also a choice of a non-zero highest weight vector $v_\La \in V(\La)$.
The module $V(\La)$ possesses a unique
compatible bar-involution
$-:V(\La) \rightarrow V(\La)$
such that $\overline{v_\La} = v_\La$.

The {\em contravariant form} $(.,.)$ on $V(\La)$ is the
unique symmetric bilinear form such that 
\begin{itemize}
\item[(1)]
$(E_i v, w) = (v, F_i w)$ and $(F_i v, w) = (v, E_i w)$
for all $v, w \in V(\La)$ and $i \in I$;
\item[(2)] 
$(v_\La, v_\La) = 1$.
\end{itemize}

The {\em Shapovalov form}
$\langle.,.\rangle$ on $V(\La)$ 
 is the unique 
sesquilinear form  (anti-linear in the first argument, linear in the second)
on $V(\La)$
such that 
\begin{itemize}
\item[(1)] $\langle u v, w\rangle = \langle v, \tau(u) w\rangle$ for
all $u \in U_q(\g)$ and $v, w \in V(\La)$;
\item[(2)] $\langle v_\La, v_\La\rangle = 1$.
\end{itemize}

Let $V(\La)_\Laurent$ denote the {\em standard $\Laurent$-form} 
for $V(\La)$, that is,
the $U_q(\g)_\Laurent$-submodule of $V(\La)$ generated by the 
highest weight vector $v_\La$. 
Let $V(\La)_\Laurent^*$ denote the {\em costandard $\Laurent$-form} for $V(\La)$,
that is, the dual lattice
\begin{align*}
V(\La)_\Laurent^* &= \{v \in V(\La)\:|\:( v,w )
\in \Laurent\text{ for all }
w\in V(\La)_\Laurent\}\\
&= \{v \in V(\La)\:|\:\langle v,w \rangle
\in \Laurent\text{ for all }
w\in V(\La)_\Laurent\}.
\end{align*}

\subsection{Fock spaces}\label{SSFock} We review the Fock space theory following \cite{BKllt}. 
Recall the tuple $\kappa=(k_1,\dots,k_l)$, with $\La=\La(\kappa)$, fixed in (\ref{EKappa}), and the definitions of $d_A(\mu)$, $d^B(\mu)$, and $d_i(\mu)$ from section~\ref{SSDeg}. Define the {\em Fock space} $F(\kappa)$ to be the $\Q(q)$-vector space on basis
$$\{M_\mu\:|\:\mu \in \Par\},$$
referred to as the {\em monomial basis}, with $U_q(\g)$-action defined by
\begin{align}\label{act1}
E_i {M}_\mu&:=\sum_A q^{d_A(\mu)}{M}_{\mu_A},\\
\label{act2}
F_i {M}_\mu&:=\sum_B q^{-d^B(\mu)}{M}_{\mu^B},\\
K_i {M}_\mu&:=q^{d_i(\mu)}{M}_\mu,\label{act3}
\end{align}
where the first sum is over all removable $i$-nodes $A$ for $\mu$, and the 
second sum is over all addable $i$-nodes $B$ for $\mu$. 

When $l=1$, this construction originates in work of Hayashi \cite{Hayashi}
and Misra and Miwa \cite{MiM}.
When $l > 1$, $F(\kappa)$ was first studied in \cite{JMMO}.
We note that the Fock space $F(\kappa)$ is simply the tensor product
of $l$ level one Fock spaces:
\begin{equation}\label{tp}
F(\kappa) = F(\La_{k_1}) \otimes \cdots \otimes F(\La_{k_l}),
\end{equation}
on which the $U_q(\g)$-structure is defined via
the comultiplication $\De$,
so that
$M_\mu$ is identified with
$M_{\mu^{(1)}} \otimes \cdots \otimes M_{\mu^{(l)}}$
for each $\mu=(\mu^{(1)},\dots,\mu^{(l)}) \in \Par$.

\begin{Theorem}\label{TBarFock} {\rm \cite[Theorem 3.26]{BKllt}} 
There is a compatible bar-involution on $F(\kappa)$ such that
$$\overline{M_\mu} = M_\mu + 
(\text{an $\Laurent$-linear combination of $M_\nu$'s
for $\nu \lex \mu$}).$$
\end{Theorem}

The Fock space in general is not irreducible as a $U_q(\g)$-module, so there could be more than one compatible bar-involution on it, but we will always work with the fixed one constructed in \cite[Theorem 3.26]{BKllt}. That construction ultimately depends on {\em Uglov's Fock spaces} \cite{Uglov}, which do have a canonical bar-involution, as well as   stability results of Yvonne \cite{Y2}, which allow us to connect Uglov's Fock spaces with the Fock space $F(\kappa)$ considered here. The idea goes back to Ariki \cite{Aclass}.

The vector $M_{\varnothing}$ is a highest weight vector
of weight $\La$.
Moreover, the $\La$-weight space of $F(\kappa)$
is one dimensional. By complete reducibility, 
it follows that there is a canonical $U_q(\g)$-module
homomorphism
\begin{equation}\label{pidef}
\pi_\kappa:F(\kappa) \twoheadrightarrow V(\La),
\qquad
M_\varnothing \mapsto v_\La.
\end{equation}
The map $\pi_\kappa$ intertwines the bar-involution
on $F(\kappa)$ with the one defined earlier on $V(\La)$.

For any $\mu \in \Par$, we define
\begin{equation}\label{stmn}
S_\mu := \pi_\kappa(M_\mu),
\end{equation}
and call this a {\em standard monomial} in $V(\La)$.We note that in general the standard monomials in $V(\La)$ are not linearly independent. 
Applying $\pi_\kappa$ to (\ref{act1}) and (\ref{act2}), we get that
\begin{align}\label{sact1}
E_i S_\mu=\sum_A q^{d_A(\mu)}S_{\mu_A},&\qquad
F_i S_\mu=\sum_B q^{-d^B(\mu)}S_{\mu^B},
\end{align}
where the first sum is over all removable $i$-nodes $A$ for $\mu$, and the second sum is over all addable $i$-nodes $B$ for $\mu$. 

Also define the {\em Shapovalov form} on $F(\kappa)$ as the sesquilinear form $\langle.,.\rangle$ on $F(\kappa)$ such that 
\begin{equation}\label{sess}
\langle M_\mu, \overline{M_\nu}\rangle = \delta_{\mu,\nu}
\end{equation}
for all $\mu, \nu \in \Par$.
By \cite[(3.41)]{BKllt}, we then have 
\begin{equation}\label{sessq}
\langle xu, v \rangle = \langle u, \tau(x) v\rangle\qquad(x \in U_q(\g),\ u, v \in F(\kappa)).
\end{equation}
The map $\pi_\kappa$ intertwines the Shapovalov form 
on $F(\kappa)$ with the Shapovalov form defined earlier on $V(\La)$.

Let $F(\kappa)_\Laurent$ denote the free $\Laurent$-submodule
of $F(\kappa)$ spanned by the $M_\mu$'s. It 
is invariant under the action of $U_q(\g)_\Laurent$. All the above definitions descend to this integral form. We also have that 
$$
V(\La)_\Laurent^* = \pi_\kappa\left(F(\kappa)_\Laurent\right).
$$

\subsection{Canonical bases}\label{ssb}
In this section, we introduce two new bases for $F(\kappa)_\Laurent$: 
the {\em dual canonical basis} 
$$\{L_\mu\:|\:\mu \in \Par\},$$
and the {\em quasi-canonical basis}
$$\{P_\mu\:|\:\mu \in \Par\}.$$
For $\mu\in\Par$, we define $L_\mu$ to be the unique
bar-invariant vector in $F(\kappa)_\Laurent$ such that
\begin{align}
L_\mu &= M_\mu + \text{(a $q\Z[q]$-linear combination of
$M_\nu$'s for $\nu \lex \mu$}).
\end{align}
The existence and uniqueness of these vectors follows from
Lusztig's lemma \cite[Lemma 24.2.1]{Lubook} and Theorem~\ref{TBarFock}.
We note that 
$\{L_\mu\:|\:\mu \in \Par\}$
is an upper global crystal basis in the sense of Kashiwara. The lower global crystal basis (=Lusztig's canonical basis) will not be used here. We refer the reader to \cite{BKllt} for an explanation of its role. 

Let us introduce notation for the transition matrices between
the monomial and dual-canonical bases: 
\begin{equation}\label{ml}
M_\mu = \sum_{\nu \in \Par} d_{\mu,\nu}(q) L_\nu,
\qquad
L_\mu = \sum_{\nu \in \Par} p_{\mu,\nu}(-q) M_\nu,
\end{equation}
for polynomials 
$$d_{\mu,\nu}(q), p_{\mu,\nu}(q) \in \Z[q].$$
Note $d_{\mu,\mu}(q) = p_{\mu,\mu}(1) = 1$ and $d_{\mu,\nu}(q)
= p_{\mu,\nu}(q) = 0$ unless $\nu \lexeq \mu$.

Finally introduce the quai-canonical basis $\{P_\mu\:|\:\mu \in \Par\}$ for
$F(\kappa)$ by transposing these transition matrices:
\begin{equation}\label{pm}
P_\mu = \sum_{\nu \in \Par} d_{\nu,\mu}(q) M_\nu,
\qquad
M_\mu = \sum_{\nu \in \Par} p_{\nu,\mu}(-q) P_\nu.
\end{equation}
The elements of the quasi-canonical basis are in general {\em not} invariant under the bar-involution.


Recall the set $\RPar$ of ($\kappa$-)restricted multipartitions, the operations $\tilde e_i, \tilde f_i$, and the functions $\eps_i, \phi_i, \wt$, defined in  section~\ref{SSCrGraph}. Set
\begin{equation}
D_\mu := \pi(L_\mu)\qquad(\mu \in \RPar).
\end{equation}

The following result collects the necessary information on the dual canonical basis in $V(\La)$. 

\begin{Theorem}\label{tdcb} \cite[sections 3.7--3.9]{BKllt}
Let $\kappa$ and $\La=\La(\kappa)$ be as in (\ref{EKappa}) and (\ref{ELa}).
\begin{itemize}
\item[(i)] The vectors
$$\{D_\mu\:|\:\mu \in \RPar\}$$
give a basis for $V(\La)_\Laurent^*$
which coincides with Lusztig's dual-canonical basis of $V(\La)$. 
In particular, each $D_\mu$ is bar-invariant.

\item[(ii)] The crystal
$$(\RPar,
\tilde e_i, \tilde f_i, \eps_i, \phi_i, \wt)$$
is the highest weight crystal associated to
$V(\La)$.
 
\item[(iii)] The vectors
$$\{L_\mu\:|\:\mu \in \Par \setminus \RPar\}$$
give a basis for $\ker \pi$ as a free $\Laurent$-module.

\item[(iv)]
For $\mu \in \RPar_\alpha$ and $i \in I$ we have that
\begin{align*}
E_i D_\mu &=
[\eps_i(\mu)] D_{\tilde e_i \mu}
+ \sum_{\substack{\nu \in \RPar_{\al-\al_i} \\
\eps_i(\nu) < \eps_i(\mu)-1
}}
x_{\mu,\nu;i}(q) D_\nu,\\
F_i D_\mu &=
[\phi_i(\mu)] D_{\tilde f_i \mu}
+ \sum_{\substack{\nu \in \RPar_{\al+\al_i} \\
\phi_i(\nu) < \phi_i(\mu)-1
}}
y_{\mu,\nu;i}(q) D_\nu,
\end{align*}
for some bar-invariant
$$x_{\mu,\nu;i}(q) \in q^{\eps_i(\mu)-2} \Z[q^{-1}],\quad 
y_{\mu,\nu;i}(q) \in q^{\phi_i(\mu)-2} \Z[q^{-1}].$$
\item[(v)] For $\mu \in \Par_\al$ and $\nu \in \RPar_\al$, we have that
$$
S_\mu = \sum_{\la \in \RPar_\al} d_{\mu,\la}(q) D_\la,
\qquad
D_\nu = \sum_{\la \in \Par_\al} p_{\nu,\la}(-q) S_\la.
$$
\item[(vi)]
The vectors
$$\{S_\mu\:|\:\mu \in \RPar\}$$
give a basis for $V(\La)_\Laurent^*$ as a free $\Laurent$-module.
Moreover, any $S_\nu$ for $\nu \in \Par_\alpha \setminus \RPar_\alpha$ can be 
expressed as a $q \Z[q]$-linear combination of
$S_\mu$'s for $\mu \in \RPar_\alpha$ with $\mu \lex \nu$.
\item[(vii)]
Given $\mu \in \RPar_\alpha$, the difference
$S_\mu - \overline{S_\mu}$ is a 
$\Z[q,q^{-1}]$-linear combination of $S_\nu$'s
for $\nu \in\RPar_\alpha$ with $\nu \lex \mu$.
\end{itemize}
\end{Theorem}

We point out that the theorem above can be strengthened in the sense that $\lex$ can be replaced with $\lhd$ everywhere. The only proof of this fact we know is indirect---it uses categorification results described in section~\ref{SCat}. 

Set
\begin{equation}
Y_\mu := \pi(P_\mu)\qquad(\mu \in \RPar).
\end{equation}
The following theorem collects the necessary information on the quasi-canonical basis of $V(\La)$. 

\begin{Theorem}\label{anotherb} \cite[Theorem 3.14]{BKllt} 
Let $\kappa$ and $\La=\La(\kappa)$ be as in (\ref{EKappa}) and (\ref{ELa}).
\begin{itemize}
\item[(i)] The vectors
$$\{Y_\mu\:|\:\mu \in \RPar\}$$
give a basis for $V(\La)_\Laurent$, which we referr to as the quasi-canonical basis of $V(\La)$. 
\item[(ii)] 
For all $\mu \in \RPar$.
\begin{equation}
Y_\mu = \sum_{\nu \in \Par} d_{\nu,\mu}(q) S_\nu
\end{equation}

\item[(iii)] For $\mu, \nu \in \RPar$, we have 
$$\langle Y_\mu, D_\nu \rangle =  \delta_{\mu,\nu}.$$
\end{itemize}
\end{Theorem}

While the quasi-canonical basis $\{Y_\mu\:|\:\mu \in \RPar\}$ of $V(\La)$ will be of importance, the {\em canonical}\, basis will not play an important role in this paper. For an interested reader, recalling the definition  of the defect from (\ref{defdef}), the relation between the  quasi-canonical basis and the canonical basis is as follows

\begin{Proposition}\label{party} \cite[Lemma 3.12]{BKllt}
The canonical basis for $V(\La)$
is $$\bigcup_{\alpha \in Q_+} \{q^{-\defect(\alpha)}Y_\mu\:|\:\mu \in \RPar_\alpha\}.$$ 
In particular, we have that
$\overline{Y_\mu} = q^{-2\defect(\alpha)} Y_\mu$
for each $\mu \in \RPar_\alpha$.
\end{Proposition}

\section{Categorifications}\label{SCat}

We continue working with the notation as in the previous sections,
so $F$ is an algebraically closed field,  $\xi \in F^\times$ is a fixed a parameter, $e$ is the corresponding quantum characteristic, and $\La$ is a fixed dominant weight. To this data and every $\alpha \in Q_+$
we have associated a block $H^\La_\alpha$ of a
cyclotomic Hecke algebra with parameter $\xi \in F^\times$, and a block $H^\infty_\alpha$ of the affine Khovanov-Lauda-Rouquier algebra.

\subsection{Categorification of $\mathbf f$}\label{sklth}
For affine Khovanov-Lauda-Rouquier algebras, just like for cyclotomic Hecke algebras, it is convenient to 
abbreviate the direct sums of all our Grothendieck groups by
\begin{equation}\label{like}
[\Proj{H^\infty}] := \bigoplus_{\alpha \in Q_+} [\Proj{H^\infty_\alpha}],
\qquad
[\Rep{H^\infty}] := \bigoplus_{\alpha \in Q_+} [\Rep{H^\infty_\alpha}].
\end{equation}
Also, for $\alpha,\beta \in Q_+$,
we identify the Grothendieck group $[\Proj{H^\infty_{\alpha,\beta}}]$ with
$[\Proj{H^\infty_\alpha}] \otimes_{\Laurent} [\Proj{H^\infty_\beta}]$
so that $[P \boxtimes Q]$ is identified with $[P] \otimes [Q]$.

The next categorification result has been proved by Khovanov and Lauda for arbitrary type \cite{KL1,KL2}. Recall the duality $\#$ from (\ref{ESharp}). 

\begin{Theorem} \cite[section~3]{KL1}\label{klthm}
There is a unique $\Laurent$-module isomorphism
$$
\gamma:{_\Laurent}\mathbf f \stackrel{\sim}{\rightarrow} [\Proj{H^\infty}]
$$
such that 
$1 \mapsto [H^\infty_0]$ (the class of the left
regular representation of
the trivial algebra $H^\infty_0$) 
and $\gamma(\theta_i^{(n)} x) = 
\theta_i^{(n)}(\gamma(x))$ for each $x \in {_\Laurent}\mathbf f$, 
$i \in I$ and $n \geq 1$.
Under this isomorphism:
\begin{itemize}
\item[(i)] ${_\Laurent} \mathbf f_\alpha$ corresponds to $[\Proj{H^\infty_\alpha}]$ for any $\al\in Q_+$;
\item[(ii)] the 
multiplication $m_{\alpha,\beta}:{_\Laurent} \mathbf f_\alpha \otimes_{\Laurent} {_\Laurent}\mathbf f_\beta
\rightarrow {_\Laurent} \mathbf f_{\alpha+\beta}$ 
corresponds to the induction product induced by the exact
induction functor 
$\Ind_{\alpha,\beta}^{\alpha+\beta}$;
\item[(iii)] the comultiplication $m_{\alpha,\beta}^*:{_\Laurent} \mathbf f_{\alpha+\beta} 
\rightarrow {_\Laurent} \mathbf f_\alpha \otimes_{\Laurent}
{_\Laurent} \mathbf f_\beta$
corresponds to the restriction coproduct induced by the 
exact restriction functor
$\Res_{\alpha,\beta}^{\alpha+\beta}$;
\item[(iv)] the bar-involution on ${_\Laurent} \mathbf f_{\alpha}$ 
corresponds to the
anti-linear involution induced by the duality $\#$.
\end{itemize}
\end{Theorem}

The following theorem from \cite[Theorem 5.19]{BKllt} shows that the isomorphism $\gamma$ from Theorem~\ref{klthm} has an additional nice property, provided $\cha F=0$. Namely, $\gamma$ identifies the Lusztig's canonical basis $\mathbf B$ in $\mathbf f$ with the basis of the Grothendieck group $[\Proj{H^\infty}]$ consisting of the classes of projective indecomposable modules. This 
proves for all type $A$ quivers (finite or affine) a conjecture of 
Khovanov and Lauda formulated in \cite[$\S$3.4]{KL1}. 
Apart from the case $e=2$, this result is also proved 
in \cite{VV3} by a very different method, which works for an arbitrary simply-laced type. The result has also been announced by Rouquier for an arbitrary type.

\begin{Theorem}\label{klcon}
Assume that $\cha F = 0$.
For every $\alpha \in Q_+$,
the isomorphism $\ga:\mathbf f_\alpha \rightarrow [\Proj{H^\infty_\alpha}]$ from Theorem~\ref{klthm}
maps 
$\mathbf B_\alpha$
to the basis of
$[\Proj{H^\infty_\alpha}]$ 
arising
from the isomorphism classes of the indecomposable projective
graded $H^\infty_\alpha$-modules normalized so that 
they are self-dual with respect to the duality
$\#$.
\end{Theorem}

\subsection{\boldmath Categorification of $V(\La)$}\label{sc}
We now connect representation theory of $H_d^\La$ and the $U_q(\g)$-module $V(\La)$. 

The exact functors 
$$E_i^{(n)},F_i^{(n)},K_i$$ 
from (\ref{EEFKRep}) and (\ref{EEFKProj})  induce
$\Laurent$-linear endomorphisms of the
Grothendieck groups 
$[\Rep{H^\La}]$ and $[\Proj{H^\La}]$ from (\ref{EGrRepHLa}) and (\ref{EGrProjHLa}). 

In view of Theorem~\ref{clas}, $[\Rep{H^\La}]$
is a free $\Laurent$-module on basis
$$\{[D(\mu)]\:|\:\mu \in \RPar\}.$$ 
Also 
let $Y(\mu)$ denote the projective cover of $D(\mu)$
in $\Rep{H^\La_\alpha}$, for each $\mu \in \RPar_\alpha$.
Thus there is a degree-preserving surjection
$$
Y(\mu) \twoheadrightarrow D(\mu).
$$
The classes $$\{[Y(\mu)]\:|\:\mu \in \RPar\}$$ give a basis
for $[\Proj{H^\La}]$ as a free $\Laurent$-module. 

First of all, recalling the map (\ref{flatmap}), we get the following explicit connection between categorifications for $\mathbf f$ and $V(\La)$.

\begin{Proposition}\label{id1} \cite[Proposition 4.16]{BKllt}
There is a unique $\Laurent$-module isomorphism
$\db$ making the following diagram commutative
$$
\begin{CD}
\mathbf {_\Laurent}\mathbf f&@>\sim>\ga>&[\Proj{H^\infty}]\\
@V\beta VV&&@VV\pr V\\
V(\La)_\Laurent&@>\sim>\db>&[\Proj{H^\La}],
\end{CD}
$$
where $\beta$ denotes the surjection
$x \mapsto x^\flat v_\La$, $\ga$ is the isomorphism
from Theorem~\ref{klthm}, and $\pr$
is the $\Laurent$-linear map induced by the additive
functor (\ref{prprpr}).
Moreover:
\begin{enumerate}
\item[{\rm (i)}] For every $i \in I$ and $n \geq 1$, the map $\db$ intertwines the action
of $F_i^{(n)} \in U_q(\g)_\Laurent$ on $V(\La)_\Laurent$
with the endomorphism of $[\Proj{H^\La}]$
induced by the divided power functor $F_i^{(n)}$.
\item[{\rm (ii)}] For every $i \in I$, the map $\db$ intertwines the action
of $K_i \in U_q(\g)_\Laurent$ on $V(\La)_\Laurent$
with the endomorphism of $[\Proj{H^\La}]$
induced by the  functor $K_i$.
\end{enumerate} 
\end{Proposition}

Now we want to state the fundamental theorem 
which makes precise a sense in which
$\Proj{H^\La}$ {categorifies} the $U_q(\g)_\Laurent$-module
$V(\La)_\Laurent$ and
$\Rep{H^\La}$ {categorifies} 
$V(\La)_\Laurent^*$.

Recall the Cartan pairing
$$\langle.,.\rangle:[\Proj{H^\La}] \times [\Rep{H^\La}]
\rightarrow \Laurent$$
from section~\ref{SSGRT}
and the Shapovalov pairing
$$\langle.,.\rangle:V(\La)_\Laurent \times V(\La)_\Laurent^*
\rightarrow \Laurent$$
from section~\ref{md}.
Let
\begin{equation}\label{dstar}
\dd: [\Rep{H^\La}] \stackrel{\sim}{\rightarrow}
V(\La)_\Laurent^*
\end{equation}
be the dual map
to the isomorphism 
$\db$ of Proposition~\ref{id1}
with respect to these pairings.

Also let
\begin{equation}
\circledast: [\Rep{H^\La}] \rightarrow [\Rep{H^\La}]
\end{equation}
be the anti-linear involution induced by the duality
$\circledast$ from section \ref{SSCHANew}.

\begin{Theorem}\label{cat} \cite[Theorem 4.18]{BKllt}
The following diagram commutes:
$$
\begin{CD}
V(\La)_\Laurent&@>\sim > \db>&[\Proj{H^\La}]\\
@Va VV&&@VVbV\\
V(\La)_\Laurent^*&@<\sim <\dd<&[\Rep{H^\La}],
\end{CD}
$$
where $$a:V(\La)_\Laurent \hookrightarrow V(\La)_\Laurent^*$$
is the canonical inclusion,
and
$$b:[\Proj{H^\La}]  \rightarrow [\Rep{H^\La}]$$
is the $\Laurent$-linear map induced by the natural inclusion
of $\Proj{H_\al^\La}$ into $\Rep{H^\La_\al}$ for each $\al \in Q_+$.
Furthermore:
\begin{itemize}
\item[(i)] $b$ is injective and becomes an isomorphism
over $\Q(q)$;
\item[(ii)] both maps $\db$ and $\dd$ 
commute with the actions of all $E_i^{(n)}, F_i^{(n)}$ and $K_i$;
\item[(iii)] 
both maps $\db$ and $\dd$ intertwine
the involution $\circledast$ coming from
duality with the bar-involution;
\item[(iv)] the isomorphism $\db$ 
identifies 
the Shapovalov form on $V(\La)_\Laurent$
with the Cartan form on $[\Proj{H^\La}]$.
\end{itemize}
\end{Theorem}

Theorem~\ref{cat} reveals a very close connection between representation theory of cyclotomic Hecke algebras and integrable modules over quantized enveloping algebras of Kac-Moody algebras. In the next sections we review more results along these lines, which indicate that the connection is really very deep. 

\subsection{Monomial bases and Specht modules}\label{SSMBSM}
We now use $\dd$ and $\db$ to identify classes in the Grothendieck group $[\Rep{H^\La}]$ of various families of modules  with special families of vectors in $V(\La)$ defined in sections~\ref{SSFock} and \ref{ssb}.  

First, the classes of Specht modules get identified with the standard monomials in $V(\La)$ 
defined in (\ref{stmn}):

\begin{Theorem}\label{ids} \cite[Theorem 5.6]{BKllt} 
For each $\mu \in \Par_\al$, we have 
$\dd([S(\mu)]) = S_\mu$.
\end{Theorem}

Now using Theorem~\ref{tdcb}(vi) and (\ref{sact1}), we get

\begin{Corollary}\label{sbase}
The classes $\{[S(\mu)]\:|\:\mu \in \RPar_\alpha\}$
give a basis for $\Rep{H^\La_\alpha}$ as a free
$\Laurent$-module.
\end{Corollary}

\begin{Corollary}\label{CSpechtCor}
For $\mu \in \Par$ and $i \in I$, the following hold in $[\Rep{H^\La}]$:
\begin{align*}
E_i [S(\mu)]&=\sum_A q^{d_A(\mu)} [S(\mu_A)],
&F_i [S(\mu)]&=\sum_B q^{-d^B(\mu)}[S(\mu^B)],
\end{align*}
where the first sum is over all removable $i$-nodes $A$ for $\mu$,
and the second sum is over all addable $i$-nodes $B$ for $\mu$.
\end{Corollary}

The first formula in Corollary~\ref{CSpechtCor} of course also follows easily from Theorem~\ref{TBr}.

\subsection{Canonical bases and graded decomposition numbers}\label{SSCBGDN}
Throughout this section we assume that the characteristic of the ground field $F$ is zero. Under this assumption we can identify the elements of the dual canonical basis with the classes of the irreducible modules and elements of the quasi-canonical basis with the classes of the projective indecomposable modules in an explicit manner in such a way that the {\em combinatorial labels match}. 

The results of this section should be viewed as a graded version of 
the Lascoux-Leclerc-Thibon conjecture (generalized to higher levels), and hence a graded version of Ariki's Categorification Theorem. 

\begin{Theorem}\label{late}  \cite[Theorem 5.14]{BKllt}
Assume that $\cha F = 0$.
For each $\mu \in \RPar_\alpha$, we have that
$$\db(Y_\mu) = [Y(\mu)]\qquad \text{and}\qquad \dd([D(\mu)]) = D_\mu,$$
where $\db$ and $\dd$ are the maps from Theorem~\ref{cat}.
\end{Theorem}

Thus in the case $\cha F=0$, we get a description of the graded decomposition numbers in terms of the polynomials $d_{\mu,\nu}(q)$ defined using dual canonical basis in (\ref{ml}). 

\begin{Theorem}\label{maindecthm}
Assume that $\cha F = 0$.
For $\mu \in \Par_\al$, we have that
$$
[S(\mu)] = \sum_{\nu \in \RPar_\alpha}
d_{\mu,\nu}(q) [D(\nu)].
$$
In other words, for $\mu \in\Par_\alpha$ and $\nu \in \RPar_\alpha$,
we have that
$$
[S(\mu):D(\nu)]_q = d_{\mu,\nu}(q).
$$
Moreover, for all such $\mu,\nu$, we have that
$d_{\mu,\nu}(q) = 0$ unless $\nu \unlhd \mu$.
\end{Theorem}

We point out again that Theorem~\ref{maindecthm} is a graded analogue (or a $q$-analogue) of the Ariki result on decomposition numbers which says that 
$$
[S(\mu):D(\nu)] = d_{\mu,\nu}(1).
$$

\begin{Corollary}\label{CMainCor}
Assume that $\cha F = 0$, $\al\in Q_+$ be of height $d$, $\mu \in \Par_\al$, and $\nu\in \RPar_d$. 
Then: 
\begin{enumerate}
\item[{\rm (i)}] $[S(\mu):D(\nu)]_q=0$ unless $\nu\in\RPar_\al$ and $\nu\unlhd\mu$.
\item[{\rm (ii)}] If $\mu\in\RPar$ then $[S(\mu):D(\mu)]_q=1$.
\item[{\rm (iii)}] If $\nu\neq\mu$, then $[S(\mu):D(\nu)]_q\in q\Z_{\geq 0}[q]$.  
\end{enumerate}
\end{Corollary}

Corollary~\ref{CMainCor} makes it possible in principle to compute the graded decomposition numbers $[S(\mu):D(\nu)]_q=0$ by induction on the dominance order on $\Par_\al$. The induction base is clear. Now, let $\mu\in\Par_\al$, and assume that we know the graded decomposition numbers $d_{\nu,\la}$ for all $\nu\lhd\mu$. Then Corollary~\ref{CCHS}  and the unitrangularity of the decomposition matrix coming from Corollary~\ref{CMainCor} allow us to compute all $\CH D(\nu)$ for $\nu\lhd\mu$. Now, we have 
\begin{equation}\label{EQCharEq}
\CH S(\mu)=\sum_{\nu\in\RPar_\al,\ \nu\unlhd\mu}d_{\mu,\nu}\,\CH D(\nu),
\end{equation}
with $d_{\mu,\mu}=1$ if $\mu$ is restricted. The left-hand side of the equation above is known from Corollary~\ref{CCHS}. 


Moreover, we know that the $q$-characters of all irreducible modules $D(\la)$ are bar-invariant and non-negative, which means that each $\bi$ appears in $\CH D(\la)$ with a bar-invariant coefficient belonging to $\Z_{\geq 0}[q,q^{-1}]$. By Corollary~\ref{CMainCor}, we also know that the graded decomposition numbers $d_{\mu,\nu}$ for $\nu\lhd\mu$ belong to $q\Z_{\geq 0}[q]$. 
It remains to note, using linear independence of characters of irreducible modules noted in Theorem~\ref{ch}, that there are unique polynomials $d_{\mu,\nu}\in q\Z_{\geq 0}[q]$ such that 
\begin{equation*}
S(\mu)-\sum_{\nu\in\RPar_\al,\ \nu\lhd\mu}d_{\mu,\nu}\,\CH D(\nu)
\end{equation*}
is bar-invariant. Moreover, since all $\CH D(\nu)$ for $\nu\lhd \mu$ are known by induction and the coefficients everywhere are non-negative, it is possible to find these $d_{\mu,\nu}$ in a finite number of steps. 

The procedure just described is quite cumbersome to perform in practice by hand. In level one, there is a much more efficient algorithm. We explain it in the next subsection. 

\subsection{\boldmath An algorithm for computing decomposition numbers for $H_d(\C,\xi)$}\label{SSAlg} In this section we describe a rather efficient algorithm 
for computing decomposition numbers for the Hecke algebra $H_d(F,\xi)$ when the characteristic of the ground field $F$ is $0$. We point out that this algorithm is {\em equivalent}\, to the Lascoux-Leclerc-Thibon algorithm \cite{LLT}, although our interpretation of it is  slightly different, cf. \cite{KN}. In particular, our algorithm proceeds `along the rows' of the decomposition matrix, rather than `along the columns', and it relies on $q$-characters of Specht modules. 

Throughout the subsection we are assuming that the characteristic of $F$ equals $0$, and that $e>0$, i.e. $\xi\in F$ is a primitive $e$th root of unity. We also assume that $\La=\La_0$, that is we work with the algebra $H_d=H_d(F,\xi)=H_d^{\La_0}(F,\xi)$. 
(Of course, when $e=0$, the parameter $\xi$ is generic, and the Specht modules over $H_d$ are all irreducible, so this case is not interesting for us here.) 

Fix $\al\in Q_+$ of height $d$. Then $\Par_\al$ (resp. $\RPar_\al$) is the set of all (resp. all $e$-restricted) partitions of $d$ of residue content $\al$, cf. Example~\ref{c3}(i) and (\ref{EContent}). 

First, we consider the following 

\vspace{2mm}
\noindent
{\bf Basic Task.}
Suppose $d(q) \in q \Z_{\geq 0}[q]$, and $m(q),r(q) \in \Z_{\geq 0}[q,q^{-1}]$ are such that $\overline{m(q)}=m(q)$, $\overline{r(q)}=r(q)$, and $r(q) \neq 0$.  If $d(q)r(q) + m(q)$ and $r(q)$ are known, find $d(q)$ and $m(q)$.

\vspace{2 mm}
We note an easy 

\vspace{2 mm}
\noindent
{\bf Algorithm for Solving the Basic Task.}  
Denote the top term of $r(q)$ by $bq^R$ (for $b>0$ and $R \geq 0$). The algorithm goes by induction on the number of non-zero terms in $d(q)r(q) + m(q)$. The induction base is when $d(q)r(q) + m(q)=0$, which implies $d(q)=m(q)=0$ since all coefficients are non-negative.  Let $d(q)r(q) + m(q) \neq 0$ and write
$$d(q)r(q) + m(q)=a_{-N}q^{-N} + \dots + a_Mq^M$$
for $M,N \geq 0$ and $a_{-N},a_M>0$.  Since $d(q) \in q \Z_{\geq 0}[q]$ and $m(q)$ and $r(q)$ are bar-invariant, we have $M \geq N$.  In order to make the inductive step, we  consider two cases.

\underline{Case 1: $N < M$.} As $m(q)$ is bar-invariant, the term $a_Mq^M$ must come from $d(q)r(q)$.  Thus $\frac{a_M}{b}q^{M-R}$ is a term in $d(q)$.  Setting $d'(q):=d(q)-\frac{a_M}{b}q^{M-R}$, we are reduced to solving the Basic Task for $d'(q)r(q) + m(q)$ which has strictly fewer terms.

\underline{Case 2: $M = N$.} As $d(q) \in q \Z_{\geq 0}[q]$, the term $a_{-N}q^{-N}$ must therefore come from $m(q)$.  Since $m(q)$ is bar-invariant, $a_{-N}q^N$ must be a term in $m(q)$ also.  Setting $m'(q):=m(q)-(a_{-N}q^{-N} + a_{-N}q^N)$, we are reduced to solving the Basic Task for $d(q)r(q) + m'(q)$ which has strictly fewer terms.

\vspace{2mm}
For $\la\in\RPar_\al$ we define 
$$\bj^\la=(j_1,\dots,j_d)\in I^\al$$ 
as follows. Let $A$ be the bottom removable node of $\la$ such that $\la_A$ is $e$-restricted. Now $\bj^\la$ is determined inductively from $j_d:=\res\, A$ and $(j_1,\dots,j_{d-1}):=\bj^{\la_A}$. 
We note that in general $\bj^\la\neq \bi^\la$, cf. (\ref{EBIMu}). 
For any finite dimensional $H_d$-module $M$ denote by 
$
m_\la(M)
$
the multiplicity of $\bj^\la$ in $\CH M$. 
A key combinatorial property of $\bj^\la$ is as follows

\begin{Proposition}\label{JLambda} \cite{KN} Let $\mu$ be a partition of $d$ and $\lambda$ be an $e$-restricted partition of $d$.  Then: 
\begin{enumerate}
\item[{\rm (i)}] $m_\la(S(\mu))=0$ unless $\mu\unrhd\la$. 
\item[{\rm (ii)}] If $\mu$ is also $e$-restricted, then $m_\la (D(\mu))=0$ unless $\mu\unrhd\la$. 
\item[{\rm (iii)}]  $m_\la (D(\la))=m_\la(S(\la))$.
\item[{\rm (iv)}] We have 
\begin{equation}\label{EMain}
m_{\lambda}(S(\mu)) = r_{\lambda} d_{\mu,\lambda} +\sum_{\nu}{d_{\mu,\nu} m_{\lambda}(D(\nu))},
\end{equation}
where the sum is over all $e$-restricted partitions $\nu$ satisfying $\mu \unrhd \nu \rhd \lambda$. 
\end{enumerate}
\end{Proposition}

The multiplicity 
$$
r_\la:=m_\la(S(\la))=m_\la (D(\la))
$$
is easy to compute explicitly as follows. Let $A$ be the bottom removable node of $\la$ such that $\la_A$ is $e$-restricted. Assume that there are $r_1$ removable nodes $A_1=A,A_{2},\dots, A_{r_1}$ weakly below $A$. We call these nodes the {\em bottom removable sequence of $\la$}. Denote $r_1$ by $r(\la)$. Now, let $\la^{(2)}$ be the ($e$-regular) partition obtained from $\la$ by removing the bottom removable sequence of $\la$. Now define $r_2:=r(\la^{(2)})$, remove the bottom removable sequence from $\la^{(2)}$ to get $\la^{(3)}$, and so on until we reach the empty partition. This defines the sequence of positive integers $r_1,r_2,\dots, r_t$ such that $r_1+r_2+\dots+r_t=d$. It is proved in \cite{KN} that
\begin{equation}\label{RRLa}
r_\la=[r_1]![r_2]!\dots[r_t]!.
\end{equation}

Let us fix $\mu\in\Par_\al$ and $\la\in\RPar_\al$. We finally exhibit an algorithm for computing the graded decomposition number $d_{\mu,\lambda}$ by induction.  In fact, this induction requires us to keep track of some extra $q$-character information, so we now carefully describe exactly how the induction goes.

Define the set
$$P_{\mu,\lambda}:=\{(\nu,\kappa) | \nu \in \Par_\al,\ \kappa \in \RPar_\al, \ \mu \unrhd \nu \unrhd \kappa \unrhd \lambda \}.$$
Let us take $``<"$ to be any total ordering on $P_{\mu,\la}$ satisfying the following two conditions:

(1) If $(\nu,\kappa), (\theta,\phi) \in P_{\mu,\lambda}$, then $\nu \lhd \theta$ implies $(\nu,\kappa)< (\theta,\phi)$;

(2) If $(\nu,\kappa),(\nu,\phi) \in P_{\mu,\lambda}$, then $\kappa \rhd \phi$ implies $(\nu,\kappa)<(\nu,\phi)$.

We now calculate $d_{\mu,\lambda}$ and $m_{\lambda}(D(\mu))$ by induction on the total order $``<"$ on $P_{\mu,\lambda}$. The induction begins at the smallest element $(\lambda,\lambda)$ in $P_{\mu,\lambda}$, where we know that $d_{\lambda,\lambda}=1$ and $m_{\lambda}(D(\lambda)) = r_{\lambda}$ by Proposition~\ref{JLambda}(iii), which is known by (\ref{RRLa}).

Now let $(\nu,\kappa) > (\lambda,\lambda)$.  If $\nu=\kappa$ then we have that $d_{\nu,\nu}=1$ and $m_{\nu}(D(\nu))=r_{\nu}$, so we may assume that $\nu \rhd \kappa$.  By induction, we know the decomposition numbers $d_{\nu,\phi}$ for all $\phi \in \RPar_\al$ satisfying $\nu \unrhd \phi \rhd \kappa$.  Also by induction we know the multiplicities $m_{\kappa}(D(\phi))$ for all $\phi \in \RPar_\al$ satisfying $\nu \rhd \phi \unrhd \kappa$.  To make the inductive step we need to compute $d_{\nu,\kappa}$ and, if $\nu$ is $e$-restricted, $m_{\kappa}(D(\nu))$.

If $\nu$ is not $e$-restricted, then by Proposition~\ref{JLambda}(iv), we have
$$ d_{\nu,\kappa}=\frac{1}{r_{\kappa}}\left(m_{\kappa}(S(\nu)) - \sum_{\phi \in \RPar_\al \text{,~} \nu \rhd \phi \rhd \kappa}{d_{\nu,\phi} m_{\kappa}(D(\phi))}\right),$$
where all terms in the right-hand side are known by induction and Corollary~\ref{CCHS}.

So let $\nu$ be $e$-restricted.  By Proposition~\ref{JLambda}(iv) again, we have 
$$m_{\kappa}(D(\nu)) + d_{\nu,\kappa} r_{\kappa} = m_{\kappa}(S(\nu)) - \sum_{\phi \in \RPar_\al \text{,~} \nu \rhd \phi \rhd \kappa}{d_{\nu,\phi} m_{\kappa}(D(\phi))},$$
where all terms in the right-hand side are known by induction and Corollary~\ref{CCHS}.

Note $r_{\kappa}$ is non-zero and bar-invariant, $d_{\nu,\kappa} \in q \Z_{\geq 0}[q]$ by Corollary~\ref{CMainCor},  and $m_{\kappa}(D(\nu))$ is bar-invariant.  Hence, we are in the assumptions of the Basic Task above, with $m(q)= m_{\kappa}(D(\nu))$, $d(q)=d_{\nu,\kappa}$, and $r(q)=r_{\kappa}$. Now we apply the algorithm for solving the Basic Task described above to calculate $m_{\kappa}(D(\nu))$ and $d_{\nu,\kappa}$ and complete the inductive step.

\section{Reduction Modulo $p$ and James' Conjecture}\label{SJC}
As usual, let us fix throughout the section a dominant weight $\La \in P_+$ of level $l$
and $\alpha \in Q_+$ of height $d$.

\subsection{Realizability over prime subfields}\label{SRed}
In this section we consider base change, so we will need to work with more than one field. Let us introduce the notation which will allow us to handle this.

\begin{Definition}
{\rm 
For any field $K$, denote by $R^\La_\al(K)$ the graded $K$-algebra given abstractly  by generators (\ref{EBlockGenerators}) of degrees given by (\ref{EBlockDeg}) subject 
only to the (homogeneous) relations (\ref{R0})--(\ref{R7}) and (\ref{ECyclotomic}) 
for 
$\bi,\bj\in I^\al$ and all 
admissible $r, s$.  The algebra $R^\La_\al(K)$ is referred to as a {\em cyclotomic Khovanov-Lauda-Rouquier algebra}. 
}
\end{Definition}

Note that the algebra $R^\La_\al(K)$ depends on the quiver $\Gamma$ introduced in section~\ref{SSDynkin}, which in turn is determined by our fixed pair $(F,\xi)$ from section~\ref{SSGrFPar}. 
In view of Theorem~\ref{TBKBlocks}, we know that $R_\al^\La(F)$ is isomorphic to the block $H_\al^\La(F,\xi)$ of the cyclotomic Hecke algebra $H_d^\La(F,\xi)$ over our fixed field $F$ with parameter $\xi$. However, for an arbitrary field $K$ there is no apparent relation between the two algebras. 

First, we explain how to descend
from the field $F$ to its prime subfield $E$. Note that the parameter $\xi$ does not have to belong to $E$. However:

\begin{Theorem}\label{EF} \cite[Theorem 6.1]{BKyoung}
Let $H_\al^\La(F)=H_\al^\La(F,\xi)$ be a block of the cyclotomic Hecke algebra $H_d^\La(F,\xi)$.  
Let $E$ be the prime subfield of $F$.
Let $H^\La_\alpha(E)$ denote the $E$-subalgebra of 
$H^\La_\alpha$ generated by the elements (\ref{EBlockGenerators}).
Then the natural map 
$$
F \otimes_{E} H^\La_\alpha(E)
\rightarrow H^\La_\alpha(F)
$$ 
is an $F$-algebra isomorphism. 
Moreover, 
there is an $E$-algebra isomorphism
$$R^\La_\alpha(E) \stackrel{\sim}{\rightarrow}
H^\La_\alpha(E)$$ 
sending
the generators of $R^\La_\alpha(E)$ to 
the generators of $H^\La_\alpha(E)$ with the same names.
\end{Theorem}

The theorem easily implies that irreducible $H_\al^\La$-modules are realizable over prime subfields in the following sense: 

\begin{Corollary}\label{mcor}
Let $D(F)$ be an irreducible $H^\La_\alpha(F)$-module, and let $E$ be the prime subfield of $F$.
Then there exists an irreducible
$H^\La_\alpha(E)$-module $D(E)$
such that 
$D(F) \cong F \otimes_{E} D(E)$.
\end{Corollary}

Another easy consequence of Theorem~\ref{EF} is a statement, conjectured in  \cite[section 4]{JamConj} and \cite[Conjecture 6.38]{MathasB}, that the decomposition matrices of the cyclotomic Hecke algebras 
$H^\La_d(F)$ depend only on $e$
and the characteristic of $F$ (but not on $\xi$ or $F$ itself), see \cite[Corollary 6.3]{BKyoung} for more details.

\subsection{Reduction modulo $p$}\label{SSRedNew}

Now suppose that $K$ is a field of characteristic zero, and 
let $\zeta \in K^\times$ be a primitive $e$th root of unity
if $e > 0$, or some element that is not a root of unity if $e = 0$.

By Theorem~\ref{TBKBlocks}, applied to the pair $(K,\zeta)$ instead of $(F,\xi)$, we know that any block $H_\al^\La(K,\zeta)$ of the cyclotomic Hecke algebra over $K$ with parameter $\zeta$ is generated by an explicitly defined set (\ref{EBlockGenerators}) of elements 
subject only to the  relations (\ref{R0})--(\ref{R7}) and (\ref{ECyclotomic}) 
for $\bi,\bj\in I^\al$ and all admissible $r, s$. 

The following theorem explains how $H^\La_\alpha(F,\xi)$ 
can be obtained from
$H^\La_\alpha(K,\zeta)$ by a base change. The idea is that, since the quantum characteristic $e$ is the same for both pairs $(K,\zeta)$ and $(F,\xi)$, they yield the  Khovanov-Lauda-Rouquier algebras of the {\em same Lie type $\Gamma$}, but over different fields $K$ and $F$. 

\begin{Theorem}\label{Red} \cite[Theorem 6.4]{BKyoung}
Let $H^\La_\alpha(\Z)$
denote the 
subring of $H^\La_\alpha(K,\zeta)$ generated by
the elements (\ref{EBlockGenerators}).
Then $H^\La_\alpha(\Z)$ is a free $\Z$-module
and there are isomorphisms
\begin{align}\label{i1}
H^\La_\alpha(K,\zeta)&\stackrel{\sim}{\rightarrow} 
K \otimes_\Z H^\La_\alpha(\Z),\\
H^\La_\alpha(F,\xi)
&\stackrel{\sim}{\rightarrow} F \otimes_\Z H^\La_\alpha(\Z),\label{i2}
\end{align}
such that $e(\bi) \mapsto 1 \otimes e(\bi)$,
$y_r\mapsto 1 \otimes y_r$ and $\psi_r\mapsto 1 \otimes \psi_r$
for each $\bi$ and $r$.
\end{Theorem}

As a consequence
we can explain how to
reduce irreducible $H^\La_\alpha(K,\zeta)$-modules
modulo $p$ to obtain well-defined $H^\La_\alpha(F,\xi)$-modules. This resembles the usual  Brauer reduction modulo $p$ in finite group theory, with an important difference  that here the algebra $H^\La_\alpha(K,\zeta)$ over the field of characteristic zero is in general not semisimple.

\begin{Theorem}\label{bct} \cite[Theorem 6.5]{BKyoung} 
If $D(K)$ is an irreducible $H^\La_\alpha(K,\zeta)$-module
and $0 \neq v \in D(K)$,
then $$D(\Z) := H^\La_\alpha(\Z) v$$ is
a lattice in $D(K)$
that is invariant under the action of $H^\La_\alpha(\Z)$.
For any such lattice,
$$
D(F) := F \otimes_\Z D(\Z)
$$ 
is naturally an $H^\La_\alpha(F,\xi)$-module
with the same $q$-character as $D(K)$.
In particular,  the class  of $D(F)$
in the Grothendieck group $[\Rep{H^\La_\alpha(F,\xi)}]$ is independent of the choice of lattice.
\end{Theorem}

So we have a well-defined decomposition matrix corresponding to the  reduction modulo $p$ procedure which we just described. To distinguish this decomposition matrix from the decomposition matrices which describe multiplicities of the irreducible modules in Specht modules, this matrix is usually called the adjustment matrix \cite{JamConj}. The adjustment matrices are systematically discussed in the next section.

\subsection{Graded adjustment matrices}\label{sgam}


We want to collect the information coming from Theorem 10.5 into a {\em graded adjustment matrix}. So let us keep working with the set up of section~\ref{SSRedNew}. 

In order to avoid any confusion, let us use the notation 
$$
S_{F,\xi}(\mu),\ D_{F,\xi}(\nu)\qquad(\mu\in\Par_\al,\ \nu\in\RPar_\al)
$$
to denote the Specht and the irreducible modules over the algebra $H_\al^\La(F,\xi)$, defined as in section~\ref{SGrMod}, and the notation
$$
S_{K,\zeta}(\mu),\ D_{K,\zeta}(\nu)\qquad(\mu\in\Par_\al,\ \nu\in\RPar_\al)
$$
to denote the Specht and the irreducible modules over the algebra $H_\al^\La(K,\zeta)$ defined similarly. Note that the labeling sets $\Par_\al$ and $\RPar_\al$ are the same in these two cases, since they only depend on $e$ and $\kappa$. 

Even though this is not important for the adjustment matrices, we want to make the reduction modulo $p$ procedure canonical as follows. Fix for the moment 
$\mu\in\RPar_\al$. Recall the cellular basis 
$
\{C_\T\mid \T\in\St(\mu)\}
$ 
of the Specht module from section~\ref{SSSMIM} and the  special standard tableaux $\T^\mu$ from section~\ref{SSTab}. 
Let
$v_\mu$ denote the image of the  vector 
$C_{\T^\mu} \in S_{K,\zeta}(\mu)$ under
some surjection 
$S_{K,\zeta}(\mu) \twoheadrightarrow D_{K,\zeta}(\mu)$.

By Theorem~\ref{PVZ}, 
$S_{K,\zeta}(\mu)$ is generated as an $H^\La_\alpha(K,\zeta)$-module
by the vector $C_{\T^\mu}$, hence $v_\mu \in D_{K,\zeta}(\mu)$ is non-zero.
Now let
$D_{\Z}(\mu)$ denote the $\Z$-span of the
vectors
$\psi_{r_1} \cdots \psi_{r_m} 
y_1^{n_1} \cdots y_d^{n_d} v_\mu$
for all $m \geq 0, 1 \leq r_1, \dots, r_m < d$ and $n_1,\dots,n_d \geq 0$.
By Theorem~\ref{bct}, $D_{\Z}(\mu)$ is a lattice in $D_{K,\zeta}(\mu)$. 

Finally, set 
\begin{equation}\label{jam}
\overline{D_{\Z}(\mu)} := F \otimes_{\Z} D_{\Z}(\mu)\qquad(\mu\in\RPar_\al)
\end{equation}
By Theorem~\ref{PVZ}  again, $\overline{D_{\Z}(\mu)}$ 
is a well-defined graded $H^\La_\alpha(F,\xi)$-module, which has the same $q$-character
as $D_{K,\zeta}(\mu)$.
Hence, recalling Theorem~\ref{late}, 
$$
\dd([\overline{D_{\Z}(\mu)}]) = D_\mu \qquad(\mu\in\RPar_\al).
$$

\begin{Theorem}\label{adjt} \cite[Theorem 5.17]{BKllt}
Let $\mu\in\RPar_\al$ and $\la\in\Par_\al$. Then
$$
[\overline{D_{\Z}(\mu)}] = [D_{F,\xi}(\mu)] + \sum_{\nu \in \RPar_\alpha,\ \nu \lhd \mu} a_{\mu,\nu}(q) [D_{F,\xi}(\nu)].
$$
for some unique
bar-invariant Laurent polynomials $a_{\mu,\nu}(q) \in \Z_{\geq 0}[q,q^{-1}]$. Moreover
$$
[S_{F,\xi}(\la): D_{F,\xi}(\mu)]_q
= 
\sum_{\nu \in \RPar_\al}  [S_{K,\zeta}(\la): D_{K,\zeta}(\nu)]_qa_{\nu,\mu}(q).
$$
\end{Theorem}

We refer to the matrix $(a_{\mu,\nu}(q))_{\mu,\nu \in \RPar_\alpha}$
as the {\em graded adjustment matrix}.
For level one and 
$\xi = 1$, our graded adjustment matrix specializes at $q=1$
to the adjustment matrix defined originally by James \cite{JamConj}.
Curiously we did not yet find an example in which $a_{\mu,\nu}(q) \notin \Z$;
this is related to a question raised by
Turner in the introduction of \cite{T}.

\subsection{James Conjecture}\label{SSJamConj}
Keep working with the notation of sections~\ref{SSRedNew} and \ref{sgam}. In particular, $F$ is a field of characteristic $p>0$, $K$ is a field of characteristic zero, and the pairs $(K,\zeta)$ and $(F,\xi)$ have the same quantum characteristic $e$. Consider the 
graded decomposition matrix for $H_\al^\La(F,\xi)$:
$$
\Delta^\La_\al(F,\xi):=\big([S_{F,\xi}(\la): D_{F,\xi}(\mu)]_q\big)_{\la\in\Par_\al,\ \mu\in\RPar_\al}
$$
and the graded decomposition matrix for $H_\al^\La(K,\zeta)$:
$$
\Delta^\La_\al(K,\zeta):=\big([S_{K,\zeta}(\la): D_{K,\zeta}(\mu)]_q\big)_{\la\in\Par_\al,\ \mu\in\RPar_\al}.
$$
We also have the graded adjustment matrix defined in section~\ref{sgam}:
$$
A^\La_\al:=\big(a_{\mu,\nu}(q)\big)_{\mu,\nu\in\RPar_\al}.
$$
Theorem~\ref{adjt} implies that 
$$
\Delta^\La_\al(F,\xi)= \Delta^\La_\al(K,\zeta) A^\La_\al.
$$

Since the matrix $\Delta^\La_\al(K,\zeta)$ is quite well understood thanks to Corollary~\ref{mcor}, the main remaining step is  to understand the adjustment matrix $A^\La_\al$. It would be particularly interesting to find
hypotheses on $\alpha$ that ensure that the adjustment matrix is the identity matrix. This is equivalent of course to the statement that all irreducible modules $D_{K,\zeta}(\mu)$ in the block corresponding to $\al$ remain irreducible when reduced modulo $p$, and so in particular $\CH D_{K,\zeta}(\mu)=\CH D_{F,\xi}(\mu)$ for all $\mu\in\RPar_\al$.

In level one, James suggested a precise conjecture for this \cite[section~4]{JamConj}. In a block form James Conjecture can be formulated as follows. 

Assume that $\La=\La_0$. Recall a well known description of the weight spaces of the {\em basic module}  $V(\La_0)$ over $\g$, see e.g. \cite{Kac}. Let $\al\in Q_+$, $W$ be the Weyl group of $\g$, $\delta$ be the null-root for $\g$, and $\nabla=\La_0-\al$ be an arbitrary weight appearing in the module $V(\La_0)$. Then there exist a unique element $w\in W$ and a unique non-negative integer $m$ such that 
$$
\nabla=w\La_0-m\delta.
$$

We point out for the reader who would prefer to switch from Lie-theoretic notation to a combinatorial one, that $w\La_0$ corresponds to the $e$-core of the partitions $\mu$ in the block corresponding to $\al$, and $m$ is what is sometimes called the {\em $e$-weight of the block}, that is the number of rim $e$-hooks which need to be removed from partitions before we reach the $e$-core.

\begin{Conjecture}
{\rm 
{\bf (James Conjecture)}  \cite[section~4]{JamConj} 
Let $F$ be of characteristic $p>0$. 
If $m<p$ then the adjustment matrix 
$A^{\La_0}_\al$ is the identity matrix. 
}
\end{Conjecture}

The reader is referred to \cite[$\S$2]{Geck} and \cite{F} for further comments and results on James Conjecture. For example Fayers \cite{F} proves that the assumption $m<p$ is essentially the best possible. For a parallel theory for algebraic groups, see \cite{AJS}.

Recalling the notion of defect $\defect(\al)$ from (\ref{defdef}), it is not hard to observe that in the case $\La=\La_0$, we have $m=\defect(\al)$, cf. \cite{FayersWeight}. 

So one could speculate the following  generalization beyond level $1$:

\begin{Conjecture}
{\rm 
{\bf (Higher Level James Conjecture)}
Let $F$ be of characteristic $p>0$. 
If $$\defect(\al)<p,$$ then the adjustment matrix 
$A^{\La}_\al$ is the identity matrix. 
}
\end{Conjecture}

A less bold speculation is obtained if we assume $\defect(\al)<p/l$ in the conjecture above. 

More developments and partial results related to James Conjecture can be found in  \cite{F,FW4,Geck,GeckKLJ,GMu}.

\section{Some other results}\label{SSBr}
In this section we review some other important results which 
are representative of how representation theory of symmetric groups has been developing in the last fifteen years. 

\subsection{Blocks of symmetric groups: Brou\'e's conjecture and Chuang-Rouquier equivalences}\label{SSCR}
Block theory of finite groups is an extremly vibrant area of research nowadays.
General conjectures of Alperin, Brou\'e, Dade and others 
on the structure of blocks and on relations between arithmetic invariants of blocks provide the main motivation for research, while symmetric groups and their double covers are often the main testing ground for general conjectures. 
In this section we discuss Brou\'e-type conjectures for symmetric groups, which were recently proved by Chuang and Rouquier in the spectacular paper \cite{CR}. 

The idea is that even though there are very many blocks around, many of them 
seem to have much in common. So one tries to look for some sort of equivalence 
between blocks whose invariants coincide. Ideally, these will
be Morita equivalences, as in the results of Scopes \cite{Scopes}, 
Kessar \cite{Ke} and Chuang-Kessar \cite{CK}.
However, it is often the case that coarser
block invariants (such as Cartan invariants, cf. \cite{BKdet,Hill}) 
coincide, but for example decomposition numbers do not.
So one cannot expect a Morita equivalence to hold in general.
Instead, Brou\'e and Rickard have suggested weaker versions of equivalence, such as perfect isometry, derived equivalence, stable equivalence, etc., 
see for example \cite{Br1, Ri1, Ri2}. 

The famous abelian defect group conjecture of Brou\'e claims that a $p$-block of a finite group $G$ with abelian defect group $D$ should be derived
equivalent to
its Brauer correspondent in $N_G(D)$ (over a suitable local ring 
with field of fractions of characteristic $0$ and 
residue field of characteristic $p$). 
Thus, in the case of a 
block $B$ of the symmetric group $S_d$ of $p$-weight $\omega<p$,
there should be a derived equivalence
between $B$ and the group algebra of the 
wreath product $(C_p \rtimes C_{p-1}) \wr S_\omega$.
We note that somewhat mysteriously the assumption on $p$ coming from the Brou\'e Conjecture agrees with the assumption on $p$ in the James Conjecture, see  section~\ref{SSJamConj}.

In view of the work of Marcus \cite{Marcus} and Chuang-Kessar \cite{CK}, this conjecture for the case of symmetric groups follows from a conjecture of Rickard (which Rickard himself has proved for $p$-weights $\leq 5$). The following result simply says that Rickard's Conjecture is true:

\begin{Theorem}\label{TCR} {\rm \cite{CR}} 
All blocks of symmetric groups with a fixed $p$-weight  
are derived equivalent.
\end{Theorem}

Note that this result has a very natural interpretation in terms of the categorification described in section~\ref{sc}. Recall that under the categorification the blocks of the symmetric groups correspond to the weight spaces of the module $V(\La_0)$. It is easy to see that  two blocks have the same $p$-weight if and only if they belong to the same $W$-orbit, where $W$ is the Weyl group of $\g$ acting naturally on the weights of $V(\La_0)$. So the Chuang-Rouquier derived equivalences simply `lift' the action of the Weyl group $W$ from Grothendieck groups to derived categories.  

Note also that, unlike in Brou\'e's conjecture, there is no restriction
 $\omega < p$ in Theorem~\ref{TCR}. We also point out that the work \cite{CR} is much more general than Theorem~\ref{TCR} might suggest. It develops a general set-up  called {\em ${\mathfrak sl}_2$-categorification} which allows one to deduce results like Theorem~\ref{TCR} by checking some rather short and elegant set of axioms. These ideas were developed even further in \cite{Ro} and \cite{KL3}.  The methods of \cite{CR} can be thought of as far-reaching generalizations of the theory developed in \cite{Gr,KBrII,KBrIII}. 

We finish this section by saying that block theory of symmetric groups and related areas have received a lot of attention recently. We refer the interested reader to the following literature for more information: \cite{T,CT1,CT2,KOR,Chuang,CPS,FW3,PuigN,PuigScopes,JLM,Paget1,Paget,MarcusAlt,TanM,JM,Gramain,ErMi,GrOno,Muller,MarTan,MarTanW3,FoHa,FrGr,R2}.

\subsection{Branching and labeling of irreducible modules}\label{SSBrLab}
In this section we work in the ungraded setting. 
A `combinatorics-free' approach to branching rules has been implemented by Grojnowski \cite{Gr}. This approach ultimately can be used to label irreducible $H_d^\La$-modules. 

To be more precise, denote by $B(\La)$ the set of the isomorphism classes of the irreducible $H_d^\La$-modules for all $d\geq 0$. Let $D$ be an irreducible $H_d^\La$-module and $i\in I$. Grojnowski and Vazirani \cite{GV} give a direct elementary argument  to show that either $E_iD$ is zero or it has a simple socle. A similar result holds for the $F_i$'s. In level $1$, these facts go back to \cite{KBrII,KBrIII}. 
Set 
$$
\tilde E_i D:=\soc (E_i D),\quad \tilde F_i D:=\soc (E_i D)\qquad(i\in I).
$$
Thus we get operations 
$$
\tilde E_i,\tilde F_i:B(\La)\to B(\La)\sqcup\{0\}.
$$
For any $D\in B(\La)$ and $i\in I$ define 
$$
\eps_i(D):=\max\{m\mid E_i^m D\neq 0\},\quad \phi_i(D):=\max\{m\mid F_i^m D\neq 0\}.
$$
Finally, there is a natural weight function $\wt$, defined in terms of central characters, which associates to any element of $B(\La)$ a weight of the form $\La-\sum_{i\in I} c_i\al_i$, see \cite[section 12]{Gr}. A key result then is as follows

\begin{Theorem}\label{TGrojnowski} 
{\rm \cite[Theorem 14.3]{Gr}} 
The tuple $$(B(\La),\tilde E_i,\tilde F_i,\eps_i,\phi_i,\wt)$$ is the highest weight crystal associated to $V(\La)$.
\end{Theorem}

Now one can bring combinatorics in: Theorem~\ref{tdcb}(ii) gives an explicit combinatorial description of the highest weight crystal associated to $V(\La)$ in terms of restricted multipartitions. So to every  $\mu\in \RPar_d$ one can associate an irreducible module $\dot{D}(\mu)$, and then
$$
\{\dot{D}(\mu)\mid\mu\in \RPar_d\}
$$ 
is a complete irredundant set of irreducible $H_d^\La$-modules.  

Note that the definition of $\dot{D}(\mu)$ is essentially inductive: we define $\dot{D}(\mu)$ as the head of $F_i \dot{D}(\tilde e_i\mu)$ for any $i$ such that $\tilde e_i\mu\neq 0$. This description is very different in nature from the definition of $D(\mu)$ as the head of the explicit Specht module $S(\mu)$ over $H_d^\La$, given in Theorem~\ref{TIrrLab}. This explains why the natural conjecture that
\begin{equation}\label{EArGr}
\dot{D}(\mu)\cong D(\mu)\qquad(\mu\in\RPar)
\end{equation}
is not easy to prove. This has been first established by Ariki in \cite{Abranch}. Another argument (in the graded setting) is given in \cite[Theorems 5.13, 5.18]{BKllt}. 

Both arguments ultimately depend on the whole power of the theory explained in this exposition. It would be interesting to find a more elementary and direct approach for identifying the two labelings of irreducible modules.

\subsection{Extremal sequences}\label{sex}
In this section we collect some graded analogues of the results proved in \cite{BKdur}, which are quite useful in many situations. For example, Corollary~\ref{C241201} has been used in \cite{BKllt} to settle the labeling problem (\ref{EArGr}). See also \cite{Lyle} for another application. 

For an $H_d^\La$-module $M$ we denote
$$
\eps_i(M):=\max\{m\mid E_i^m M\neq 0\}.
$$
For example, $\eps_i(D(\mu)) =\eps_i(\mu)$ by Theorem~\ref{gg}(iii).

The following result provides an inductive tool for finding certain composition multiplicities. 

\begin{Lemma} \cite[Lemma 2.13]{BKdur}
\label{L291201}
Let $M\in\Rep{H_d^\La}$, $\eps=\eps_i(M)$, and $\mu\in\RPar_{d-\eps}$. 
If $$[E_i^{(\eps)} M: D(\mu)]\neq 0$$ then  
$\tilde f_i^\eps \mu\neq 0$ and 
$$[M : D(\tilde f_i^\eps \mu)]=[E_i^{(\eps)} M: D(\mu)].$$
\end{Lemma}

Given $\bi=(i_1,\dots,i_d)\in I^d$ we can gather consecutive 
equal terms to
write it in the form 
\begin{equation}
\label{E241201}
\bi=(j_1^{m_1}\dots j_r^{m_n})
\end{equation}
where $j_s\neq j_{s+1}$ for all $1\leq s<n$. For example 
$(2,2,2,1,1)=(2^3 1^2)$. Now, for an $H_d^\La$-module $M$, the tuple
(\ref{E241201}) is called {\em extremal for $M$} if 
$$
m_s=\eps_{j_s}(E_{j_{s+1}}^{m_{s+1}}\dots E_{j_n}^{m_n} M)
$$
for all $s = n, n-1, \dots, 1$. 
By definition $e(\bi)M\neq 0$ if $\bi$ is extremal for $M$. 

The main result about extremal tuples is

\begin{Theorem} \cite[Theorem 2.16]{BKdur} 
\label{T241201}
Let $$\bi=(i_1,\dots,i_d)=(j_1^{m_1}\dots j_n^{m_n})$$ 
be an extremal tuple for an irreducible $H_d^\La$-module $D(\mu)$. Then 
$$\mu=\tilde f_{i_d}\dots\tilde f_{i_1} \emptyset,$$ and 
$$\dim e(\bi)D(\mu)=[m_1]!\dots [m_n]!.$$ 
In particular, the tuple $\bi$ is not extremal for any irreducible 
$D(\nu)\not\cong D(\mu)$. 
\end{Theorem}

This result can be used to compute certain composition multiplicities as follows:

\begin{Corollary}[{\cite[Corollary 2.17]{BKdur}}]
\label{C241201}
If $\bi = (i_1,\dots,i_d)$ is an extremal sequence for $M \in \Rep{R^\La_d}$
 of the form (\ref{E241201}), then 
$\mu := \tilde f_{i_d} \cdots \tilde f_{i_1} \varnothing$
is a well-defined element of $\RPar_d$,
and 
$$
[M:D(\mu)]_q = (\qdim\ e(\bi) M)/([m_1]!\dots [m_n]!).
$$ 
\end{Corollary}

Since the $q$-characters of Specht modules are known, Corollary~\ref{C241201} is especially easy to use to compute some decomposition numbers by induction, see \cite{BKdur} for details.

\subsection{More on branching for symmetric groups}\label{SSMoreBr}

We draw the reader's attention to the fact that in level one, the original modular branching rules of \cite{KBrI}--\cite{KBrIV} and their generalizations \cite{BBr, BKtr}, give much more information than just a description of the socles of the restrictions of irreducible modules. It would be interesting to have similar more general results in higher levels. 

Here we state the results for the symmetric groups only, and in the ungraded setting. Assume that the characteristic of the ground field $F$ is $p>0$, $\La=\La_0$ and $\xi=1$. Then and $e=p$. We know from Theorem~\ref{TIrrLab} and Example~\ref{c3}(i) that the irreducible $F\Sigma_d$-modules $D(\mu)$ are labeled by the $p$-restricted partitions $\mu$ of~$d$. 

We now slightly generalize the combinatorics of section~\ref{SSCrGraph} in this particular case. 
Let $i \in I=\Z/p\Z$, and $\mathbf \sigma$ be the reduced $i$-signature of $\mu$. Recall that the $i$-node corresponding to leftmost $-$ in $\mathbf \sigma$ is called the {\em good} $i$-node of $\mu$, while the addable $i$-node corresponding to rightmost $+$ in $\mathbf \sigma$ is called the {\em cogood} $i$-node for $\mu$. Moreover, the removable $i$-nodes corresponding to  the $-$'s in $\mathbf \sigma$ are called the {\em normal} $i$-nodes of $\mu$, while the addable $i$-nodes corresponding to the $+$'s in $\mathbf \sigma$ are called the {\em conormal} $i$-nodes for $\mu$. Thus the good $i$-node is the top normal $i$-node, while the cogood $i$-node is the bottom conormal $i$-node. Also, $\eps_i(\mu)$ is the number of the normal $i$-nodes of $\mu$, while $\phi_i(\mu)$ is the number of the conormal $i$-nodes for $\mu$. 

The additional branching results not covered by Theorem~\ref{gg},  are now as follows.

\begin{Theorem}\label{t4} \cite{KBrII,KDec}
Let $\mu$ be a $p$-restricted partition of $d$, and $\nu$ be a $p$-restricted partition of $d-1$. 
\begin{itemize}
\item[(i)] $Hom_{\Sigma_{d-1}}(S(\nu),E_i D(\mu))\neq 0$ if and only if $\nu=\mu_A$ for some normal $i$-node $A$ of $\mu$. 
\item[(ii)] Suppose $\nu=\mu_A$ for some removable $i$-node $A$ of $\mu$. Then the multiplicity 
$
[E_i  D(\mu) : D(\nu)]\neq 0
$
if and only if $A$ is normal for $\mu$, in which case 
$
[E_i  D(\mu) : D(\nu)]
$  is the number of 
normal $i$-nodes weakly below $A$ 
(counting $A$ itself).
\end{itemize}
\end{Theorem}

\begin{Theorem}\label{t4'} \cite{BKtr}
Let $\mu$ be a $p$-restricted partition of $d$, and $\nu$ be a $p$-restricted partition of $d+1$. 
\begin{itemize}
\item[(i)] $Hom_{\Sigma_{d+1}}(S(\nu),F_i D(\mu))\neq 0$ if and only if $\nu=\mu^B$ for some conormal $i$-node $B$ for $\mu$. 
\item[(ii)] Suppose $\nu=\mu^B$ for some addable $i$-node $B$ for $\mu$. Then the multiplicity 
$
[F_i  D(\mu) : D(\nu)]\neq 0
$
if and only if $B$ is conormal for $\mu$, in which case 
$
[F_i  D(\mu) : D(\nu)]
$  is the number of 
conormal $i$-nodes weakly above $B$ 
(counting $B$ itself).
\end{itemize}
\end{Theorem}

The additional branching information provided by Theorems~\ref{t4}, \ref{t4'} turns out very useful in applications, see for example \cite{BarKZ,BarK,BKIrrRes,El,KZ,JM, Shch, FT,Wi,JaWi,KuNe}. 

\subsection{Mullineux Involution}\label{SSMull}

The algebra $H_d$ has an automotphism 
$$
\si:H_d\to H_d,\ T_r\mapsto -T_r+\xi-1\qquad(1\leq r<d).
$$
Twisting an irreducible $H_d$-module $D(\mu)$  yields an irreducible module $D(\mu)^\si$. For the symmetric groups we simply have 
$$
D(\mu)^\si\cong D(\mu)\otimes \sign,
$$
where $\sign$ is the one-dimensional sign representation of $\Sigma_d$. 

The involution $\si$ arises in many natural situations. For example, one needs to deal with it when studying representation theory of alternating groups. Another illustration comes from Remark~\ref{RJamD}. 


The next result provides an explicit combinatorial algorithm for computing the involution $\si$ in terms of the combinatorics of good nodes. 

\begin{Theorem}\label{TMul}%
 {\rm \cite{KBrIII}} 
Let $\mu$ be an $e$-restricted partition of $d$. Write 
\begin{equation}\label{EPres}
\mu=\tilde f_{i_1}\tilde f_{i_2}\dots \tilde f_{i_d}\emptyset
\end{equation}
for some $i_1,i_2,\dots,i_d\in I$. Then $D(\mu)^\si\cong D(\nu)$, where
\begin{equation}\label{EPresNu}
\nu=\tilde f_{-i_1}\tilde f_{-i_2}\dots \tilde f_{-i_d}\emptyset.
\end{equation}
\end{Theorem}

It is implicit in the theorem above that there always exists a presentation of $\mu$ in the form (\ref{EPres}) starting with the empty partition $\emptyset$ and that the formula (\ref{EPresNu}) determines a well-defined $e$-restricted partition. Indeed, to find a presentation of $\mu$ in the form (\ref{EPres}) one should consecutively remove good nodes from $\mu$ and record their residues in order of removal: $i_1,i_2,\dots,i_d$. Then one should build $\nu$ by consecutively adding  good nodes of residues $-i_d,\dots,-i_2,-i_1$. 

Originally, Mullineux \cite{Mul} conjectured a different description of the involution $D(\mu)\mapsto D(\mu)^\si$. It has been shown in \cite{FK} that that description is equivalent to the one given in Theorem~\ref{TMul}. Questions related to the Mullineux involution received quite a lot of attention. We refer the reader to \cite{BBr,BO,Xu,BrKuj,LTV,FMull,HuMull,Paget,BOX,BOMJS,JacLec,FRunMul,ShuWang} as well as 
\cite[Remarks 3.17, 3.18]{BKllt} and \cite[section~3.7]{BKariki} for further information.

\subsection{Higher level Schur-Weyl duality, $W$-algebras, and category $\mathcal O$}\label{SSHLSWD}
Throughout this section we assume that $\xi=1$ and $F=\C$. In this case, there is quite a different approach to representation theory of $H_d^\La$ based on a generalization of Schur-Weyl duality. This duality connects representation theory of the cyclotomic Hecke algebra $H_d^\La$ with that of a {\em finite $W$-algebra} and ultimately with a {\em parabolic category $\mathcal O$}. We review the main features of that theory here referring the reader to \cite{BKschur,BKariki} for details.

Using shifts of the elements $X_r\in H_d^\La$ by the same scalar $k_1$, we may assume without loss of generality that in (\ref{EKappa}) we have $0=k_1\geq k_2\geq\dots\geq k_l$. Let $n\geq -k_l$, and define 
$$
q_m:=n+k_m\qquad(1\leq m\leq l). 
$$
Then $(q_1,\dots,q_l)$ is a partition. Let $q_1+\dots+q_l=N$, and let $\lambda = (p_1 \leq \cdots \leq p_n)$ be a partition of $N$ transposed to $(q_1,\dots,q_l)$.

Unlike in the main body of the paper, it is convenient to identify $\lambda$ with its Young diagram in an unusual way, numbering
rows by $1,2,\dots,n$ from top to bottom
and columns 
by $1,2,\dots,l$ from left to right, 
so that 
there are $p_i$ boxes in the $i$th row and $q_j$ boxes in the $j$th column.
For example, for $\La=2 \Lambda_0 + \Lambda_{-1}+\Lambda_{-2}$ and $n=3$, we have  $\lambda = (p_1,p_2,p_3) = (2,3,4)$, 
$(q_1,q_2,q_3,q_4)=(3,3,2,1)$, 
$N= 9$.
The Young diagram is
$$
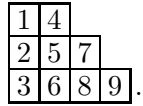
\Diagram{1&4\cr2&5&7\cr3&6&8&9\cr}
\begin{picture}(0,0)
\put(0,0){\makebox(7,-25){.}}
\end{picture}
$$
We always number the boxes of the diagram $1,2,\dots,N$ down columns starting
from the first column,
and write $\row(i)$ and 
$\col(i)$ for the row and column
numbers of the $i$th box. 
This identifies the boxes
with 
the standard basis $v_1,\dots, v_N$ of the natural $\mathfrak{gl}_N(\C)$-module
$V$. 
Define
 $e \in \mathfrak{gl}_N(\C)$ to be the nilpotent matrix of Jordan type $\lambda$
which maps the basis vector corresponding 
the $i$th box to the one immediately to its left, or to zero if there
is no such box;
in our example, $e = e_{1,4} + e_{2,5}+e_{5,7}+e_{3,6}+e_{6,8}+e_{8,9}$.
Finally, define a $\Z$-grading 
$$\mathfrak{gl}_N(\C) = \bigoplus_{r \in \Z} \mathfrak{gl}_N(\C)_r$$ 
on $\mathfrak{gl}_N(\C)$ by
declaring that $e_{i,j}$ is of degree $\col(j)-\col(i)$
for each $i,j=1,\dots,N$, and set
$$\mathfrak{p} = \bigoplus_{r \geq 0} \mathfrak{gl}_N(\C)_r,\quad
\mathfrak{h} = \mathfrak{gl}_N(\C)_0,\quad
\mathfrak{m} = \bigoplus_{r < 0} \mathfrak{gl}_N(\C)_r.
$$
The 
centralizer $\mathfrak{z}_e$ of $e$ in $\mathfrak{gl}_N(\C)$
is a graded subalgebra of $\mathfrak{p}$.
Hence its universal enveloping algebra $U(\mathfrak{z}_e)$ is a graded
subalgebra of $U(\mathfrak{p})$.

We now define the {\em finite $W$-algebra} $W(\lambda)$
associated to the partition $\lambda$, 
following \cite{BKrep};
see also \cite{P, GG}.
Let $\eta:U(\mathfrak{p}) \rightarrow U(\mathfrak{p})$
be the algebra automorphism defined by
\begin{align*}
\eta(e_{i,j}) &= 
e_{i,j} + \delta_{i,j} (n-q_{\col(j)} - q_{\col(j)+1} - \cdots - q_l)
\end{align*}
for each $e_{i,j} \in \mathfrak{p}$.
Let $I_\chi$ be the kernel of the homomorphism $\chi:U(\mathfrak{m}) 
\rightarrow \C$ defined by
$x \mapsto (x,e)$ for all $x \in \mathfrak{m}$, where $(.,.)$ is the trace form
on $\mathfrak{gl}_N(\C)$. 
Then 
$W(\lambda)$
is defined to be the following subalgebra of $U(\mathfrak{p})$:
\begin{align*}
W(\lambda) 
&= \{u \in U(\mathfrak{p})\mid
[x, \eta(u)] \in U(\mathfrak{g}) I_\chi \text{ for all }x \in \mathfrak{m}\}.
\end{align*}

We now introduce a structure of a 
$(W(\lambda), H_d^\La)$-bimodule on 
$V^{\otimes d}$.
The left action of $W(\lambda)$ on $V^{\otimes d}$ is simply the restriction
of the natural action of $U(\mathfrak{p})$.
To define the right action of $H_d^\La$,
let 
$\Sigma_d$ act on the right 
by place permutation as usual.
Let $X_1$ act as the endomorphism
$$
\bigg(e + \sum_{j=1}^N (q_{\col(j)}-n) e_{j,j} \bigg)
\otimes 1^{\otimes (d-1)}
-
\sum_{k=2}^d \sum_{\substack{i,j=1 \\ \col(i) < \col(j)}}^N
e_{i,j} \otimes 1^{\otimes(k-2)} \otimes e_{j,i} \otimes 1^{\otimes (d-k)}.
$$
It turns out \cite{BKschur} that this extends uniquely to make 
$V^{\otimes d}$ into a 
 $(W(\lambda),H_d^\Lambda)$-bimodule.

We have now defined a homomorphism $\Phi_d$ and an antihomomorphism
$\Psi_d$
$$ 
W(\lambda) \stackrel{\Phi_d}{\longrightarrow} \End_{\C}(V^{\otimes d})
\stackrel{\Psi_d}{\longleftarrow} H_d^\Lambda.
$$
Let $W_d(\lambda)$ denote the image of 
$\Phi_d$.
This finite dimensional algebra is
a natural
analogue of the classical {\em Schur algebra} for higher levels.
Let $H_d(\lambda)$ denote
the image of the homomorphism
$\Psi_d:H_d(\Lambda) \rightarrow \End_{\C}(V^{\otimes d})^{\op}$,
so that $V^{\otimes d}$ is also a $(W_d(\lambda), H_d(\lambda))$-bimodule.
Actually, if at least $d$ parts of $\lambda$ are equal to $l$,
then the map $\Psi_d$ is injective so $H_d(\lambda) = H_d^\Lambda$.
In general $H_d(\lambda)$ is 
a sum of certain blocks of $H_d^\Lambda$.

\begin{Theorem} \label{TA} \cite[Theorem A]{BKschur}
The maps $\Phi_d$ and $\Psi_d$ satisfy the double centralizer property, i.e.
\begin{align*}
W_d(\lambda)=\End_{H_d(\lambda)}(V^{\otimes d}),
\qquad&\End_{W_d(\lambda)}(V^{\otimes d})^{\op} = H_d(\lambda).
\end{align*}
Moreover, 
the functor
\begin{align*}
\Hom_{W_d(\lambda)}(V^{\otimes d}, ?):
&\,W_d(\lambda)\text{\rm-mod}
\rightarrow H_d(\lambda)\text{\rm-mod}
\end{align*}
is an equivalence of categories.
\end{Theorem}

Let $P$ denote the
module $U(\mathfrak{gl}_N(\C)) \otimes_{U(\mathfrak{p})} 
\C_{-\rho}$
induced from the one dimensional $\mathfrak{p}$-module $\C_{-\rho}$
on which each $e_{i,j} \in \mathfrak{p}$ acts as 
$\delta_{i,j}(q_1+q_2+\cdots+q_{\col(j)} - n)$.
This is an irreducible projective 
module in  a parabolic category $\mathcal{O}$ for $\mathfrak{gl}_N(\C)$ with respect to $\mathfrak{p}$.
Let $\mathcal{O}^d(\lambda)$ denote the 
Serre subcategory of the parabolic category $\mathcal{O}$
generated by the module
$P \otimes V^{\otimes d}$.
We note that $\mathcal{O}^d(\lambda)$  is a sum of certain integral blocks of the  parabolic category
$\mathcal{O}$, and every integral block is equivalent to a block of
$\mathcal{O}^d(\lambda)$ for sufficiently large $d$.

Moreover, 
the module $P \otimes V^{\otimes d}$
is a {\em self-dual projective module} in $\mathcal{O}^d(\lambda)$, and every
self-dual projective indecomposable module in $\mathcal{O}^d(\lambda)$
is a summand of $P \otimes V^{\otimes d}$. 
Applying the construction from \cite[$\S$2.2]{AS}, we can view
$P \otimes V^{\otimes d}$ as a $(\mathfrak{gl}_N(\C), H_d^\aff)$-bimodule.
It turns out that the
 right action of $H_d^\aff$ on
$P \otimes V^{\otimes d}$ 
factors through the quotient $H_d(\lambda)$ of $H_d$ to make
$P \otimes V^{\otimes d}$ into a faithful right 
$H_d(\lambda)$-module,
i.e. $H_d(\lambda) \hookrightarrow \End_{\C}(P \otimes V^{\otimes d})^{\op}$.

\begin{Theorem} \label{TB} \cite[Theorem B]{BKschur} We have 
$\End_{\mathfrak{gl}_N(\C)}(P \otimes V^{\otimes d})^{\operatorname{op}}
= H_d(\lambda).
$
\end{Theorem}

The link between Theorems \ref{TA}  and \ref{TB} is provided by the
{\em Whittaker functor} 
$$
\mathbb{V}: \mathcal{O}^d(\lambda)
\rightarrow W_d(\lambda)\text{-mod}
$$ 
introduced originally by
Kostant and Lynch \cite{Lynch} and studied recently in \cite[$\S$8.5]{BKrep} (cf. \cite{Soergel}):
the $W_d(\lambda)$-module
$V^{\otimes d}$ from above is isomorphic to 
$\V(P \otimes V^{\otimes d})$.
We actually show in \cite{BKschur}
that
$\Hom_{\mathfrak{g}}(P \otimes V^{\otimes d}, ?)
\cong
\Hom_{W_d(\lambda)}(V^{\otimes d}, ?) \circ \V$,
i.e. the following diagram of functors commutes up to isomorphism:
$$
\begin{CD}
&\:\:\mathcal{O}^d(\lambda) &\\
\\\\
W_d(\lambda)\text{-mod}@>\sim > \Hom_{W_d(\lambda)}(V^{\otimes d}, ?)> &
H_d(\lambda)\text{-mod.}
\end{CD}
\begin{picture}(0,0)
\put(-147,-8){\makebox(0,0){$\swarrow$}}
\put(-139,12){\makebox(0,0){$\scriptstyle\mathbb{V}$}}
\put(-144,-5){\line(1,1){33}}
\put(-48,-8){\makebox(0,0){$\searrow$}}
\put(-51,-5){\line(-1,1){33}}
\put(-26,12){\makebox(0,0){$\scriptstyle\Hom_{\mathfrak{g}}(P \otimes V^{\otimes d}, ?)$}}
\end{picture}
$$
The categories $W_d(\lambda)\text{-mod}$ and $H_d(\lambda)\text{-mod}$ 
thus give two
different realizations of a
natural quotient of the category $\mathcal{O}^d(\lambda)$
in the general sense of \cite[$\S$III.1]{Gab}, 
the respective quotient functors being
the Whittaker functor $\mathbb{V}$ and the functor
$\Hom_{\mathfrak{g}}(P \otimes V^{\otimes d}, ?)$.

In many circumstances, the Whittaker functor 
turns out to be easier to 
work with than the functor $\Hom_{\mathfrak{g}}(P \otimes V^{\otimes d}, ?)$,
so this point of view
facilitates various other important computations
regarding the relationship between $\mathcal{O}^d(\lambda)$
and $H_d(\lambda)$-mod.
For example, 
we use it to identify
the images of arbitrary projective indecomposable modules in $\mathcal{O}^d(\lambda)$ with  the indecomposable summands
of the degenerate analogues of the
permutation modules introduced by Dipper, James and 
Mathas \cite{DJM}.
Also, we identify the images of parabolic Verma modules with {\em Specht modules}, thus recovering formulae for 
the latter's composition multiplicities directly from the Kazhdan-Lusztig
conjecture for $\mathfrak{gl}_N(\C)$.
We remark that the degenerate analogue of Ariki's categorification theorem \cite{Ariki} follows
as an easy consequence of these results, as is explained 
in  \cite{BKariki}.

\subsection{Projective representations}\label{SSProj}
An analogue of Grojnowski's theory \cite{Gr}  for {\em projective representations} of symmetric and alternating groups has been established in \cite{BKHCl}. Alternatively, the  theory developed in  \cite{BKHCl} can be thought of as a generalization of Grojnowski's work \cite{Gr} from the symmetric groups $\Sigma_d$ to their Schur double covers $\hat \Sigma_d$. 

For the double covers, the role of the Kac-Moody algebra $\g$ is played by the twisted Kac-Moody algebra of type $A_{2e}^{(2)}$. Both degenerate and non-degenerate cases are treated in \cite{BKHCl}, as well as higher level analogues. Level $1$ theory has been treated before in \cite{BKSerg} using Schur-Weyl type duality due to Sergeev \cite{Se85}. In \cite{BKSerg} we also suggest a reasonable notion of `Specht modules' for $\hat\Sigma_d$. See however \cite{BKReg}, where we suggest that certain complexes of modules might be a better choice for `Specht modules'. 

We point out that although socle branching rules and the corresponding classification of irreducible modules in spirit of section~\ref{SSBrLab} work for $\hat\Sigma_d$ just like for $\Sigma_d$ in \cite{Gr}, the connections with geometry and canonical bases are not available at the moment. We also note that the analogue of the labeling problem (\ref{EArGr}) for $\hat\Sigma_d$ remains open. 

For various other recent developments in the theory of projective representations of symmetric groups we refer the reader to \cite{Se99,Ya1,Ya2,Wang,Jos,JN,Naz1,Naz2,Ph,JW,Jo1,Jo2,BKdet,BKReg,KS,Ke}.

\subsection{Problems on symmetric groups related to Aschbacher-Scott program}\label{SSASP}
In this section we indicate some applications of representation theory of symmetric groups to structure questions in finite group theory. We briefly describe the set up, referring the interested reader to \cite{KlLi}. 

Let  $C$ be 
a finite classical group, and $H<C$ be a maximal proper subgroup. 
Aschbacher's theorem \cite{A}, which starts the {\em Aschbacher-Scott program}, (cf. also \cite{Sc}) claims that  then 
\begin{equation}\label{EM}
\textstyle H \in \SCL\cup \bigcup^{8}_{i=1}\CL_{i},
\end{equation}
where $\CL_{i}$, $i = 1, \ldots ,8$, are
collections of certain explicit natural subgroups of $C$, and $\SCL$ is the collection
of all almost quasi-simple subgroups that act absolutely irreducibly on the natural module for the classical group $C$. 

However, the converse to Aschbacher's theorem does not hold in general. 
So, to understand maximal subgroups of finite classical groups, one needs to determine
when a subgroup $H$ as in (\ref{EM}) is actually maximal in
$C$. For $H \in \cup^{8}_{i=1}\CL_{i}$, this has been done by
Kleidman and Liebeck in \cite{KlLi}. 
Let $H \in \SCL$. If $H$ is {\em not} maximal then $H < G < C$ for a certain maximal subgroup  $G$ in $C$. The most difficult case to handle is when $G \in \SCL$ as well. This motivates the following problem. 

\begin{Problem}\label{Prestr}
Classify all triples $(G,V,H)$ where 
$G$ is an almost quasi-simple finite group, $V$ 
is an $FG$-module of dimension greater than one, and $H$ is a proper subgroup of $G$ such that 
the restriction $\res^G_H V$ is irreducible.
\end{Problem}

Under the assumption that the characteristic of the ground field is greater than $3$, Problem~\ref{Prestr} has been solved for $G$ of {\em alternating type}, i.e. $G=A_d,\Sigma_d,\hat A_d$ or $\hat \Sigma_d$, see \cite{Saxl,KlWa,KJS,BKIrrRes,KS2,KT}, and  partial results are available even in the cases $\ell=3$ and $2$, see e.g. \cite{KS1}. A variety of techniques has been used, but modular branching rules played an important  role. 

Another part of the Aschbacher-Scott program yields the following problem.

\begin{Problem}\label{Ptens}
Let  
$G$ be an almost quasi-simple finite group. Classify all pairs $(V,W)$ 
of irreducible $FG$-modules of dimension greater than one such that $V\otimes W$ is irreducible.
\end{Problem}

Under the assumption that the characteristic of the ground field is greater than $2$, Problem~\ref{Ptens} has been mainly solved for $G$ of alternating type, see \cite{BessK0,BessK0Spin,BessK,BessKAlt,KT}, and  partial results are available even in the case $\ell=2$, see e.g. \cite{GowK,GrJa}.


\begin{thebibliography}{MO}

\bibitem
{AJS}
H.~H. Andersen, J.~C. Jantzen, and W. Soergel, Representations of quantum groups at a $p$th root of unity and of semisimple groups in characteristic $p$: independence of $p$, {\em Ast\'erisque} {\bf 220} (1994), 321 pp.


\bibitem
{AS}
T. Arakawa and T. Suzuki,
Duality between $\mathfrak{sl}_n(\mathbb{C})$ and the degenerate affine Hecke algebra, {\em  J. Algebra} {\bf  209} (1998), 288--304; {\tt q-alg/9710037}.



\bibitem
{Ariki}
S. Ariki, On the decomposition numbers of the Hecke algebra of $G(m,1,n)$, {\em J. Math. Kyoto Univ.} {\bf 36} (1996), 789--808. 

\bibitem
{Aclass}
S. Ariki,
On the classification of simple modules for the cyclotomic Hecke
algebra of type $G(m,1,n)$ and Kleshchev multipartitions,
{\em Osaka J. Math.} {\bf 38} (2001), 827--837.

\bibitem
{Abook}
S. Ariki, {\em Representations of Quantum Algebras and Combinatorics of Young 
Tableaux}, University Lecture Series 26, American Mathematical Society, Providence, RI, 2002. 

\bibitem
{Abranch}
S. Ariki,
 Proof of the modular branching rule for cyclotomic Hecke algebras, 
{\em J. Algebra} {\bf 306} (2006), 290--300.

\bibitem
{ArikiJac}
S. Ariki and N. Jacon, Dipper-James-Murphy's conjecture for Hecke algebras of type $B$, {\tt arXiv:math/0703447}.

 
\bibitem
{AJL}
S. Ariki, N. Jacon and C. Lecouvey, The modular branching rule for affine Hecke algebras of type A, {\tt arXiv:0808.3915}


\bibitem
{AK}
S. Ariki and K. Koike, 
A Hecke algebra of $(\Z/r \Z)\wr S\sb n$ 
and construction of its irreducible representations,
{\em Advances Math.} {\bf 106} (1994), 216--243.

\bibitem
{ArKrTs}
S. Ariki, V. Kreiman and S. Tsuchioka, On the tensor product of two basic representations 
of $U_v(\widehat{\sl})$, {\em Adv. Math.}, {\bf 218} (2008), 28Ð86. 

\bibitem
{AM}
S. Ariki and A. Mathas, The number of simple modules of the Hecke algebra of type $G(r,1,n)$, {\em Math. Z.} {\bf 233} (2000), 601--623.


\bibitem
{AMR}
S. Ariki, A. Mathas and H. Rui, Cyclotomic Nazarov-Wenzl algebras,  {\em Nagoya Math. J.} {\bf 182} (2006), 47--134. 


\bibitem
{A}
  M. Aschbacher, On the maximal subgroups of the finite classical groups,
{\em Invent. Math.} {\bf 76} $(1984)$, 469--514.

\bibitem
{BarK}
A.A. Baranov and  A.S. Kleshchev, Maximal ideals in modular group algebras of the finitary symmetric and alternating groups, {\em Trans. Amer. Math. Soc.} {\bf 351} (1999), 595--617. 

\bibitem
{BarKZ}
A.A. Baranov, A.S. Kleshchev and A.E. Zalesskii, Asymptotic results on modular representations of symmetric groups and almost simple modular group algebras, {\em J. Algebra} {\bf 219} (1999), 506--530.

\bibitem
{BOMJS}
C. Bessenrodt and J.B. Olsson,  Residue symbols and Jantzen-Seitz partitions, {\em J. Combin. Theory Ser. A} {\bf 81} (1998), 201--230. 

\bibitem
{BO}
C. Bessenrodt and J.B. Olsson, On residue symbols and the Mullineux conjecture, {\em J. Algebraic Combin.} {\bf 7} (1998), 227--251.

\bibitem
{BOBr}
C. Bessenrodt and J.B. Olsson, Branching of modular representations of the alternating groups,  {\em J. Algebra} {\bf 209} (1998), 143--174.

\bibitem
{BOX}
C. Bessenrodt, J.B. Olsson and M. Xu, On properties of the Mullineux map with an application to Schur modules, {\em Math. Proc. Cambridge Philos. Soc.} {\bf 126} (1999),   443--459. 

\bibitem
{BessK0} C. Bessenrodt and A. Kleshchev, On Kronecker products of complex representations of the symmetric and alternating groups, {\em Pacific J. Math.}  {\bf 190} (1999), 201--223. 

\bibitem
{BessK} C. Bessenrodt and A. Kleshchev, On tensor products of modular representations of symmetric groups, {\em Bull. Lond. Math. Soc.} {\bf 32} (2000), 292--296. 

\bibitem
{BessKAlt} C. Bessenrodt and A. Kleshchev, Irreducible tensor products over alternating groups, {\em J. Algebra} {\bf 228} (2000), 536--550.

\bibitem
{BessK0Spin} C. Bessenrodt and A. Kleshchev, On Kronecker products of spin characters of the double covers of the symmetric groups, {\em Pacific J. Math.} {\bf 198} (2001), 295--305. 



\bibitem
{BGS}
A. Beilinson, V. Ginzburg and W. Soergel,
Koszul duality patterns in representation theory,
{\em J. Amer. Math. Soc.} {\bf 9} (1996), 473--527.

\bibitem
{Br1} 
M. Brou\'e, Equivalences of blocks of group algebras, in {\em  Finite-dimensional algebras and related topics}, NATO ASI Ser. 424, Kluwer, 1994, 1--26.

\bibitem
{BM} 
M. Brou\'e and G. Malle,  Zyklotomische Heckealgebren, {\em  Ast\'erisque} {\bf 212} (1993), 119--189. 

\bibitem
{BBr}
J. Brundan, Modular branching rules and the {M}ullineux map for {H}ecke
  algebras of type {$\mathbf A$}, 
{\em Proc. London Math. Soc.} {\bf 77} (1998), 551--581.

\bibitem
{cyclo}
J. Brundan, Centers of degenerate cyclotomic Hecke algebras and parabolic category $\mathcal O$, {\em Represent. Theory} {\bf 12} (2008), 236-259.

\bibitem
{BKtr}
J. Brundan and A. Kleshchev,  On translation functors for general linear and symmetric groups, {\em Proc. London Math. Soc. (3)} {\bf 80} (2000), 75--106.

\bibitem
{BKIrrRes}
J. Brundan and A. Kleshchev,   Representations of the symmetric group which are irreducible over subgroups, {\em J. Reine Angew. Math.} {\bf 530} (2001), 145--190. 

\bibitem
{BKHCl}
J. Brundan and A. Kleshchev, Hecke-Clifford superalgebras, crystals of type $A_{2\ell}^{(2)}$ and modular branching rules for $\hat S_n$, {\em Represent. Theory} {\bf 5} (2001), 317--403. 

\bibitem
{BKSerg}
J. Brundan and A. Kleshchev,
Projective representations of symmetric groups via Sergeev duality, {\em Math. Z.} {\bf 239} (2002), 27--68.

\bibitem
{BKdet}
J. Brundan and A. Kleshchev, Cartan determinants and Shapovalov forms, {\em Math. Ann.} {\bf 324} (2002), 431--449.

\bibitem
{BKdur}
J. Brundan and A. Kleshchev,
Representation theory of symmetric groups and their double covers,
in: {\em Groups, Combinatorics and Geometry (Durham, 2001)},
pp.31--53,
World Scientific Publishing, River Edge, NJ, 2003.

\bibitem
{BKReg}
J. Brundan and A. Kleshchev, James' regularization theorem for double covers of symmetric groups, {\em J. Algebra} {\bf 306} (2006), 128--137. 

\bibitem
{BKrep}
J. Brundan and A. Kleshchev,
Representations of shifted Yangians and finite $W$-algebras,
{\em Mem. Amer. Math. Soc.}
{\bf 196} (2008), no. 918, 107 pp..

\bibitem
{BKschur}
J. Brundan and A. Kleshchev,
Schur-Weyl duality for higher levels,
{\em Selecta Math. (N.S.)} {\bf 14} (2008), 1--57.

\bibitem
{BKyoung}
J. Brundan and A. Kleshchev,
Blocks of cyclotomic Hecke algebras and Khovanov-Lauda algebras, {\em Invent. Math.}, to appear; 
{\tt arXiv:0808.2032}.

\bibitem
{BKariki}
J. Brundan and A. Kleshchev,
The degenerate analogue of Ariki's categorification theorem, {\em Math. Z.}, to appear; 
 {\tt arXiv:0901.0057}.


\bibitem
{BKllt}
J. Brundan and A. Kleshchev,
Graded decomposition numbers for cyclotomic Hecke algebras, {\em Adv. in Math.}, to appear; 
 {\tt arXiv:0901.4450}.

\bibitem
{BKS}
J. Brundan, A. Kleshchev and I. Suprunenko, Semisimple restrictions from ${\rm GL}(n)$ to ${\rm GL}(n-1)$, {\em J. Reine Angew. Math.} {\bf 500} (1998), 83--112.

\bibitem
{BKW}
J. Brundan, A. Kleshchev and W. Wang,
Graded Specht modules,
 {\tt arXiv:0901.0218}.

\bibitem
{BrKuj}
J. Brundan and J. Kujawa, A new proof of the Mullineux conjecture, {\em J. Algebraic Combin.} {\bf 18} (2003), 13--39.

\bibitem
{Ch}
I. Cherednik, A new interpretation of Gelfand-Tzetlin bases, 
{\em Duke Math. J.} {\bf 54} (1987), 563--577. 

\bibitem
{CG}
N. Chriss and V. Ginzburg, {\em Representation Theory and Complex Geometry},
Birkh\"auser, 1997.

\bibitem
{Chuang}
J. Chuang, The derived categories of some blocks of symmetric groups and a conjecture of Brou\'e, {\em J. Algebra} {\bf 217} (1999), 114--155. 

\bibitem
{CK}
J. Chuang and R. Kessar, Symmetric groups, wreath products, Morita equivalences, and BrouŽ's abelian defect group conjecture, {\em Bull. London Math. Soc.} {\bf 34} (2002),  174--184.

\bibitem
{CR}
J. Chuang and R. Rouquier,
Derived equivalences for symmetric groups and $\mathfrak{sl}_2$-categorification, {\em Ann. of Math.} {\bf 167} (2008), 245--298.

\bibitem
{CT1}
J. Chuang and Kai Meng Tan, Filtrations in Rouquier blocks of symmetric groups and Schur algebras, {\em Proc. London Math. Soc. (3)} {\bf 86} (2003), 685--706. 


\bibitem
{CT2}
J. Chuang and Kai Meng Tan,  Some canonical basis vectors in the basic $U\sb q(\widehat{\mathfrak sl}\sb n)$-module, {\em J. Algebra} {\bf 248} (2002), 765--779. 

\bibitem
{CPS}
E. Cline, B. Parshall and L. Scott, 
Stratifying endomorphism algebras, 
{\em Mem. Amer. Math. Soc.} {\bf 124} (1996), no. 591, viii+119 pp.

\bibitem
{DG}
P. Diaconis and C. Green, Applications of Murphys elements, {\em Stanford University Technical Report}, (1989), no. 335.

\bibitem
{DJ}
R. Dipper and G. D. James, Representations of Hecke algebras of  
general linear groups, {\em Proc. London Math. Soc.} {\bf 52} (1986), 20--52. 

\bibitem
{DJBI}
R.Dipper and G.D. James, Blocks and idempotents of Hecke algebras of general linear groups, {\em Proc. London Math. Soc. (3)}, {\bf 54} (1987), 57--82. 

\bibitem
{DJM}
R. Dipper, G. D. James and A. Mathas, Cyclotomic $q$-Schur algebras, {\em Math. Z.} {\bf 229} (1998), 385--416. 

\bibitem
{DuRui}
J. Du and H. Rui, Specht modules for Ariki-Koike algebras, {\em Comm. Algebra} {\bf 29}  (2001), 4701--4719. 

\bibitem
{El}
H. Ellers, Searching for more general weight conjectures, using the symmetric group as an example, {\em J. Algebra} {\bf 225} (2000),  602--629. 

\bibitem
{EM}
H. Ellers and J. Murray, 
Branching rules for Specht modules, 
{\em J. Algebra} {\bf 307} (2007),  278--286.


\bibitem
{ErMi} 
K. Erdmann, G.O. Michler, Blocks for symmetric groups and their covering groups and quadratic forms, {\em BeitrŠge Algebra Geom.} {\bf 37} (1996), 103--118.

\bibitem
{FayersWeight}
M. Fayers, Weights of multipartitions and representations of 
Ariki-Koike algebras, {\em Adv. in Math.} {\bf 206} (2006), 112--144. 

\bibitem
{F}
M. Fayers, 
An extension of James' conjecture,
{\em Int. Math. Res. Not. IMRN} {\bf 10} (2007), 24 pp.

\bibitem
{FW4}
M. Fayers,  James's Conjecture holds for weight four blocks of Iwahori-Hecke algebras. {\em J. Algebra} {\bf 317} (2007), 593--633. 

\bibitem
{FW3} 
M. Fayers,  Decomposition numbers for weight three blocks of symmetric groups and Iwahori-Hecke algebras, {\em Trans. Amer. Math. Soc.} {\bf 360} (2008), 1341--1376.



\bibitem
{FMull}
M. Fayers,   Regularisation and the Mullineux map, {\em Electron. J. Combin.} {\bf 15}  (2008), no. 1, Research Paper 142, 15 pp.

\bibitem
{FRunMul}
M. Fayers, General runner removal and the Mullineux map, {\tt  arXiv:0712.2390}.


\bibitem
{FT}
M. Fayers and Kai Meng Tan, 
The ordinary quiver of a weight three block of the symmetric group is bipartite, {\em Adv. Math.} {\bf 209} (2007), 69--98. 


\bibitem
{FoHa} 
P. Fong and M.E. Harris, On perfect isometries and isotypies in alternating groups, {\em Trans. Amer. Math. Soc.} {\bf 349} (1997), 3469--3516.


\bibitem
{FK} 
B. Ford and A. Kleshchev, A proof of the Mullineux conjecture, {\em Math. Z.}  {\bf 226} (1997), 267--308.

\bibitem
{FrGr}
A.R. Francis and J.J. Graham, Centres of Hecke algebras: the Dipper-James conjecture,  {\em J. Algebra} {\bf 306} (2006), 244--267. 


\bibitem
{Gab}
P. Gabriel,
Des cat\'egories Ab\'eliennes,
{\em Bull. Sci. Math. France} {\bf 90} (1962), 323--448.




\bibitem
{GG}
W. L. Gan and V. Ginzburg,
Quantization of Slodowy slices,
{\em Internat. Math. Res. Notices} {\bf 5} (2002), 243--255.



\bibitem
{Geck}
M. Geck, 
Representations of Hecke algebras at roots of unity, Ast\'erisque {\bf 252} 
(1998), 33--55.

\bibitem
{GeckKLJ} 
M. Geck, Kazhdan-Lusztig cells, $q$-Schur algebras and James' conjecture, {\em J. London Math. Soc. (2)} {\bf 63} (2001), 336--352. 

\bibitem
{GMu}
M. Geck and J. Mueller, James' Conjecture for Hecke algebras of exceptional type, I, {\em  J. Algebra} {\bf 321} (2009), 3274--3298. 

\bibitem
{GowK}
R. Gow and A. Kleshchev, 
Connections between the representations of the symmetric group and the symplectic group in characteristic $2$,  
{\em J. Algebra} {\bf 221} (1999), 60--89. 

\bibitem
{GrJa}
 J. J. Graham and G. James, On a conjecture of Gow and Kleshchev concerning tensor products, {\em J. Algebra} {\bf 227} (2000), 767--782.

\bibitem
{GL}
 J. J. Graham and G. I. Lehrer,
Cellular algebras, {\em Invent. Math.} {\bf 123} (1996), 1--34.

\bibitem
{Gramain}
J.-B. Gramain,  On defect groups for generalized blocks of the symmetric group, {\em  J. Lond. Math. Soc. (2)} {\bf 78} (2008), 155--171.

\bibitem
{GrOno}
A. Granville and K. Ono, Defect zero $p$-blocks for finite simple groups, {\em  Trans. Amer. Math. Soc.} {\bf 348} (1996), 331--347. 

\bibitem
{G1}
I. Grojnowski,
Representations of affine Hecke algebras (and affine quantum $GL_n$)
at roots of unity, {\em Internat. Math. Res. Not.} {\bf 5} (1994), 215--217.

\bibitem
{Gr}
I. Grojnowski, Affine $\mathfrak{sl}_p$ controls the representation theory of the symmetric group and related Hecke algebras, {\tt arXiv:math.RT/9907129}. 


\bibitem
{GV}
I. Grojnowski and M. Vazirani, Strong multiplicity one theorem for affine
Hecke algebras of type $A$, {\em Transf. Groups.} {\bf 6} (2001), 143--155. 



\bibitem
{Hayashi}
T. Hayashi,
$q$-Analogues of Clifford and Weyl algebras:
spinor and oscillator representations of quantum enveloping algebras,
{\em Comm. Math. Phys.} {\bf 127} (1990), 129--144.

\bibitem
{Hill} 
D. Hill, Elementary divisors of the Shapovalov form on the basic representation of Kac-Moody algebras, {\em J. Algebra} {\bf 319} (2008), 5208--5246.

\bibitem
{H}
P. Hoefsmit, {\em Representations of Hecke Algebras of Finite Groups with ${\rm BN}$-Pairs of Classical Type}, Ph.D. thesis, University of British Columbia, 1974.

\bibitem
{Hu}
J. Hu, Branching rules for Hecke algebras of type $D\sb n$, {\em Math. Nachr.} {\bf 280}  (2007), 93--104. 

\bibitem
{HuMull}
J. Hu, Mullineux involution and twisted affine Lie algebras, {\em J. Algebra} {\bf 304} (2006),  557--576.

\bibitem{HuMathas}
J. Hu and A. Mathas, Graded cellular bases for the cyclotomic Khovanov-Lauda-Rouquier algebras of type A; {\tt  arXiv:0907.2985}.


\bibitem
{Drinfeld}
V. Drinfeld, Degenerate affine Hecke algebras and Yangians, {\em Func. Anal. Appl.} {\bf 20} (1986), 56--58.

\bibitem
{JacLec}
N. Jacon and C. Lecouvey, On the Mullineux involution for Ariki-Koike algebras, {\em  J. Algebra} {\bf 321} (2009), 2156--2170.  


\bibitem
{Jbook} G.D. James, {\em  The representation theory of the symmetric groups}, Lecture Notes in Mathematics 682, Springer, Berlin, 1978.

\bibitem
{JamConj}
G.D. James, The decomposition matrices of ${\rm GL}\sb n(q)$ for $n\le 10$, 
{\em Proc. London Math. Soc. (3)} {\bf 60} (1990), 225--265. 

\bibitem
{JLM}
G.D. James, S. Lyle, A. Mathas, Rouquier blocks, {\em Math. Z.} {\bf 252} (2006),  511--531. 
 
\bibitem
{JM}
G. D. James and A. Mathas, 
Symmetric group blocks of small defect,  
{\em J. Algebra} {\bf 279} (2004), 566--612. 


\bibitem
{JaWi}
G. D. James and A. Williams, 
Decomposition numbers of symmetric groups by induction,  
{\em J. Algebra} {\bf 228} (2000), 119--142. 


\bibitem
{JMMO}
M. Jimbo, K. Misra, T. Miwa and M. Okado, 
Combinatorics of representations of $U_q(\widehat{\mathfrak{sl}}(n))$
at $q=0$, {\em Comm. Math. Phys.} {\bf 136} (1991), 543--566. 

\bibitem
{JW}
N. Jing and W. Wang, Twisted vertex representations and spin characters, {\em Math. Z.} {\bf 239} (2002), 715--746.

\bibitem
{Jo1} 
A.R. Jones, The structure of the Young symmetrizers for spin representations of the symmetric group. I, {\em J. Algebra} {\bf 205} (1998), 626--660.

\bibitem
{Jo2} 
A.R. Jones, The structure of the Young symmetrizers for spin representations of the symmetric group. II, {\em J. Algebra} {\bf 213} (1999), 381--404. 

\bibitem
{JN} 
A.R. Jones and M.L. Nazarov, Affine Sergeev algebra and $q$-analogues of the Young symmetrizers for projective representations of the symmetric group, {\em Proc. London Math. Soc. (3)} {\bf 78} (1999), 481--512.

\bibitem
{Jos}
T. Jozefiak, Relating spin representations of symmetric and hyperoctahedral groups, {\em J. Pure Appl. Algebra} {\bf 152} (2000), 187--193.

\bibitem
{Jucys1}
A. Jucys, On the Young operators of symmetric groups, {\em Litovsk. Fiz. Sb.} {\bf 6}  (1966), 163--180. (in Russian)

\bibitem
{Jucys2}
A. Jucys, Factorization of Young's projection operators for symmetric groups, {\em Litovsk. Fiz. Sb.} {\bf 11} (1971), 1--10. (in Russian)


\bibitem
{Jucys3}
A. Jucys, Symmetric polynomials and the center of the symmetric group ring, {\em Reports Math. Phys.} {\bf 5} (1974), 107--112.


\bibitem
{Kac} V. G. Kac, {\em Infinite Dimensional Lie Algebras}, Cambridge University Press, 1990. 

\bibitem
{Ka1}
M. Kashiwara,
Crystallizing the $q$-analogue of universal enveloping algebras,
{\em Commun. Math. Phys.} {\bf 133} (1990), 249--260.

\bibitem
{Ka2}
M. Kashiwara,
On crystal bases of the $q$-analogue of universal enveloping algebras,
{\em Duke Math. J.} {\bf 63} (1991), 465--516.

\bibitem
{KaG}
M. Kashiwara, 
Global crystal bases of quantum groups,
{\em Duke Math. J.} {\bf 69} (1993), 455--485.

\bibitem
{Kas}
M. Kashiwara,
On crystal bases, 
{\em CMS Conf. Proc.} {\bf 16} (1995), 155--197.

\bibitem
{KL}
D. Kazhdan and G. Lusztig,
Proof of the Deligne-Langlands conjecture for Hecke algebras,
{\em Invent. Math.} {\bf 87} (1987), 153--215.

\bibitem
{Ke} 
R. Kessar, Blocks and source algebras for the double covers of the symmetric and alternating groups, {\em J. Algebra} {\bf 186} (1996), 872--933.

\bibitem
{KS} 
R. Kessar and M. Schaps, Crossover Morita equivalences for blocks of the covering groups of the symmetric and alternating groups, {\em J. Group Theory} {\bf 9} (2006),  715--730.

\bibitem
{KL1}
M. Khovanov and A. Lauda,
A diagrammatic approach to categorification of quantum groups I, 
{\em Represent. Theory} {\bf 13} (2009), 309--347.

\bibitem
{KL2}
M. Khovanov and A. Lauda,
A diagrammatic approach to categorification of quantum groups II;
{\tt  arXiv:0804.2080}

\bibitem
{KL3}
M. Khovanov and A. Lauda,
A diagrammatic approach to categorification of quantum groups III;
{\tt  arXiv:0807.3250}.

\bibitem
{KlLi}
  P. B. Kleidman and M. W. Liebeck, `{\it The Subgroup Structure of the
Finite Classical Groups}', London Math. Soc. Lecture Note Ser. no.
$129$, Camb. Univ. Press, $1990$.

\bibitem
{KlWa}
 P. B. Kleidman and D.B. Wales, The projective characters of the symmetric groups that remain irreducible on subgroups, {\em J. Algebra} {\bf 138} (1991), 440--478. 


\bibitem
{KJS} 
A. Kleshchev, On restrictions of irreducible modular representations of semisimple algebraic groups and symmetric groups to some natural subgroups. I, {\em Proc. London Math. Soc. (3)} {\bf 69} (1994), 515--540.

\bibitem
{KBrI} 
A. Kleshchev, Branching rules for modular representations of symmetric groups. I, {\em J. Algebra} {\bf 178} (1995), 493--511.

\bibitem
{KBrII} 
A. Kleshchev, Branching rules for modular representations of symmetric groups 
II, {\em J. Reine Angew. Math.} {\bf 459} (1995), 163--212.

\bibitem
{KBrIII}
A. Kleshchev,
Branching rules for modular representations of symmetric groups
  {III}: some corollaries and a problem of {M}ullineux,
{\em J. London Math. Soc.} {\bf 54} (1996), 25--38.

\bibitem
{KCS}
A. Kleshchev, Completely splittable representations of symmetric groups, {\em J. Algebra} {\bf 181} (1996), 584--592. 

\bibitem
{KDec}
A. Kleshchev, On decomposition numbers and branching coefficients for symmetric and special linear groups, {\em Proc. London Math. Soc. (3)} {\bf 75} (1997), 497--558. 

\bibitem
{KBrIV}
A. Kleshchev, Branching rules for modular representations of symmetric groups. IV, {\em J. Algebra} {\bf 201} (1998), 547--572.

\bibitem
{Kbook}
A. Kleshchev, {\em Linear and Projective Representations of Symmetric Groups}, Cambridge University Press, Cambridge, 2005. 

\bibitem
{KN}
A. Kleshchev and D. Nash, An interpretation of the LLT algorithm, {\em preprint}, University of Oregon, 2009. 

\bibitem
{KR}
A. Kleshchev and A. Ram, Homogeneous representations of Khovanov-Lauda algebras, {\em J. European Math. Soc.}, to appear. 

\bibitem
{KS1}
 A. Kleshchev and J. Sheth, Representations of the symmetric group are reducible over simply transitive subgroups, {\em Math. Z.} {\bf 235} (2000), 99--109.


\bibitem
{KS2}
 A. Kleshchev and J. Sheth, Representations of the alternating group which are irreducible over subgroups, {\em Proc. London Math. Soc.} {\bf 84} (2002), 194--212.
 
\bibitem
{KT}
A. Kleshchev and Pham Huu Tiep, On restrictions of modular spin representations of symmetric and alternating groups, {\em Trans. Amer. Math. Soc.} {\bf 356} (2004), 1971--1999.

\bibitem
{KZ}
A. Kleshchev and A. S. Zalesski, Minimal polynomials of elements of order $p$ in $p$-modular projective representations of alternating groups, {\em Proc. Amer. Math. Soc.}  {\bf 132} (2004), 1605--1612.


\bibitem
{Ku}
J. Kujawa, Crystal structures arising from representations of ${\rm GL}(m\vert n)$, {\em Represent. Theory} {\bf 10} (2006), 49--85.

\bibitem
{KOR}
B. K\"ulshammer, J.B. Olsson, G.R. Robinson, Generalized blocks for symmetric groups,  {\em Invent. Math.} {\bf 151} (2003), 513--552.

\bibitem
{KuNe}
M.K\"unzer and G. Nebe, Elementary divisors of Gram matrices of certain Specht modules, {\em Comm. Algebra} {\bf 31} (2003), 3377--3427.

\bibitem
{LLT}
A. Lascoux, B. Leclerc and J.-Y. Thibon, Hecke algebras at roots of unity and crystal bases of quantum affine algebras, {\em Commun. Math. Phys.} {\bf 181} (1996), 205-263.

\bibitem
{Lec}
B. Leclerc,
Dual canonical bases,
quantum shuffles and $q$-characters,
{\em Math. Z.} {\bf 246} (2004), 691--732.

\bibitem
{LT1}
B. Leclerc and J.-Y. Thibon, Canonical bases of $q$-deformed Fock spaces, {\em Internat. Math. Res. Notices} {\bf 9} (1996), 447--456.

\bibitem
{LT2}
B. Leclerc and J.-Y. Thibon, Littlewood-Richardson coefficients and Kazhdan-Lusztig polynomials, in: 
{\em Combinatorial methods in representation theory (Kyoto, 1998)}, pp. 155--220, Adv. Stud. Pure Math., 28, Kinokuniya, Tokyo, 2000. 

\bibitem
{LTV}
B. Leclerc, J.-Y. Thibon, and E. Vasserot, Zelevinsky's involution at roots of unity, {\em J. Reine Angew. Math.} {\bf 513} (1999), 33--51. 


\bibitem
{Lu}
G. Lusztig, Affine Hecke algebras and their graded version, {\em J. Amer. Math. Soc.} {\bf 2} (1989), 599--635.


\bibitem
{Lubook}
G. Lusztig, {\em Introduction to Quantum Groups}, Birkh\"auser, 1993.

\bibitem
{Lyle}
S. Lyle, Some reducible Specht modules, {\em J. Algebra} {\bf 269} (2003), 536--543.

\bibitem
{LM}
S. Lyle and A. Mathas, Blocks of cyclotomic Hecke algebras, {\em Advances Math.} {\bf 216}  (2007), 854--878.


\bibitem
{Lynch}
T. E. Lynch,
{\em Generalized Whittaker vectors and representation theory},
PhD thesis, M.I.T., 1979.




\bibitem
{MM}
G. Malle and A. Mathas, Symmetric cyclotomic Hecke algebras, {\em J. Algebra}
{\bf 205} (1998), 275--293.

\bibitem
{Marcus}
A. Marcus, On equivalences between blocks of group algebras: reduction to the simple components, {\em J. Algebra} {\bf 184} (1996), 372--396. 

\bibitem
{MarcusAlt}
A. Marcus, Brou\'e's abelian defect group conjecture for alternating groups, {\em Proc. Amer. Math. Soc.} {\bf 132} (2004), 7--14.

\bibitem
{MarTanW3}
S. Martin and Kai Meng Tan, Defect 3 blocks of symmetric group algebras. I, {\em J. Algebra} {\bf 237} (2001),  95--120.

\bibitem
{MarTan}
S. Martin and Kai Meng Tan, $[3:2]$-pairs of symmetric group algebras and their intermediate defect 4 blocks, {\em  J. Algebra} {\bf 288} (2005), 505--526. 



\bibitem
{MathasB}
A. Mathas, {\em Iwahori-Hecke algebras and Schur algebras of the symmetric group}, University Lecture Series 15,  American Mathematical Society, Providence, RI, 1999.

\bibitem
{Ma} A. Mathas, Seminormal forms and Gram determinants for cellular algebras, {\em J. reine angew. Math.} {\bf 619} (2008), 141--173.

\bibitem
{MiM}
K.~C. Misra and T. Miwa,
Crystal base of the basic representation of $U_q(\widehat{\mathfrak{sl}}_n)$,
{\em Comm. Math. Phys.} {\bf 134} (1990), 79--88.



\bibitem
{Muller}
J. M\"uller, Brauer trees for the Schur cover of the symmetric group, {\em J. Algebra} {\bf  266} (2003), 427--445.

\bibitem
{Mul}
G. Mullineux, Bijections of $p$-regular partitions and $p$-modular irreducibles of the symmetric groups, {\em J. London Math. Soc. (2)} {\bf 20} (1979), 60--66.

\bibitem
{MurphyYoung}
G. E. Murphy, 
A new construction of Young's seminormal representation of the symmetric group, 
{\em J. Algebra} {\bf 69} (1981), 287--291. 

\bibitem
{MurphyId}
G. E. Murphy, 
The idempotents of the symmetric group and Nakayama's conjecture, 
{\em J. Algebra} {\bf 81} (1983), 258--265. 

\bibitem
{MurphyCell}
G. E. Murphy, 
The representations of Hecke algebras of type $A_n$,  
{\em J. Algebra} {\bf 173} (1995), 97--121. 




\bibitem
{NO}
C. N\v ast\v asescu and F. Van Oystaeyen,
{\em Methods of Graded Rings}, Lecture Notes in Math. 1836,
Springer, 2004.

\bibitem
{Naz1}
M.L. Nazarov, Young's orthogonal form of irreducible projective representations of the symmetric group, {\em J. London Math. Soc. (2)} {\bf 42} (1990), 437--451. 

\bibitem
{Naz2}
M.L. Nazarov, Young's symmetrizers for projective representations of the symmetric group. Adv. Math. 127 (1997), no. 2, 190--257.

\bibitem
{OV}
A. Okounkov and A. Vershik, A new approach to representation theory of symmetric groups. {\em Selecta Math. (N.S.)} {\bf 2} (1996), 581--605. 

\bibitem
{Paget1}
R. Paget, Induction and decomposition numbers for RoCK blocks, {\em Q. J. Math.} {\bf 56} (2005), 251--262. 

\bibitem
{Paget}
R. Paget, The Mullineux map for RoCK blocks, {\em Comm. Algebra} {\bf 34} (2006),  3245--3253. 

\bibitem
{Ph}
A.M. Phillips, Restricting modular spin representations of symmetric and alternating groups to Young-type subgroups, {\em Proc. London Math. Soc. (3)} {\bf 89} (2004),  623--654. 

\bibitem
{P}
A. Premet,
Special transverse slices and their enveloping algebras,
{\em Adv. Math.}
{\bf 170} (2002), 1--55.

\bibitem
{PuigN}
L. Puig, The Nakayama conjecture and the Brauer pairs, in {\em  S\'eminaire sur les groupes finis, Tome III, ii,} pp. 171--189, 
Publ. Math. Univ. Paris VII, 25, Univ. Paris VII, Paris, 1986. 

\bibitem
{PuigScopes} 
L. Puig, On Joanna Scopes' criterion of equivalence for blocks of symmetric groups, {\em Algebra Colloq.} {\bf 1} (1994), 25--55. 


\bibitem
{Ram}
A. Ram, Seminormal representations of Weyl groups and Iwahori-Hecke algebras. 
{\em Proc. London Math. Soc.} {\bf 75} (1997), 99--133. 

\bibitem
{Ri1}
J. Rickard, Morita theory for derived categories, {\em J. Lond. Math. Soc.} {\bf 39} (1989), 436--456.

 
\bibitem
{Ri2}
J. Rickard, Derived equivalences as derived functors, {\em J. Lond. Math. Soc.} {\bf 43} (1991), 37--48.


\bibitem
{Ric}
J. Rickard,
Equivalences of derived categories for symmetric algebras,
{\em J. Algebra} {\bf 257} (2002), 460--481.


\bibitem
{R2}
R. Rouquier, Isom\'etries parfaites dans les blocs \'a d\'efaut ab\'elien des groupes sym\'etriques et sporadiques, {\em J. Algebra} {\bf 168} (1994), 648--694.

\bibitem
{R1}
R. Rouquier,
Derived equivalences and finite dimensional algebras,
{\em Proc. ICM (Madrid 2006)} {\bf 2}, 191--221, EMS Publishing House, 2006.



\bibitem
{Ro}
R. Rouquier, $2$-Kac-Moody algebras, {\tt arXiv:0812.5023}.


\bibitem
{Saxl}
J. Saxl, The complex characters of the symmetric groups that remain irreducible in subgroups, {\em J. Algebra} {\bf 111} (1987), 210--219.



\bibitem
{Sh}
V. Shchigolev, Generalization of modular lowering operators for ${\rm GL}\sb n$, {\em Comm. Algebra} {\bf 36} (2008), 1250--1288. 

\bibitem
{Shch}
V. Shchigolev, On extensions and branching rules of modules that are close to completely splittable, {\em Sb. Math.} {\bf 196} (2005), no. 7-8, 1209--1249.  

\bibitem
{Schur}
I. Schur, 
{\em Gesammelte Abhandlungen}, Band I. 
Herausgegeben von Alfred Brauer und Hans Rohrbach. Springer-Verlag, Berlin-New York, 1973. 

\bibitem
{Scopes}
J. Scopes, Cartan matrices and Morita equivalence for blocks of the symmetric groups, {\em J. Algebra} {\bf 142} (1991), 441--455.

\bibitem
{Sc}
L.L. Scott, Representations in characteristic $p$, {\em The Santa Cruz Conference on Finite Groups (Univ. California, Santa Cruz, Calif., 1979)}, pp. 319--331, 
Proc. Sympos. Pure Math., 37, Amer. Math. Soc., Providence, R.I., 1980. 

\bibitem
{Se85}
A.N. Sergeev, Tensor algebra of the identity representation as a module over the Lie 
superalgebras $GL(n, m)$ and $Q(n)$, {\em Math. USSR Sbornik} {\bf 51} (1985), 419Ð427

\bibitem
{Se99}
A.N. Sergeev, The Howe duality and the projective representations of symmetric groups,  {\em Represent. Theory} {\bf 3} (1999), 416--434. 


\bibitem
{ShuWang}
B. Shu and W. Wang, Modular representations of the ortho-symplectic supergroups, {\em  Proc. Lond. Math. Soc. (3)} {\bf 96} (2008), 251--271. 

\bibitem
{Stern}
E. Stern,
Semi-infinite wedges and vertex operators,
{\em Internat. Math. Res. Notices} {\bf 1995}, 201--220.

\bibitem
{Soergel}
W. Soergel,
Kategorie $\mathcal O$, perverse Garben und Moduln \"uber den Koinvarianten zur Weylgruppe, {\em J. Amer. Math. Soc.} {\bf 3} (1990), 421--445.







\bibitem
{TanM}
Kai Meng Tan,  Martin's conjecture holds for weight 3 blocks of symmetric groups, {\em J. Algebra} {\bf 320} (2008), 1115--1132.

\bibitem
{Ts}
S. Tsuchioka,  A modular branching rule for the generalized symmetric groups, {\em J. Algebra} {\bf 316} (2007), 459--470.


\bibitem
{T}
W. Turner, 
Rock blocks, to appear in {\em Mem. Amer. Math. Soc.};
{\tt arXiv:0710.5462}.


\bibitem
{Uglov}
D. Uglov, 
Canonical bases of higher level $q$-deformed Fock spaces and Kazhdan-Lusztig 
polynomials,
in: {\em Physical Combinatorics (Kyoto, 1999)}, pp. 249--299,
Progress in Math., 191, 
Birkh\"auser, 2000.

\bibitem
{VV1}
M. Varagnolo and E. Vasserot,
On the decomposition matrices of the quantized Schur algebra,
{\em Duke Math. J.} {\bf 100} (1999), 267--297.


\bibitem
{VV3}
M. Varagnolo and E. Vasserot,
Canonical bases and Khovanov-Lauda algebras,
{\tt arXiv:0901.3992}.



\bibitem
{Vaz} M. Vazirani, 
Parameterizing Hecke algebra modules: Bernstein-Zelevinsky multisegments, Kleshchev multipartitions, and crystal graphs, {\em Transform. Groups} {\bf 7} (2002), 267--303.





\bibitem
{WW}
J. Wan and W. Wang, 
Modular representations and branching rules for wreath Hecke algebras, 
{\em Int. Math. Res. Not. IMRN} (2008), Art. ID rnn128, 31 pp. 

\bibitem
{Wang}
W. Wang, Spin Hecke algebras of finite and affine types, {\em Adv. Math.} {\bf 212}  (2007), 723--748. 

\bibitem
{W}
H. Wenzl, Hecke algebras of Type $A_n$ and subfactors, {\em Invent. Math.} {\bf 92} (1988), 349--383. 

\bibitem
{Wi}
A. Williams, 
Symmetric group decomposition numbers for some three-part partitions, {\em Comm. Algebra} {\bf 34} (2006), 1599--1613. 

\bibitem
{Xu}
M. Xu, On $p$-series and the Mullineux conjecture, {\em Comm. Algebra} {\bf 27} (1999), 5255--5265. 

\bibitem
{Ya1}
M. Yamaguchi, A duality of a twisted group algebra of the hyperoctahedral group and the queer Lie superalgebra, in {\em  Combinatorial methods in representation theory (Kyoto, 1998)}, 401--422, Adv. Stud. Pure Math., 28, Kinokuniya, Tokyo, 2000. 

\bibitem
{Ya2}
M. Yamaguchi,  A duality of the twisted group algebra of the symmetric group and a Lie superalgebra, {\em J. Algebra} {\bf 222} (1999), 301--327.

\bibitem
{Y2}
X. Yvonne,
Canonical bases of higher-level $q$-deformed Fock spaces, 
{\em J. Algebraic Combin.} {\bf 26}  (2007), 383--414.




\end{thebibliography}
\end{document}